\ifx\shlhetal\undefinedcontrolsequence\let\shlhetal\relax\fi

\documentstyle[11pt]{article}

\input amssym.def
\input amssym.tex

\newcommand{\inv}{{\rm inv}}
\newcommand{\Inv}{{\rm Inv}}
\newcommand{\cfinv}{{\rm cf{-}inv}}
\newcommand{\cfInv}{{\rm cf{-}Inv}}
\newcommand{\cf}{{\rm cf}}
\newcommand{\Depth}{{\rm Depth}}
\newcommand{\Ult}{{\rm Ult}\;}
\newcommand{\rest}{{\restriction}}
\newcommand{\ve}{\varepsilon}
\newcommand{\dom}{{\rm dom}}
\newcommand{\base}{{\rm base}}
\newcommand{\otp}{{\rm otp}}
\newcommand{\rng}{{\rm rng}}
\newcommand{\ind}{{\rm ind}}

\newcommand{\hcof}{{\rm h{-}cof}}
\newcommand{\mod}{{\rm mod}\;}
\newcommand{\cast}{\mathop{\circledast}}
\newcommand{\op}{\mathop{{\bf O}}}
\newcommand{\hd}{{\rm hd}}
\newcommand{\hL}{{\rm hL}}
\newcommand{\ut}{{\rm ut}}
\newcommand{\QED}{\hfill\hspace{0.2in}\vrule width 6pt height 6pt depth 0pt
\vspace{0.1in}}

\newtheorem{theorem}{Theorem}[section]
\newtheorem{claim}{Claim}[theorem]
\newtheorem{fact}[theorem]{Fact}
\newtheorem{proposition}[theorem]{Proposition}
\newtheorem{definition}[theorem]{Definition}
\newtheorem{example}[theorem]{Example}
\newtheorem{corollary}[theorem]{Corollary}
\newtheorem{lemma}[theorem]{Lemma}
\newtheorem{conclusion}[theorem]{Conclusion}
\newtheorem{problem}[theorem]{Problem}

\def\lessdot{\mathrel{\mathord{<}\!\!\raise
0.8 pt\hbox{$\scriptstyle\circ$}}}

\def\lesseqdot{\mathrel{\mathord{\leq}\!\!\raise
0.8 pt\hbox{$\scriptstyle\circ$}}}

\def\mid{\|}

\author{
{\bf Andrzej Ros\l anowski}\\
Institute of Mathematics\\
Hebrew University of Jerusalem\\
91904 Jerusalem, Israel\\
and\\
Mathematical Institute\\
Wroc\l aw University\\
50 384 Wroc\l aw, Poland
\and
{\bf Saharon Shelah}\thanks{\ \ The research was partially supported by DFG
grant Ko 490/7-1. This is publication number 534 of the second author}\\
Institute of Mathematics\\
Hebrew University of Jerusalem\\
91904 Jerusalem, Israel\\
and\\
Department of Mathematics\\
Rutgers University\\
New Brunswick, NJ 08903, USA
}
\title{Cardinal invariants of ultraproducts of Boolean algebras}
\date{\today}

\bigskip 

\begin{document}
\setcounter{page}{0}

\maketitle

\begin{abstract}
\noindent We deal with some of problems posed by Monk \cite{M}, \cite{M3}
and related to cardinal invariant of ultraproducts of Boolean algebras. We
also introduce and investigate several new cardinal invariants.
\end{abstract}
\bigskip

\noindent{\em AMS 1991 Mathematics Subject Classification}\\
primary:\ 03G05, 03E05,\quad secondary: 06E15, 03E35
\bigskip

\setcounter{section}{-1}

\section{Introduction}
In the present paper we deal with cardinal invariants of Boolean
algebras and ultraproducts. Several questions in this area were  posed by
Monk (\cite{M}, \cite{M2}, \cite{M3}) and we address some of them.
General schema of these problems can be presented in the following
fashion. Let $\inv$ be a cardinal function on Boolean algebras. Suppose
that $B_i$ are Boolean algebras (for $i<\kappa$) and that $D$ is an
ultrafilter on the cardinal $\kappa$. We ask what is the relation
between $\inv(\prod\limits_{i<\kappa}B_i/D)$ and
$\prod\limits_{i<\kappa}\inv(B_i)/D$? For each invariant $\inv$ we may
consider two questions:
\begin{description}
\item[$(<)_{\inv}$]\hskip2cm is $\inv(\prod\limits_{i<\kappa}B_i/D)<
\prod\limits_{i<\kappa}\inv(B_i)/D$ possible?
\item[$(>)_{\inv}$]\hskip2cm is $\inv(\prod\limits_{i<\kappa}B_i/D)>
\prod\limits_{i<\kappa}\inv(B_i)/D$ possible?
\end{description}
We deal with these questions for several cardinal invariants. We find it
helpful to introduce ``finite'' versions $\inv_n$ of the invariants.
This helps us in some problems as $\inv^+(\prod\limits_{i<\kappa} B_i/D)
\geq\prod\limits_{i<\kappa}\inv^+_{f(i)}(B_i)/D$ for each function
$f:\kappa\longrightarrow\omega$ such that $\lim\limits_{D} f=\omega$.

In section 1 we will give a general setting of the subject. These
results were known much earlier (at least to the second author). We
present them here to establish a uniform approach to the invariants and
show how the \L o\'s theorem applies. In the last part of this section
we present a simple method which uses the main result
of \cite{MgSh 433} to show the consistency of the inequality
$\inv(\prod\limits_{i<\kappa}B_i/D)<\prod\limits_{i<\kappa}\inv(B_i)/D$ for
several invariants $\inv$. These problems will be fully studied and presented
in \cite{MgSh 433}.

Section 2 is devoted to the (topological) density of Boolean algebras.
We show here that, in ZFC, the answer to the question $(<)_d$ is
``yes''. This improves Theorem A of \cite{KoSh 415} (a
consistency result) and answers (negatively) Problem J of \cite{M3}. It
should be remarked here that the answer to $(>)_d$ is ``no'' (see
\cite{M2}).

In the third section we introduce strong $\lambda$-systems which are one of
tools for our constructions. Then we apply them to build
Boolean algebras which (under some set-theoretical assumptions) show that
the inequalities $(>)_{\hcof}$ and $(>)_{\rm inc}$ are possible
(a consistency). These results seem to be new, the second one can be
considered as a partial answer to Problem X of \cite{M3}. We get similar
constructions for spread, hereditary Lindel\"of degree and hereditary
density. However they are not sufficient to give in ZFC positive answers to
the corresponding questions $(>)_{\inv}$. These investigations are continued
in \cite{Sh 620}, where the respective Boolean algebras are built in ZFC. The
consistencies of the reverse inequalities will be presented in \cite{MgSh 433}.

The fourth section deals with the independence number and the tightness.  It
has been known that both questions $(>)_{\ind}$ and $(>)_t$ have the answer
``yes''. In coming paper \cite{MgSh 433} it will be shown that $(<)_{\ind}$,
$(<)_t$ may be answered positively (a consistency result; see section 1
too). Our results here were inspired by other sections of this paper and
\cite{Sh 503}. We introduce and study ``finite'' versions of the
independence number getting a surprising asymmetry between odd and even
cases. A completely new cardinal invariant appears naturally here. It has
some reflection in what we can show for the tightness. Finally we
re-present and re-formulate the main result of \cite{Sh 503} (on
products of interval Boolean algebras) 
putting it in our general schema and showing explicitly its heart. 
\bigskip

\noindent{\bf History:} A regular study of cardinal invariants of
Boolean algebras was initialized in \cite{M}, where several problems were
posed. Those problems stimulated and directed the work in the area. Some
of the problems were naturally related to the behaviour of the
invariants in ultraproducts and that found a reflection in papers coming
later. Several bounds, constructions and consistency results were proved
in \cite{Pe}, \cite{Sh 345}, \cite{KoSh 415}, \cite{MgSh 433}, \cite{Sh
479}, \cite{Sh 503}. New techniques of constructions of Boolean algebras
were developed in \cite{Sh 462} (though the relevance of the methods for
ultraproducts was not stated explicitly there).

This paper is, in a sense, a development of the notes ``F99: Notes on
cardinal invariants\ldots'' which the second author wrote in January 1993.
A part of these notes is incorporated here, other results will be presented in
\cite{MgSh 433} and \cite{RoSh 599}. 

The methods and tools for building Boolean algebras which we present here
will be applied in a coming paper to deal with the problems of attainment in
different representations of cardinal invariants.
\bigskip

\noindent{\bf Notation:} Our notation is rather standard. All cardinals
are assumed to be infinite and usually they are denoted by $\lambda$,
$\kappa$, $\theta$, $\Theta$ (with possible indexes).

We say that a family $\{\langle s^\alpha_0,\ldots,s^\alpha_{m-1}\rangle:
\alpha<\lambda\}$ of finite sequences forms a $\Delta$-system with the
root $\{0,\ldots,m^*-1\}$ (for some $m^*\leq m$) if the sets
$\{s^\alpha_{m^*},\ldots,s^\alpha_{m-1}\}$ (for $\alpha<\lambda$) are
pairwise disjoint and 
\[(\forall\alpha<\lambda)(\forall l<m^*)(s^\alpha_l=s^0_l).\] 

In Boolean algebras we use $\vee$ (and $\bigvee$), $\wedge$ (and
$\bigwedge$) and $-$ for the Boolean operations. If $B$ is a Boolean
algebra, $x\in B$ then $x^0=x$, $x^1=-x$.

The sign $\cast$ stands for the operation of the free product of
Boolean algebras (see \cite{HoBA}, def.11.1) and  $\prod^w$ denotes the
weak product of Boolean algebras (as defined in \cite{HoBA}, p.112).

All Boolean algebras we consider are assumed to be infinite (and we will
not repeat this assumption). Similarly whenever we consider a cardinal
invariant $\inv(B)$ we assume that it is infinite.
\bigskip

\noindent{\bf Acknowledgment:} We would like to thank Professor Donald Monk
for his very helpful comments at various stages of preparation of the paper as
well as for many corrections and improvements.

\section{Invariants and ultraproducts}
\subsection{Definable cardinal invariants.}
In this section we try to systematize the definition of cardinal invariants
and we define what is a def.car.invariant (definable cardinal invariant) of
Boolean algebras. Then we get immediate consequences of this approach for
ultraproducts. Actually, Boolean algebras are irrelevant in this section and
can be replaced by any structures.

\begin{definition}
\label{basdef}
\begin{enumerate}
\item For a (not necessary first order)
theory $T$ in the language of Boolean algebras plus one distinguished
predicate $P=P_0$ (unary if not said otherwise) plus, possibly, some others
$P_1,P_2,\ldots$ we define cardinal invariants $\inv_T$, $\inv^+_T$ of Boolean
algebras by (for a Boolean algebra $B$):
\begin{quotation}
$\inv_T (B) \stackrel{\rm def}{=}\sup \{\mid P \mid : (B,P_n)_n \mbox{ is a
model of } T \}$

$\inv^+_T (B) \stackrel{\rm def}{=}\sup \{\mid P \mid^+ : (B,P_n)_n\mbox{ is
a model of } T \}$

$\Inv_T(B)\stackrel{\rm def}{=}\{\|P\|: (B,P_n)_n\mbox{ is a model of } T\}$
\end{quotation}
We call $\inv^{(+)}_T$ a {\em def.car.~invariant} (definable cardinal
invariant).

\item If in 1., $T$ is first order, we call such cardinal invariant a
{\em def.f.o.car.~invariant} (definable first order cardinal invariant).

\item A theory $T$ is {\em $n$-universal in $(P_0,P_1)$} if all sentences
$\phi\in T$ are of the form
\[(\forall x_1,\ldots,x_n\in P_0)\psi(\bar{x}),\]
where all occurrences of $x_1,\dots,x_n$ in $\psi$ are free and $P_0$ does not
appear there and any appearance of $P_1$ in $\psi$ is in the form
$P_1(x_{i_0},\ldots,x_{i_k})$ with no more complicated terms.

If we allow all $n$ then $T$ is said to be {\em universal in $(P_0,P_1)$}.

Note: quantifiers can still occur in $\psi(\bar{x})$ on other variables.

\item If in 1., $T$ is universal in $(P_0,P_1)$, first order except
the demand that $P_1$ is a well ordering of $P_0$ we call such cardinal
invariant {\em def.u.w.o.car.~invariant} (definable universal well ordered
cardinal invariant).

\item If in 1., $P_1$ is a linear order on $P$ (i.e. $T$ says so)
and in the definition of $\inv_T(B)$, $\inv^+ (B)$ we replace ``$\mid P
\mid$'' by the cofinality of $(P,P_1)$ then we call those cardinal
invariants {\em def.cof. invariant} (definable cofinality invariant,
$\cfinv_T$); we can have the f.o. and the u.w.o. versions. We define similarly
$\cfInv_T (B)$ as the set of such cofinalities. To use $\cfinv$ we can put it
in ${}^+$ (we may omit ``cf-'' if the context allows it). We can use the order
type instead of the cofinality and cardinality writing {\rm ot}-\inv. For the
cardinality we may use {\em car}-\inv.

\item For a theory $T$ as in 2.,  the {\em minimal definable first order
cardinal invariant} of $B$ (determined by $T$) is $\min~\Inv_T(B)$.
\end{enumerate}
\end{definition}

To avoid a long sequence of definitions we refer the reader to \cite{M},
\cite{M2} for definitions of the cardinal functions below. Those invariants
which are studied in this paper are defined in the respective sections.

\begin{proposition}
\begin{enumerate}
\item The following cardinal invariants of Boolean algebras are
def.f.o.car.~invariants (of course each has two versions: $\inv$ and
$\inv^+$):
\begin{quotation}
\noindent $c$ (cellularity), Length, irr (irredundance), cardinality, ind
(independence), $s$ (spread), Inc (incomparability).
\end{quotation}

\item The following cardinal invariants of Boolean algebras are def.f.o.cof.
invariants: 
\begin{quotation}
\noindent hL (hereditary Lindelof), hd (hereditary density).
\end{quotation}

\item The following cardinal invariants of Boolean algebras are def.u.w.o.car.
invariants 
\begin{quotation}
\noindent Depth, $t$ (tightness), h-cof (hereditary cofinality), hL, hd.
\end{quotation}

\item $\pi$ (algebraic density) and $d$ (topological density) are minimal 
def.f.o.card. invariants.
\end{enumerate}
\end{proposition}

\noindent PROOF: All unclear cases are presented in next sections.\QED

\begin{proposition}
\label{0.4}
\begin{enumerate}
\item If $\inv^+_T (B)$ is a limit cardinal then the $\sup$ in the
definition of $\inv _T (B)$ is not obtained and $\inv_T (B)=\inv^+_T (B)$. 

\item  If $\inv^+_T (B)$ is not a limit cardinal then it is $(\inv_T (B))^+$
and the $\sup$ in the definition of $\inv_T (B)$ is obtained.\QED
\end{enumerate}
\end{proposition}

\begin{definition}
A linear order $(I,<)$ is $\Theta$-like if
\[\|I\|=\Theta\ \ \mbox{ and }\ \ (\forall a\in I)(\|\{b\in I:
b<a\}\|<\Theta).\] 
\end{definition}

\begin{proposition}
\label{fo}
Assume that $\inv^{(+)}_T$ is a definable first order
cardinal invariant. Assume further that: $D$ is an ultrafilter on a cardinal
$\kappa$, for $i < \kappa$, $B_i$ is a Boolean algebra and $B \stackrel{\rm
def}{=}\prod\limits_{i < \kappa } B_i/D$. Then
\begin{description}
\item[{(a)}] if $\lambda_i<\inv^+_T(B_i)$ for $i < \kappa$ then
$\prod\limits_{i<\kappa } \lambda_i/D<\inv^+_T (B)$, 
\item[{(b)}] $\prod\limits_{i<\kappa}\inv^+_T(B_i)/D\leq\inv^+_T (B)$,
\item[{(c)}] if $\inv_T(B)<\prod\limits_{i<\kappa }\inv_T(B_i)/D$ then for the
$D$-majority of $i<\kappa$ we have:\\
$\lambda_i\stackrel{\rm def}{=}\inv_T(B_i)$ is a limit cardinal and the linear
order $\prod\limits_{i<\kappa}(\lambda_i,<)/D$ is $(\inv_T(B))^+$--like; hence
for the $D$-majority of the $i<\kappa$ we have: $\lambda_i$ is a regular limit
cardinal (i.e.~weakly inaccessible),
\item[{(d)}] $\min\Inv_T(B)\leq\prod\limits_{i<\kappa}\min\Inv_T(B_i)/D$.
\end{description}
\end{proposition}

\noindent PROOF: (a)\ \ \ This is an immediate consequence of \L o\'s theorem.
\medskip

\noindent (b)\ \ \  For $i<\kappa$ let $\lambda_i=\inv^+_T(B_i)$. Suppose 
$\lambda<\prod\limits_{i<\kappa}\lambda_i/D$. As $\prod\limits_{i<\kappa}
(\lambda_i,<)/D$ is a linear order of cardinality $>\lambda$ we find $f\in
\prod\limits_{i<\kappa}\lambda_i$ such that 
\[\|\{g/D\in\prod_{i<\kappa}\lambda_i/D: g/D < f/D\}\|\geq\lambda.\]
Since $f(i)<\inv^+_T(B_i)$ (for $i<\kappa$) we may apply (a) to conclude that
\[\lambda\leq\|\prod_{i<\kappa}f(i)/D\|<\inv^+_T(B).\]

\noindent (c)\ \ \ Let $\lambda=\inv_T(B)$, $\lambda_i=\inv_T(B_i)$ and assume
that $\lambda<\prod\limits_{i<\kappa}\lambda_i/D$. By part (b) we conclude
that then  
\begin{description}
\item[{$(*)$}]\ \ \ \ \ $\lambda^+=\prod\limits_{i<\kappa}\inv^+_T(B_i)/D
=\prod\limits_{i<\kappa}\inv_T(B_i)/D= \inv^+_T(B)$.
\end{description}
Let $A=\{i<\kappa:\inv_T(B_i)<\inv^+_T(B_i)\}$. Note that $A\notin D$: if
not then we may assume $A=\kappa$ and for each $i<\kappa$ we have $\lambda_i
<\inv^+_T(B_i)$. By part (a) and (*) above we get $\lambda^+=\prod\limits_{i<
\kappa}\lambda_i/D<\inv^+_T(B)$, a contradiction. Consequently we may assume
that $A=\emptyset$. Thus for each $i<\kappa$ we have $\lambda_i=\inv_T(B_i)=
\inv^+_T(B_i)$ and $\lambda_i$ is a limit cardinal, $\lambda_i=\sup\Inv_T(B_i)
\notin\Inv_T(B_i)$ (by \ref{0.4}).

The linear order $\prod\limits_{i<\kappa}(\lambda_i,<)/D$ is of the cardinality
$\lambda^+$ (by $(*)$). Suppose that $f\in\prod\limits_{i<\kappa}\lambda_i$ and
choose $\mu_i\in\Inv_T(B_i)$ such that $f(i)\leq\mu_i$ for $i<\kappa$. Then 
$\|\prod\limits_{i<\kappa}f(i)/D\|\leq\prod\limits_{i<\kappa}\mu_i/D\in\Inv_T(
\prod\limits_{i<\kappa}B_i/D)\subseteq \lambda^+$. Hence the order
$\prod\limits_{i<\kappa}(\lambda_i,<)/D$ is $\lambda^+$-like. 

Finally assume that $A=\{i<\kappa: \lambda_i\mbox{ is singular }\}\in D$, so
w.l.o.g.~$A=\kappa$. Choose cofinal subsets $Q_i$ of $\lambda_i$ such that 
$Q_i\subseteq\lambda_i=\sup Q_i$, $\|Q_i\|=\cf(\lambda_i)$ (for $i<\kappa$)
and let $M_i=(\lambda_i,<,Q_i,\ldots)$. Take the ultrapower
$M=\prod\limits_{i<\kappa}M_i/D$ and note that $M\models$``$Q^M$ is unbounded
in $<^M$''. As earlier, $\|Q^M\|=\prod\limits_{i<\kappa}\|Q_i\|/D\leq\lambda$
so $\cf(\prod\limits_{i<\kappa}(\lambda_i,<)/D)\leq\lambda$ what contradicts
$\lambda^+$-likeness of the product order.
\medskip

\noindent (d)\ \ \ It follows from (a).\QED

\begin{definition}
Let $(I,<)$ be a partial order.
\begin{enumerate}
\item The depth $\Depth(I)$ of the order $I$ is the supremum of
cardinalities of well ordered (by $<$) subsets of $I$.
\item $I$ is $\Theta$-Depth-like if $I$ is a linear ordering which contains
a well ordered cofinal subset of length $\Theta$ but 
$\Depth^+(\{b\in I: b<a\},<)\leq\Theta$ for each $a\in I$.
\end{enumerate}
\end{definition}

\begin{lemma}
\label{depth}
Let $D$ be an ultrafilter on a cardinal $\kappa$, $\lambda_i$ (for $i<\kappa$)
be cardinals. Then:
\begin{enumerate}
\item if there is a $<_D$-increasing sequence $\langle f_\alpha/D: \alpha
\leq\mu_0\rangle\subseteq\prod\limits_{i<\kappa}(\lambda_i^+,<)/D$, $\mu_0$ is
a cardinal then $\mu_0<\Depth^+(\prod\limits_{i<\kappa}(\lambda_i,<)/D)$,
\item $\Depth(\prod\limits_{i<\kappa}(\lambda_i^+,<)/D)\leq
\Depth^+(\prod\limits_{i<\kappa}(\lambda_i,<)/D)$ and hence

$\Depth^+(\prod\limits_{i<\kappa}(\lambda_i^+,<)/D)\leq
(\Depth^+(\prod\limits_{i<\kappa}(\lambda_i,<)/D))^+$.
\end{enumerate}
\end{lemma}

\noindent PROOF: 1)\ \ \ Let $\mu_1=\cf(\prod\limits_{i<\kappa}(\lambda_i,<)
/D)$, so $\mu_1<\Depth^+(\prod\limits_{i<\kappa}(\lambda_i,<)/D)$. If
$\mu_0\leq\mu_1$ then we are done, so let us assume that $\mu_0>\mu_1$ and let
us consider two cases.
\medskip

\noindent{\sc Case A:}\ \ \ $\cf(\mu_0)\neq\mu_1$.
\medskip

\noindent Let $\langle g_\beta/D:\beta<\mu_1\rangle$ be an increasing sequence
cofinal in $\prod\limits_{i<\kappa}(\lambda_i,<)/D$. For each $i<\kappa$
choose an increasing sequence $\langle A^i_\xi: \xi<\lambda_i\rangle$ of
subsets of $f_{\mu_0}(i)$ such that $f_{\mu_0}(i)=\bigcup\limits_{\xi<
\lambda_i} A^i_\xi$ and $\|A^i_\xi\|<\lambda_i$. Then
\[(\forall\alpha<\mu_0)(\exists\beta<\mu_1)(\{i<\kappa: f_\alpha(i)\in 
A^i_{g_\beta(i)}\}\in D)\]
and, passing to a subsequence of $\langle f_\alpha/D:\alpha<\mu_0\rangle$ if
necessary, we may assume that for some $\beta_0<\mu_1$ for all $\alpha<\mu_0$ 
\[\{i<\kappa: f_\alpha(i)\in A^i_{g_{\beta_0}(i)}\}\in D\]
(this is the place we use the additional assumption $\cf(\mu_0)\neq\mu_1)$.
Each set $A^i_{g_{\beta_0}(i)}$ is order--isomorphic to some ordinal $g(i)<
\lambda_i$ (as $\|A^i_{g_{\beta_0}(i)}\|<\lambda_i$). These isomorphisms
give us a ``copy'' of the sequence $\langle f_\alpha/D:\alpha<\mu_0\rangle$
below some $g/D\in\prod\limits_{i<\kappa}\lambda_i/D$, witnessing $\mu_0<
\Depth^+ (\prod\limits_{i<\kappa}(\lambda_i,<)/D)$.
\medskip

\noindent{\sc Case B:}\ \ \ $\cf(\mu_0)=\mu_1<\mu_0$.
\medskip

\noindent For each regular cardinal $\mu\in (\cf(\mu_0),\mu_0)$ we may apply
Case A to $\mu$ and the sequence $\langle f_\alpha/D:\alpha\leq\mu\rangle$ and
conclude $\mu<\Depth^+(\prod\limits_{i<\kappa}(\lambda_i,<)/D)$. Hence $\mu_0
\leq\Depth^+(\prod\limits_{i<\kappa}(\lambda_i,<)/D)$. Let $\langle\mu^\xi:\xi
<\cf(\mu_0)\rangle\subseteq(\cf(\mu_0),\mu_0)$ be an increasing cofinal in
$\mu_0$ sequence of regular cardinals. Note that for each $\xi<\cf(\mu_0)$ and
a function $f\in\prod\limits_{i<\kappa}\lambda_i$ we can find a
$<_D$-increasing sequence $\langle h^*_\alpha:\alpha<\mu^\xi\rangle\subseteq
\prod\limits_{i<\kappa}\lambda_i$ such that $f<_D h^*_0$. Using this fact we
construct inductively a $<_D$--increasing sequence $\langle h_\alpha/D:\alpha
<\mu_0\rangle\subseteq\prod\limits_{i<\kappa}\lambda_i/D$ (which will show
that $\mu_0<\Depth^+(\prod\limits_{i<\kappa}(\lambda_i,<)/D)$):

\noindent Suppose we have defined $h_\alpha$ for $\alpha<\mu^\xi$ (for some 
$\xi<\cf(\mu_0)$). Since $\mu^\xi$ is regular and $\mu^\xi\neq\mu_1$ the
sequence $\langle h_\alpha/D:\alpha<\mu^\xi\rangle$ cannot be cofinal in
$\prod\limits_{i<\kappa}(\lambda_i,<)/D$. Take $f/D\in\prod\limits_{i<\kappa}
\lambda_i/D$ which $<_D$-bounds the sequence. By the previous remark we find a
$<_D$-increasing sequence $\langle h_\alpha/D:\mu^\xi\leq\alpha<\mu^{\xi+1}
\rangle\subseteq\prod\limits_{i<\kappa}\lambda_i/D$  such that $f<_D
h_{\mu^\xi}$. So the sequence $\langle h_\alpha: \alpha<\mu^{\xi+1}\rangle$ is
increasing. 

\noindent Now suppose that we have defined $h_\alpha/D$ for $\alpha<
\sup\limits_{\xi<\xi_0}\mu^\xi$ for some limit ordinal $\xi_0<\cf(\mu_0)$.
The cofinality of the sequence $\langle h_\alpha/D:\alpha<\sup\limits_{\xi<
\xi_0}\mu^\xi\rangle$ is $\cf(\xi_0)<\mu_1$. Consequently the sequence is
bounded in $\prod\limits_{i<\kappa}\lambda_i/D$ and we may proceed as in the
successor case and define $h_\alpha/D$ for $\alpha\in[\sup\limits_{\xi<\xi_0}
\mu^\xi,\mu^{\xi_0})$. 
\medskip

\noindent 2)\ \ \ It follows immediately from 1).\QED

\begin{proposition}
\label{0.6}
Assume that $\inv^{(+)}_T$ is a definable universal well ordered cardinal
invariant. Assume further that: $D$ is an ultrafilter on a cardinal $\kappa$,
for $i < \kappa$, $B_i$ is a Boolean algebra and $B \stackrel{\rm def}{=}
\prod\limits_{i<\kappa}B_i/D$. Then 
\begin{description}
\item[{(a)}] if $\lambda_i<\inv^+_T(B_i)$ for $i<\kappa$ then $\Depth^+
\prod\limits_{i<\kappa}(\lambda_i,<)/D\leq\inv^+_T (B)$, 

\item[{(b)}] $\Depth(\prod\limits_{i<\kappa}(\inv^+_T(B_i),<)/D)\leq\inv^+_T
(B)$, 

\item[{(c)}] if $\inv_T (B)<\Depth\prod\limits_{i<\kappa}(\inv_T(B_i),<)/D$
then for the $D$-majority of $i<\kappa$ we have: $\lambda_i \stackrel{\rm def}
{=}\inv_T(B_i)$ is a limit cardinal and the linear order $\prod\limits_{i<
\kappa}(\lambda_i, <)/D$ is $(\inv_T(B))^+$-Depth-like; hence for the
$D$-majority of the $i<\kappa$ we have: $\lambda_i$ is a regular limit
cardinal, i.e.~weakly inaccesible. 
\end{description}
\end{proposition}

\noindent PROOF: (a)\ \ \ Suppose that $\mu<\Depth^+\prod\limits_{i<\kappa}
(\lambda_i,<)/D$. As $\lambda_i<\inv^+_T(B_i)$ we find $P^i_0,P^i_1,\ldots$
such that $M_i\stackrel{\rm def}{=}(B_i,P_0^i,P^i_1,\ldots)\models T$,
$\|P^i_0\|\geq\lambda_i$. Look at $M\stackrel{\rm def}{=}\prod\limits_{i<
\kappa}M_i/D$. Note that $(P^M_0,P^M_1)$ is a linear ordering such that
\[\Depth^+(P^M_0,P^M_1)\geq\Depth^+\prod_{i<\kappa}(\lambda_i,<)/D.\]
Thus we find $P^*_0\subseteq P^M_0$ such that $\|P^*_0\|=\mu$ and
$(P^*_0,P^M_1)$ is a well ordering. As formulas of $T$ are universal in
$(P_0,P_1)$, first order except the demand that $P_1$ is a well order on $P_0$
we conclude $M^*\stackrel{\rm def}{=}(B,P^*_0,P^M_1,\ldots)\models T$. Hence
$\mu=\|P^*_0\|<\inv^+_T(B)$.
\medskip

\noindent (b)\ \ \ We consider two cases here
\medskip

\noindent{\sc Case 1:} For $D$-majority of $i<\kappa$ we have
$\inv_T(B_i)<\inv_T^+(B_i).$
\medskip

\noindent Then we may assume that for each $i<\kappa$
\[\lambda_i\stackrel{\rm def}{=}\inv_T(B_i)<\inv^+_T(B_i)=\lambda^+_i.\]
By lemma~\ref{depth}(2) we have
\[\Depth(\prod\limits_{i<\kappa}(\lambda_i^+,<)/D)\leq
\Depth^+(\prod\limits_{i<\kappa}(\lambda_i,<)/D).\]
On the other hand, it follows from (a) that
\[\Depth^+(\prod\limits_{i<\kappa}(\lambda_i,<)/D)\leq\inv_T^+(B)\]
and consequently we are done (in this case).
\medskip

\noindent{\sc Case 2:} For $D$-majority of $i<\kappa$ we have
$\inv_T(B_i)=\inv_T^+(B_i).$
\medskip

\noindent So suppose that $\inv_T(B_i)=\inv^+_T(B_i)$ for each $i<\kappa$.
Suppose that 
\[\bar{g}=\langle g_\alpha/D:\alpha<\mu\rangle\subseteq\prod\limits_{i<\kappa}
\inv^+_T(B_i)/D\]
is a $<_D$--increasing sequence.

If $\bar{g}$ is bounded then we apply (a) to conclude that $\mu<\inv_T^+(B)$.
If $\bar{g}$ is unbounded (so cofinal) then there are two possibilities:
either $\mu$ is a limit cardinal or it is a successor. In the first case we
apply the previous argument to initial segments of $\bar{g}$ and we conclude
$\mu\leq\inv^+_T(B)$. In the second case we necessarily have $\mu=\cf(
\prod\limits_{i<\kappa}({\rm \inv}^+_T(B_i),<)/D)= \mu_0^+$ (for some $\mu_0$)
and $\mu_0<\inv^+_T(B)$. Thus $\mu\leq\inv^+_T(B)$. 

Consequently, if there is an increasing (well ordered) sequence of the
length $\mu$ in $\prod\limits_{i<\kappa}(\inv^+_T(B_i),<)/D$ then $\mu\leq
\inv^+_T(B)$ and the case 2 is done too.
\medskip

\noindent (c)\ \ \ Assume that $\lambda\stackrel{\rm def}{=}\inv_T(B)<
\Depth\prod\limits_{i<\kappa}(\inv_T(B_i),<)/D$. By (b) we get that then 
\begin{description}
\item[{$(**)$}]\qquad $\lambda^+=\Depth\prod\limits_{i<\kappa}(\inv_T(B_i),<)/D
=\Depth\prod\limits_{i<\kappa}(\inv^+_T(B_i),<)/D=\inv^+_T(B)$.
\end{description}
Suppose that $\{i<\kappa:\inv_T(B_i)<\inv_T^+(B_i)\}\in D$. Then by (a) we have
\[\Depth^+\prod\limits_{i<\kappa}(\inv_T(B_i),<)/D\leq\inv^+_T(B),\]
but (by $(**)$ and \ref{0.4}) we know that
$$\Depth^+\prod\limits_{i<\kappa}(\inv_T(B_i),<)/D=\lambda^{++}>\inv^+_T(B),$$
a contradiction. Consequently for the $D$-majority of $i<\kappa$ we have
$\lambda_i=\inv_T(B_i)=\inv^+_T(B_i)$ and $\lambda_i$ is a limit cardinal.

Note that if $f\in\prod\limits_{i<\kappa}\inv_T(B_i)$ then
$\Depth^+\prod\limits_{i<\kappa}(f(i),<)/D\leq\lambda^+$ (this is because of
the previous remark, $(**)$ and (a)). Moreover, $(**)$ implies that there is
an increasing sequence $\langle f_\alpha/D:\alpha<\lambda^+\rangle\subseteq
\prod\limits_{i<\kappa}(\inv_T(B_i),<)/D$. By what we noted earlier the
sequence has to be unbounded (so cofinal). Consequently the linear order
$\prod\limits_{i<\kappa}(\inv_T(B_i),<)/D$ is $\lambda^+$--Depth--like. Now
assume that $A=\{i<\kappa:\lambda_i\mbox{ is singular}\}\in D$. Let $Q_i
\subseteq\lambda_i$ be a cofinal subset of $\lambda_i$ of the size
$\cf(\lambda_i)$ (for $i<\kappa$). Then $\Depth^+\prod\limits_{i<\kappa} 
(Q_i,<)/D\leq\lambda^+$ but $\prod\limits_{i<\kappa}Q_i/D$ is cofinal in
$\prod\limits_{i<\kappa}\inv_T(B_i)/D$ - a contradiction, as the last order
has the cofinality $\lambda^+$. \QED

\begin{proposition}
Assume that $\inv^{(+)}_T$ is a definable first order cofinality invariant.
Assume further that: $D$ is an ultrafilter on a cardinal $\kappa$, for
$i<\kappa$, $B_i$ is a Boolean algebra and $B\stackrel{\rm df}{=}
\prod\limits_{i<\kappa} B_i/D$. Then
\begin{description}
\item[{(a)}] if $\lambda_i\in\Inv_T(B_i)$ for $i<\kappa$ and $\lambda=
\cf(\prod\limits_{i<\kappa}(\lambda_i,<)/D)$ then $\lambda\in\Inv_T(B)$, 

\item[{(b)}] if $\inv^+_T(B)\leq\cf(\prod\limits_{i<\kappa }\inv_T(B_i)/D)$
then for the $D$-majority of $i<\kappa$ we have:\ \ \ $\inv_T(B_i)$ is a limit
cardinal. 
\end{description}
\end{proposition}

\noindent PROOF: should be clear. \QED

\begin{proposition}
Suppose that $T$ is a finite $n$-universal in $(P_0,P_1)$ theory in the
language of Boolean algebras plus two predicates $P_0$, $P_1$ and the theory
says that $P_1$ is a linear ordering on $P_0$. Let $\inv^{(+)}_T$ be the
respective cardinality invariant. Assume further that: $D$ is an ultrafilter
on a cardinal $\kappa$, $B_i$ is a Boolean algebra (for $i<\kappa$) and
$B\stackrel{\rm def}{=}\prod\limits_{i<\kappa} B_i/D$. Lastly assume
$\lambda\longrightarrow (\mu)^n_\kappa$, $n\geq 2$ and $\lambda\in
\Inv(B)$.

\noindent Then for the $D$-majority of the $i<\kappa$, $\mu<\inv^+_T(B_i)$.
\end{proposition}

\noindent PROOF: We may assume that $T=\{\psi_0,\psi\}$, where the
sentence $\psi_0$ says ``$P_1$ is a linear ordering of $P_0$'' (and we
denote this ordering by $<$) and
\[\psi=(\forall x_0<\ldots<x_{n-1})(\phi(\bar{x}))\]
where $\phi$ is a formula in the language of Boolean algebras. Note that a
formula 
\[(\forall x_0,\ldots,x_{n-1}\in P_0)(\phi(\bar{x}))\]
as in \ref{basdef}(3) is equivalent to the formula
\[\bigwedge_{f\in {}^n n}(\forall x_0,\ldots,x_{n-1}\in P_0)\big(
[\bigwedge_{f(k)=f(l)} x_k=x_l\ \&\ \bigwedge_{f(k)<f(l)} x_k<x_l]\ \
\Rightarrow\ \ \phi^f_0(\bar{x})\big),\]
where, for any $f:n\longrightarrow n$, the formula $\phi_0^f$ is obtained from
$\phi$ by respective replacing appearances of $P_1(x_i,x_j)$ by either
$x_i=x_i$ or $x_i\neq x_i$. Consequently the above assumption is easily
justified. 

Let $A=\{i<\kappa:\mu<{\rm \inv}^+_T(B_i)\}$. Assume that $A\notin D$. As
$\lambda\in {\rm Inv}_T(B)$ we find $P_0,P_1$ such that $\|P_0\|=\lambda$ and
$P_1=<$ is a linear ordering of $P_0$ and $(B,P_0,P_1)\models\psi$. For each
element of $\prod\limits_{i<\kappa} B_i/D$ we fix a representative of this
equivalence class (so we will freely pass from $f/D$ to $f$ with no additional
comments). Now, we define a colouring $F:[P_0]^n\longrightarrow\kappa$ by 
\[\begin{array}{ll}
F(f_0/D,\ldots,f_{n-1}/D)=&\mbox{the first }i\in\kappa\setminus A\mbox{ such
that}\\
\ &\mbox{if }f_0/D<\ldots<f_{n-1}/D\\
\ &\mbox{then }f_0(i),\ldots,f_{n-1}(i)\mbox{ are pairwise distinct and}\\
\ &B_i\models\phi[f_0(i),\ldots,f_{n-1}(i)]\\
\end{array}\]
The respective $i$ exists since $A\notin D$ and
\[B\models\mbox{``}f_0/D,\ldots,f_{n-1}/D\mbox{ are distinct and } \phi[f_0/D,
\ldots,f_{n-1}/D]\mbox{''}.\]
By the assumption $\lambda\longrightarrow(\mu)^n_\kappa$ we find $W\in
[P_0]^\mu$ which is homogeneous for $F$. Let $i$ be the constant value of $F$
on $W$ and put\quad $P^i_0=\{f(i): f/D\in W\}$\quad (recall that we fixed
representatives of the equivalence classes). Now we may introduce $P^i_1$ as
the linear ordering of $P^i_0$ induced by $P_1$. 

\noindent Note that $f(i)\neq f'(i)$ for distinct $f/D,f'/D\in W$ and if
$f_0(i),\ldots,f_{n-1}(i)\in P^i_0$, $f_0(i)<_{P^i_1} f_1(i)<_{P^i_1}\ldots
<_{P^i_1} f_{n-1}(i)$ then $f_0/D<\ldots<f_{n-1}/D$ and 
hence 
\[B_i\models\phi[f_0(i),\ldots,f_{n-1}(i)].\]
As $\|P^i_0\|=\mu$, $(B_i,P^i_0, P^i_1)\models\psi\wedge\psi_0$ we
conclude that $\mu<{\rm \inv}^+_T(B_i)$ what contradicts $i\notin
A$. \QED
\medskip

One of the tools in study the invariants are ``finite'' versions of them (for
invariants determined by infinite theories). Suppose $T=\{\phi_n:n<\omega\}$
and if $T$ is supposed to describe a def.u.w.o.car.~invariant then $\phi_0$
already says that $P_1$ is a well ordering of $P_0$. Let $T^n=\{\phi_m: m<n\}$
for $n<\omega$. 

\begin{conclusion}
\label{finite}
Suppose that $D$ is a uniform ultrafilter on $\kappa$,
$f:\kappa\longrightarrow\omega$ is such that $\lim\limits_D f=\omega$. Let
$B_i$ (for $i<\kappa$) be Boolean algebras, $B=\prod\limits_{i<\kappa} B_i/D$.
\begin{enumerate}
\item If $T$ is first order then:
\begin{description}
\item[a)\ ] if $\lambda_i\in\Inv_{T^{f(i)}}(B_i)$ (for $i<\kappa$) then
$\prod\limits_{i<\kappa}\lambda_i/D\in\Inv_T(B)$, 
\item[b)\ ] $\prod\limits_{i<\kappa}\inv^+_{T^{f(i)}}(B_i)/D\leq \inv^+_T(B)$.
\end{description}
\item If $T$ is u.w.o. then:
\begin{description}
\item[a)\ ] if $\lambda_i\in\Inv_{T^{f(i)}}(B_i)$ (for $i<\kappa$)
and $\lambda<\Depth^+\prod\limits_{i<\kappa}(\lambda_i,<)/D$\\
then $\lambda\in\Inv_T(B)$,
\item[b)\ ] $\Depth\prod\limits_{i<\kappa}(\inv^+_{T^{f(i)}}(B_i),<)/D\leq
\inv^+_T(B)$.
\end{description}
\end{enumerate}
\end{conclusion}

\noindent PROOF: Like \ref{fo} and \ref{0.6}. \QED

\subsection{An example concerning the question $(<)_{\inv}$.}
Now we are going to show how the main result of \cite{MgSh 433} may be used to
give affirmative answers to the questions $(<)_{\inv}$ for several cardinal
invariants.
\begin{proposition}
\label{product}
Suppose that $D$ is an $\aleph_1$-complete ultrafilter on $\kappa$,
$B_{i,\alpha}$ are Boolean algebras (for $\alpha<\lambda_i$, $i<\kappa$). Let
$C:\prod\limits_{i<\kappa}\lambda_i/D \longrightarrow\prod\limits_{i<\kappa}
\lambda_i$ be a choice function (so $C(x)\in x$ for an equivalence class
$x\in\prod\limits_{i<\kappa}\lambda_i/D$). 
\begin{enumerate}
\item If $B_i=\cast\limits_{\alpha<\lambda_i}B_{i,\alpha}$ then
\[\prod\limits_{i<\kappa}B_i/D\simeq\cast\{\prod\limits_{i<\kappa}
B_{i,C(x)(i)}/D: x\in\prod\limits_{i<\kappa}\lambda_i/D\}.\]
\item If $B_i=\prod_{\alpha<\lambda_i}^w B_{i,\alpha}$ then
\[\prod\limits_{i<\kappa}B_i/D\simeq{\prod}^w\{\prod\limits_{i<\kappa}
B_{i,C(x)(i)}/D: x\in\prod\limits_{i<\kappa}\lambda_i/D\}.\quad \QED\]
\end{enumerate}
\end{proposition}

\begin{definition}
Let $\op$ be an operation on Boolean algebras.
\begin{enumerate}
\item For a theory $T$ we define the property $\square^T_{\op}$:
\begin{description}
\item[$\square^T_{\op}$] if $\mu$ is a cardinal, $B_i$ are Boolean algebras
for $i<\mu^+$ then
\[\hspace{-1cm}\sup_{i<\mu}\inv_T(B_i)\leq\inv_T(\op_{i<\mu}B_i)\ \ \mbox{ and
}\ \ \inv_T(\op_{i<\mu^+}B_i) \leq\mu+\sup_{i<\mu^+}\inv_T(B_i).\]
\end{description}
\item Of course we may define the respective property for any cardinal
invariant (not necessary of the form $\inv_T$). But then we additionally
demand that $\tau(B)\leq\|B\|$ (where $\tau$ is the considered invariant). 
\end{enumerate}
\end{definition}

\begin{proposition}
\label{general}
Suppose that a def.car.invariant $\inv_T$ (or just an invariant $\tau$)
satisfies either $\square^T_{\cast}$ or $\square^T_{\prod^w}$ and suppose
that for each cardinal $\chi$ there is a Boolean algebra $B$ such that
$\chi\leq\inv_T(B)$ and there is no weakly inaccessible cardinal in the
interval $(\chi,\|B\|]$. Assume further that 
\begin{description}
\item[$(\odot)$] $\langle\lambda_i:i<\kappa\rangle$ is a sequence of
weakly inaccessible cardinals, $\lambda_i>\kappa^+$, $D$ is an
$\aleph_1$-complete ultrafilter on $\kappa$ and
$\prod\limits_{i<\kappa}(\lambda_i,<)/D$ is $\mu^+$-like (for some cardinal
$\mu$).  
\end{description}
Then there exist Boolean algebras $B_i$ for $i<\kappa$ such that
$\inv_T(B_i)=\lambda_i$ (for $i<\kappa$) and $\inv_T(\prod\limits_{i<\kappa}
B_i/D)\leq\mu$. So we have 
\[\prod_{i<\kappa}\inv_T(B_i)/D=\mu^+>\inv_T(\prod_{i<\kappa}B_i/D).\]
\end{proposition}

\noindent PROOF: Assume that $\inv_T$ satisfies $\square^T_{\cast}$. For
$i<\kappa$ and $\alpha<\lambda_i$ fix an algebra $B_{i,\alpha}$ such that 
\[\|\alpha\|\leq\inv_T(B_{i,\alpha})\leq\|B_{i,\alpha}\|<\lambda_i\]
(possible by our assumptions on $\inv_T$) and let $B_i=\cast\limits_{\alpha<
\lambda_i}B_{i,\alpha}$. By \ref{product} we have 
\[\prod\limits_{i<\kappa}B_i/D=\cast\{\prod\limits_{i<\kappa}B_{i,C(x)(i)}/D:
x\in\prod\limits_{i<\kappa}\lambda_i/D\},\]
where $C:\prod\limits_{i<\kappa}\lambda_i/D\longrightarrow\prod\limits_{i<
\kappa}\lambda_i$ is a choice function. So by $\square^T_{\cast}$ (the second
inequality):  
\[\inv_T(\prod\limits_{i<\kappa}B_i/D)\leq\mu+\sup\{\inv_T(\prod\limits_{i<
\kappa}B_{i,C(x)(i)}/D): x\in\prod\limits_{i<\kappa}\lambda_i/D\}.\] 
Since $\|B_{i,\alpha}\|<\lambda_i$ and $\prod\limits_{i<\kappa}(\lambda_i,<
)/D$ is $\mu^+$-like for each $x\in\prod\limits_{i<\kappa}\lambda_i/D$ we have
\[\inv_T(\prod\limits_{i<\kappa}B_{i,C(x)(i)}/D)\leq\prod\limits_{i<\kappa}
\|B_{i,C(x)(i)}\|/D\leq\mu.\]  
Moreover, by the first inequality of $\square^T_{\cast}$, for each
$\alpha<\lambda_i$ 
\[\|\alpha\|\leq\inv_T(B_{i,\alpha})\leq \inv_T(B_i)\leq\|B_i\|=\lambda_i\]
and thus $\inv_T(B_i)=\lambda_i$. \QED
\medskip

\noindent{\bf Remark:}\ \ \ 1.\ \ The consistency of $(\odot)$ is the
main result of \cite{MgSh 433}, where several variants of it and their
applications are presented.
 
\noindent 2.\ \ If $\inv_T$ is either def.f.o.car invariant or def.u.w.o.car
invariant then we may apply \ref{fo}.c or \ref{0.6}.c respectively to
conclude that for D-majority of $i<\kappa$ we have $\inv(B_i)=\inv^+(B_i)$.
Consequently in these cases we may slightly modify the construction in
\ref{general} to get additionally $\inv(B_i)=\inv^+(B_i)$ for each $i<\kappa$.

\noindent 3.\ \ Proposition~\ref{general} applies to several cardinal
invariants. For example the condition $\square^T_{\prod^w}$ is satisfied by:
\begin{quotation}
\noindent Depth (see \S 4 of \cite{M2}), Length (\S 7 of \cite{M2}), Ind
(\S 10 of \cite{M2}), $\pi$-character (\S 11 of \cite{M2}) and the
tightness t (\S 12 of \cite{M2}). 
\end{quotation}
Moreover, \ref{general} can be applied to the topological density $d$, as this
cardinal invariant satisfies the corresponding condition $\square^d_{\cast}$.
[Note that $d(\cast\limits_{i<\mu^+}B_i)=\max\{\lambda,\sup\limits_{i<\mu^+}
d(B_i)\}$, where $\lambda$ is the first cardinal such that $\mu^+\leq
2^\lambda$, so $\lambda\leq\mu$; see \S 5 of \cite{M2}.]

\section{Topological density} 
The topological density of a Boolean algebra $B$ (i.e. the density of its
Stone space $\Ult B$) equals to $\min\{\kappa:B\mbox{ is }\kappa
\mbox{-centered}\}$. To describe it as a minimal definable first order
cardinality invariant we use the theory defined below. 

\begin{definition}
\begin{enumerate}
\item  For $n<\omega$ define the formulas $\phi^d_n$ by:
\[\hspace{-0.9cm}\phi^d_0=(\forall x)(\exists y\!\in\! P_0)(x\neq 0\
\Rightarrow\ P_1(y,x))\ \&\ (\forall x)(\forall y\!\in\! P_0)(P_1(y,x)\
\Rightarrow\ x\neq 0)\]  
and for $n>0$:
\[\hspace{-0.9cm}\phi^d_n=(\forall x_0,\ldots,x_n)(\forall y\in P_0)(
P_1(y,x_0)\ \&\ldots\&\ P_1(y,x_n)\ \Rightarrow\ x_0\wedge\ldots\wedge x_n\neq
0).\] 
\item For $n\leq\omega$ let $T_d^n=\{\phi_k: k< n\}$.
\item For a Boolean algebra $B$, $n\leq \omega$ we put $d_n(B)=\min\Inv_{
T_d^n}(B)$. 
\item For $1\leq n<\omega$, a subset $X$ of a Boolean algebra $B$ has {\em the
$n$--intersection property} provided that the meet of any $n$ elements of $X$
is nonzero; if $X$ has the $n$--intersection property for all $n$, then $X$ is
{\em centered}, or has {\em the finite intersection property}. 
\end{enumerate}
\end{definition}
Note that $d_\omega(B)$ is the topological density $d(B)$ of $B$. Since $T^0_d
=\emptyset$, the invariant $d_0(B)$ is just $0$. The theory $T^{n+1}_d$ says
that for each $y\in P_0$ the set $X_y\stackrel{\rm def}{=}\{x: P_1(y,x)\}$ has
the $n+1$--intersection property and $\bigcup\limits_{y\in P_0}X_y=B\setminus
\{0\}$. Thus, for $1\leq n<\omega$, $d_n(B)$ is the smallest cardinal $\kappa$
such that $B\setminus\{0\}$ is the union of $\kappa$ sets having the
$n$--intersection property. 
\medskip

We easily get (like \ref{finite}):

\begin{fact}
\label{fa2.1}
\begin{enumerate}
\item For a Boolean algebra $B$, the sequence $\langle d_n(B): 1\leq n\leq
\omega\rangle$ is increasing and $d(B)\leq\prod\limits_{1\leq n<\omega}d_n(B)$.
\item If $D$ is an ultrafilter on a cardinal $\kappa$, $f:\kappa
\longrightarrow\omega$ is a function such that $\lim\limits_{D}f=\omega$ and
$B_i$ (for $i<\kappa$) are Boolean algebras\\
then $d(\prod\limits_{i<\kappa}B_i/D)\leq\prod\limits_{i<\kappa}d_{f(i)}
(B_i)/D$. \QED 
\end{enumerate}
\end{fact}

\begin{fact}
\label{factdense}
\begin{enumerate}
\item If $1\leq n<\omega$ and $X$ is a dense subset of $B\setminus\{0\}$, then
$d_n(B)$ is the least cardinal $\kappa$ such that $X$ can be written as a
union of $\kappa$ sets each with the $n$--intersection property.
\item If $X$ is a dense subset of $B\setminus\{0\}$, then $d_\omega(B)$ is the
least cardinal $\kappa$ such that $X$ can be written as a union of $\kappa$
sets each with the finite intersection property.
\item If $B$ is an interval Boolean algebra then $d_2(B)=d(B)$.
\end{enumerate}
\end{fact}

\noindent PROOF: Suppose $X\subseteq B\setminus\{0\}$ is dense, $1\leq n<
\omega$. Obviously $X$ can be written as a union of $d_n(B)$ sets each with
the $n$--intersection property. If $X=\bigcup\limits_{i<\kappa}Y_i$, where
$Y_i$ have the $n$--intersection property, let $Z_i\stackrel{\rm def}{=}\{b\in
B:(\exists y\in Y_i)(y\leq b)\}$. Then each $Z_i$ has the $n$--intersection
property and $B\setminus\{0\}=\bigcup\limits_{i<\kappa}$. This proves
condition 1; condition 2 is proved similarly. Condition 3 follows since for an
interval algebra $B$ intervals are dense in $B$ and if $a_1,\ldots,a_k$ are
intervals such that $a_i\wedge a_j\neq 0$ then $\bigwedge\limits^{k}_{i=1}
a_i\neq 0$. \QED 
\medskip

A natural question that arises here is if we can distinguish the invariants
$d_n$. The positive answer is given by the examples below. 

\begin{example}
\label{ex2.3}
Let $\kappa$ be an infinite cardinal, $n>2$. There is a Boolean
algebra $B$ such that $d_n(B)>\kappa$, $d_{n-1}(B)\leq 2^{<\kappa}$.
\end{example}

\noindent PROOF: Let $B$ be the Boolean algebra generated freely by
$\{x_\eta:\eta\in {}^{\kappa}n\}$ except:
\begin{quotation}
\noindent if $\nu\in {}^{\kappa>} n$, $\nu\hat{\ }\langle l\rangle\subseteq
\eta_l\in{}^\kappa n$ (for $l<n$) 

\noindent then $x_{\eta_0}\wedge\ldots\wedge x_{\eta_{n-1}}=0$.
\end{quotation}
Suppose that $B^+=\bigcup\limits_{i<\kappa} D_i$. For $\eta\in {}^\kappa n$
let $i(\eta)<\kappa$ be such that $x_\eta\in D_{i(\eta)}$. Now we inductively
try to define $\eta^*\in {}^\kappa n$:
\begin{quotation}
\noindent assume that we have defined $\eta^*\rest i$ ($i<\kappa$) and we want
to choose $\eta^*(i)$. If there is $l<n$ such that $i(\eta)\neq i$ for each
$\eta\supseteq\eta^*\rest i\hat{\ }\langle l\rangle$ then we choose one such
$l$ and put $\eta^*(i)=l$. If there is no such $l$ then we stop our
construction. 
\end{quotation}
If the construction was stopped at stage $i<\kappa$ (i.e.~we were not able to
choose $\eta^*(i)$) then for each $l<n$ we have a sequence $\eta_l\in
{}^\kappa n$ such that $\eta^*\rest i\hat{\ }\langle l\rangle\subseteq\eta_l$ 
and $i(\eta_l)=i$. Thus $x_{\eta_0},\ldots,x_{\eta_{n-1}}\in D_i$ and
$x_{\eta_0}\wedge\ldots\wedge x_{\eta_{n-1}}=0$, so that $D_i$ does not
satisfy the $n$--intersection property. If we could carry our construction up
to $\kappa$ then we would get $\eta^*\in{}^\kappa n$ such that $x_{\eta^*}
\notin\bigcup\limits_{i<\kappa} D_i$. Consequently the procedure had to stop
and we have proved that $d_n(B)>\kappa$. 

Now we are going to show that $d_{n-1}(B)\leq 2^{<\kappa}$. Let $X$ be the set
of all nonzero elements of $B$ of the form
\[x_{\eta_0}\wedge\ldots\wedge x_{\eta_l}\wedge (-x_{\eta_{l+1}})\wedge\ldots
\wedge (-x_{\eta_k})\]
in which the sequences $\eta_0,\ldots,\eta_k\in {}^\kappa n$ are pairwise
distinct, $0<l<k<\omega$. Clearly $X$ is dense in $B$. We are going to apply
fact \ref{factdense}(1). To this end, if $0<l<k$, $\alpha<\kappa$, and
$\langle\nu_0,\ldots,\nu_k\rangle$ is a sequence of distinct members of
${}^\alpha n$, let $D^{l,k,\alpha}_{\langle\nu_0,\ldots,\nu_k\rangle}$ be the
set
\[\{x_{\eta_0}\wedge\ldots\wedge x_{\eta_l}\wedge(-x_{\eta_{l+1}})\wedge\ldots
\wedge (-x_{\eta_k}):\nu_0\subseteq\eta_0\in {}^\kappa n,\ldots,\nu_k\subseteq
\eta_k\in {}^\kappa n\}\setminus\{0\}.\]
Note that $X$ is the union of all these sets. There are $2^{<\kappa}$
possibilities for the parameters, so it suffices to show that each of the sets
$D^{l,k,\alpha}_{\langle\nu_0,\ldots,\nu_k\rangle}$ has the
$(n-1)$--intersection property. 

Before beginning on this, note that if $\eta_0,\ldots,\eta_k\in {}^\kappa n$
are such that $\eta_i\neq\eta_j$ when $i\leq l<j\leq k$ and
\[B\models x_{\eta_0}\wedge\ldots\wedge x_{\eta_l}\wedge (-x_{\eta_{l+1}})
\wedge\ldots\wedge (-x_{\eta_k})=0,\]
then necessarily there is $\nu\in {}^{<\kappa} n$ such that 
\[(\forall m<n)(\exists i\leq l)(\nu\hat{\ }\langle
m\rangle\subseteq\eta_l).\] 
Now we check that $D^{l,k,\alpha}_{\langle\nu_0,\ldots,\nu_k\rangle}$ has the
$(n-1)$--intersection property, where $0<l<k<\omega$, $\alpha<\kappa$, and
$\nu_0,\ldots,\nu_k$ are pairwise distinct elements of ${}^\alpha n$. Thus
suppose that
\[x_{\eta^j_0}\wedge\ldots\wedge x_{\eta^j_l}\wedge (-x_{\eta^j_{l+1}})\wedge
\ldots\wedge (-x_{\eta^j_k})\]
are members of $D^{l,k,\alpha}_{\langle\nu_0,\ldots,\nu_k\rangle}$ for each
$j<n-1$; and suppose that 
\[B\models\bigwedge_{j<n-1} x_{\eta^j_0}\wedge\ldots\wedge\bigwedge_{j<n-1}
x_{\eta^j_l}\wedge\bigwedge_{j<n-1}(-x_{\eta^j_{l+1}})\wedge\ldots\wedge
\bigwedge_{j<n-1}(-x_{\eta^j_k})=0.\]
By the above remark, choose $\nu\in {}^{<\kappa} n$ such that for all $m<n$
there exist an $i(m)\leq l$ and a $j(m)<n-1$ such that $\nu\hat{\ }\langle
m\rangle\subseteq\eta^{j(m)}_{i(m)}$ (note that if $j_0,j_1<n-1$, $i_0\leq l$
and $l+1\leq i_1\leq k$ then $\eta^{j_0}_{i_0}\neq\eta^{j_1}_{i_1}$ as $\nu_0,
\ldots,\nu_k$ are pairwise distinct). 
\medskip

\noindent{\sc Case 1:}\ \ \ $\nu_i\subseteq\nu$ for some $i\leq k$.

\noindent Then for each $m<n$ we have $\nu_i\subseteq\nu\subseteq\nu\hat{\ }
\langle m\rangle\subseteq\eta^{j(m)}_{i(m)}$ and consequently $i(m)=i$ (for
$m<n$). As $j(m)<n-1$ for $m<n$ we find $m_0<m_1<n-1$ such that
$j(m_0)=j(m_1)=j$. Then $\nu\hat{\ }\langle m_0\rangle\subseteq\eta^j_i$,
$\nu\hat{\ }\langle m_1\rangle\subseteq\eta^j_i$ give a contradiction.
\medskip

\noindent{\sc Case 2:}\ \ \ $\nu_i\not\subseteq\nu$ for all $i\leq k$.

\noindent Note that for all $m<n$ the sequences $\nu\hat{\ }\langle m\rangle$
and $\nu_{i(m)}$ are compatible. By the case we are in, it follows that $\nu$
is shorter than $\nu_{i(m)}$. So $\nu\hat{\ }\langle m\rangle\subseteq
\nu_{i(m)}$, $i(m)<l$. But then by construction, $D^{l,k,\alpha}_{\langle
\nu_0,\ldots\nu_k\rangle}$ is empty, a contradiction. \QED

\begin{example}
\label{ex2.4}
Let $\lambda_i$ be cardinals (for $i<\kappa$) such that $2^\kappa<
\prod\limits_{i<\kappa}\lambda_i$, $2<n<\omega$. Then there is a Boolean
algebra $B$ such that 
\[d_{n-1}(B)\leq\sum_{\alpha<\kappa}\prod\limits_{i<\alpha}\lambda_i\ \ \mbox{
and }\ \ d_n(B)=\|B\|=\prod\limits_{i<\kappa}\lambda_i.\] 
In particular, if $\lambda$ is a strong limit cardinal, $\cf(\lambda)<
\lambda$, $2< n<\omega$ then there is a Boolean algebra $B$ such that  
$d_n(B)=\|B\|=2^\lambda$, $d_{n-1}(B)\leq\lambda$.
\end{example}

\noindent PROOF: Let $B$ be the Boolean algebra generated freely by
$\{x_\eta: \eta\in\prod\limits_{i<\kappa}\lambda_i\}$ except that:
\begin{quotation}
\noindent if $\alpha<\kappa$, $v\in\prod\limits_{i<\alpha}\lambda_i$,
$v\subseteq\eta_l\in\prod\limits_{i<\kappa}\lambda_i$, $\|\{\eta_l(\alpha):
l<n\}\|=n$

\noindent then $x_{\eta_0}\wedge\ldots\wedge x_{\eta_{n-1}}=0$
\end{quotation}
The same arguments as in the previous example show that
\[d_{n-1}(B)\leq\sum_{\alpha<\kappa}\prod\limits_{i<\alpha}\lambda_i.\]
Suppose now that $\prod\limits_{i<\kappa}\lambda_i=\bigcup\{D_j: j<\theta\}$,
$\theta<\prod\limits_{i<\kappa}\lambda_i$ and if $\eta_0,\ldots,\eta_{n-1}\in
D_j$, $j<\theta$ then $x_{\eta_0}\wedge\ldots\wedge x_{\eta_{n-1}}\neq 0$.
Thus the trees $T_j=\{\eta\rest\alpha: \alpha<\kappa,\eta\in D_j\}$ have no
splitting into more than $n-1$ points and hence $\|D_j\|\leq n^\kappa<
\prod\limits_{i<\kappa}\lambda_i$ for all $j<\theta$ and we get a
contradiction, proving $d_n(B)=\prod\limits_{i<\kappa}\lambda_i$. \QED 

\begin{corollary}
\label{prJ}
Let $\lambda$ be a strong limit cardinal, $\kappa<\cf(\lambda)<\lambda$.
Suppose that $D$ is an ultrafilter on $\kappa$ which is not
$\aleph_1$-complete. Then there exist Boolean algebras $B_i$ (for $i<\kappa$)
such that 
\[d(\prod\limits_{i<\kappa}B_i/D)\leq\lambda<2^\lambda=\prod\limits_{i<\kappa}
d(B_i)/D.\]
\end{corollary}

\noindent PROOF: As $D$ is not $\aleph_1$-complete we find a function $f:
\kappa\longrightarrow\omega\setminus 2$ such that $\lim\limits_{D}f=\omega$.
Let $B_i$ be such that $\|B_i\|=d_{f(i)+1}(B_i)=2^\lambda$, $d_{f(i)}(B_i)\leq
\lambda$ (see \ref{ex2.4}). Then, by \ref{fa2.1}, we have $d(\prod\limits_{i<
\kappa}B_i/D)\leq\prod\limits_{i<\kappa}d_{f(i)}(B_i)/D\leq\lambda^\kappa=
\lambda$. As $d(B_i)=d_{f(i)+1}(B_i)= 2^\lambda$ we have $\prod\limits_{i<
\kappa}d(B_i)/D=2^\lambda$. \QED
\medskip

\noindent{\bf Remark:}\ \ \ 1.\ Corollary \ref{prJ} applied e.g. to
$\lambda=\beth_{\omega_1}$, $\kappa=\omega$ gives a negative answer to
Problem J of \cite{M3}.

\noindent 2.\ The algebras $B_i$ in \ref{prJ} are of a quite large
size: $\|B_i\|=2^\lambda$, $\lambda$ strong limit of the cofinality
$>\kappa$. Moreover the cardinal $\lambda$ had to be singular. The
natural question if these are real limitations is answered by the
theorem below. This example, though more complicated than the
previous ones, has several nice properties. E.g. it produces algebras of
the size $2^{(2^{\aleph_0})^+}$ already.

\begin{theorem}
\label{th2.5}
Assume that $\theta=\cf(\theta)$, $\theta^\kappa=\theta$. Then there are
Boolean algebras $B_\gamma$ for $\gamma<\kappa$ such that $d(B_\gamma)=
d_2(B_\gamma)=\theta^+$ and $d(\prod\limits_{\gamma<\kappa}B_\gamma/D)\leq
\theta$ for every non-principal ultrafilter $D$ on $\kappa$. 
\end{theorem}

\noindent PROOF: Let $\lambda=2^\theta$. Choose $\eta_{\alpha,i}\in {}^\theta
\theta$ for $\alpha<\lambda$, $i<\theta^+$ such that
\begin{description}
\item[1. ] if $\eta_{\alpha_1,i_1}=\eta_{\alpha_2,i_2}$ then
$(\alpha_1,i_1)=(\alpha_2,i_2)$,
\item[2. ] for each $f\in {}^\theta\theta$ and $i<\theta^+$ the set
$\{\alpha<\lambda: (\forall\varepsilon<\theta)(f(\varepsilon)\leq
\eta_{\alpha,i}(\varepsilon))\}$ is of the size $\lambda$
\end{description}
(the choice is possible as there is $2^\theta=\lambda$ pairs $(f,i)$ to take
care of and for each such pair we have $2^\theta$ candidates for
$\eta_{\alpha,i}$). 

For two functions $f,g\in {}^\theta\theta$ we write $f<^*g$ if and only if
\[\|\{\varepsilon<\theta: f(\varepsilon)\geq g(\varepsilon)\}\|<\theta.\]
We say that a set $A\subseteq\lambda\times\theta^+$ is $i$-large (for $i<
\theta^+$) if for every $f\in {}^\theta\theta$ we have $\|\{\alpha<\lambda:
(\alpha,i)\in A\ \&\ f<^*\eta_{\alpha,i}\}\|=\lambda$ and we say that $A$ is
large if $\sup\{i<\theta^+: A \mbox{ is } i\mbox{-large}\}=\theta^+$.
\begin{claim}
\label{cl1}
The union of at most $\theta$ sets which are not large is not large.
\end{claim}
Proof of the claim: Should be clear as $\cf(\theta)=\theta<\cf(\lambda)$.
\medskip

Now we are going to describe the construction of the Boolean algebras we need.
First suppose that $S\subseteq\{j<\theta^+:\cf(j)=\theta\}$ is a stationary
set and let $S^+=\{(\alpha,j)\in\lambda\times\theta^+: j\in S\}$. Now choose a
sequence $\bar{F}=\langle F_\varepsilon:\ve<\theta\rangle$ such that:
\begin{description}
\item[3. ] $F_\ve$ is a function with the domain $\dom(F_\ve)=S^+$,
\item[4. ] if $i\in S$, $\alpha<\lambda$ then $F_\ve(\alpha,i)=(F_{\ve,1}
(\alpha,i),F_{\ve,2}(\alpha,i))\in\lambda\times i$,
\item[5. ] if $i\in S$, $\alpha<\lambda$ then the sequence $\langle F_{\ve,2}
(\alpha,i):\ve<\theta\rangle$ is strictly increasing with the limit $i$,
\item[6. ] if $\langle A_\ve:\ve<\theta\rangle$ is a sequence of large
subsets of $\lambda\times\theta^+$ then for some stationary set $S'\subseteq
S$ for each $i\in S'$, $f\in {}^\theta\theta$ we have 
\[\|\{\alpha<\lambda: f<^*\eta_{\alpha,i}\ \&\ (\forall\ve<\theta)
(F_\ve(\alpha,i)\in A_\ve)\}\|=\lambda,\] 
\item[7. ] if $\ve<\zeta<\theta$ then $\eta_{F_\ve(\alpha,i)}<^*\eta_{F_\zeta
(\alpha,i)}$.
\end{description}
To construct the sequence $\bar{F}$ fix $i\in S$. Let $\{(f_\alpha,g_\alpha,
\bar{j}_\alpha):\alpha<\lambda\}$ enumerate with $\lambda$-repetitions all
triples $(f,g,\bar{j})$ such that $f\in {}^\theta\theta$, $\bar{j}=\langle
j_\ve:\ve<\theta\rangle$ is an increasing cofinal sequence in $i$ and $g\in
{}^\theta\lambda$ is such that 
\begin{description}
\item[$(*)$]\qquad\qquad $\ve<\zeta<\theta\ \ \Rightarrow\ \
\eta_{g(\ve),j_\ve}<^*\eta_{g(\zeta),j_\zeta}$
\end{description}
(recall that $\cf(i)=\theta$, $\lambda=2^\theta$). Now we inductively choose
$\langle\beta_\alpha:\alpha<\lambda\rangle\subseteq\lambda$ such that
$\beta_\alpha\notin\{\beta_\delta:\delta<\alpha\}$,
$f_\alpha<^*\eta_{\beta_\alpha,i}$ (this is possible by (2)). Finally for
$\alpha<\lambda$ and $\ve<\theta$ define $F_\ve(\alpha,i)$ by:
\begin{quotation}
\noindent if $\alpha=\beta_\delta$ for some $\delta<\lambda$ then
$F_\ve(\alpha,i)=(g_\delta(\ve),j^\delta_\ve)$,

\noindent if $\alpha\notin\{\beta_\delta:\delta<\lambda\}$ then
$F_\ve(\alpha,i)=(g_\alpha(\ve),j^\alpha_\ve)$
\end{quotation}
(where $\bar{j}_\delta=\langle j^\delta_\ve:\ve<\theta\rangle$). Easily
conditions (3)--(5) and (7) are satisfied. To check clause (6) suppose that
$\langle A_\ve:\ve<\theta\rangle$ is a sequence of large sets and let $S'$ be
the set of all $i\in S$ such that there exists an increasing cofinal sequence
$\langle j_\ve:\ve<\theta\rangle\subseteq i$ such that $A_\ve$ is
$j_\ve$-large (for each $\ve<\theta$). The set $S'$ is stationary. [Why? For 
$\ve<\theta$ let $C_\ve$ be the set of all points in $\theta^+$ which are
limits of increasing sequences from $\{j<\theta^+: A_\ve$ is $j$--large$\}$.
Clearly each $C_\ve$ is a club of $\theta^+$ and thus $\bigcap\limits_{\ve<
\theta} C_\ve$ is a club of $\theta^+$. Now one easily checks that $S\cap
\bigcap\limits_{\ve<\theta} C_\ve \subseteq S'$.] 

\noindent We are going to show that $S'$ works for $\langle A_\ve:\ve<\theta
\rangle$. Take $i\in S'$ and suppose that $f\in {}^\theta\theta$. Let
$\bar{j}=\langle j_\ve:\ve<\theta\rangle\subseteq i$ be an increasing cofinal
sequence witnessing $i\in S'$. Take $g\in {}^\theta\lambda$ such that
\[\ve<\zeta<\theta\ \Rightarrow\ [\eta_{g(\ve),j_\ve}<^*\eta_{g(\zeta),
j_\zeta}\ \&\ (g(\ve),j_\ve)\in A_\ve]\]
(possible by the $j_\ve$-largeness of $A_\ve$ and the regularity of $\theta$).
When we defined $F_\ve(\alpha,i)$ (for $\ve<\theta$, $\alpha<\lambda$), the
triple $(f,g,\bar{j})$ appeared $\lambda$ times in the enumeration
$\{(f_\alpha,g_\alpha,\bar{j}_\alpha):\alpha<\lambda\}$. Whenever
$(f,g,\bar{j})=(f_\alpha,g_\alpha,\bar{j}_\alpha)$ we had $F_\ve(\beta_\alpha,
i)=(g_\alpha(\ve),j^\alpha_\ve)=(g(\ve),j_\ve)\in A_\ve$ and $f=f_\alpha<^*
\eta_{\beta_\alpha,i}$. Consequently if $i\in S'$, $f\in {}^\theta\theta$ then
\[\|\{\alpha<\lambda:f<^*\eta_{\alpha,i}\ \&\ (\forall\ve<\theta)(F_\ve(
\alpha,i)\in A_\ve)\}\|=\lambda\] 
and condition (6) holds.
\medskip

For the sequence $\bar{F}$ we define a Boolean algebra $B_{\bar{F}}$: it is
freely generated by $\{x_{\alpha,i}:\alpha<\lambda, i<\theta^+\}$ except that
\begin{quotation}
\noindent if $F_\ve(\alpha_1,i_1)=(\alpha_2,i_2)$ for some $\ve<\theta$ 

\noindent then $x_{\alpha_1,i_1}\wedge x_{\alpha_2,i_2}=0$.
\end{quotation}
Now fix a sequence $\langle S_\gamma:\gamma<\kappa\rangle$ of pairwise
disjoint stationary subsets of $\{j\in\theta^+: \cf(j)=\theta\}$ and for each
$\gamma<\kappa$ fix a sequence $\bar{F}_\gamma=\langle F^\gamma_\ve:\ve<\theta
\rangle$ satisfying conditions (3)--(7) above (for $S_\gamma$). 
\begin{claim}
For each $\gamma<\kappa$, $d_2(B_{\bar{F}_\gamma})>\theta$.
\end{claim}
Proof of the claim: Let $\bar{F}=\bar{F}_\gamma$ and suppose that
$B^+_{\bar{F}}=\bigcup\limits_{\ve<\theta}D_\ve$. Let $A_\ve=\{(\alpha,i): 
x_{\alpha,i}\in D_\ve\}$ and let $A_\ve^\prime=A_\ve$ if $A_\ve$ is large and
$A_\ve^\prime=\lambda\times\theta^+$ otherwise. So the sets $A_\ve^\prime$ are
large (for $\ve<\theta$) and by condition (6) the set 
\[A\stackrel{\rm def}{=}\{(\alpha,i)\in\lambda\times\theta^+:(\forall\ve<
\theta)(F_\ve(\alpha,i)\in A_\ve^\prime)\}\]
is large too. Since $A_\ve\neq A_\ve^\prime$ implies that $A_\ve$ is not large
we get (by \ref{cl1}) that 
\[A\setminus\bigcup\{A_\ve: A_\ve\neq A_\ve^\prime\ \&\ \ve<\theta\}\neq
\emptyset.\] 
So take $(\alpha,i)\in A\setminus\bigcup\{A_\ve: A_\ve\neq A_\ve^\prime\ \&\
\ve<\theta\}$. We find $\ve<\theta$ such that $x_{\alpha,i}\in D_\ve$ (so
$(\alpha,i)\in A_\ve$). Then $A_\ve=A_\ve^\prime$ and we get $F_\ve(\alpha,i)
\in A_\ve$. Hence $x_{\alpha,i}, x_{F_\ve(\alpha,i)}\in D_\ve$ and
$x_{\alpha,i}\wedge x_{F_\ve(\alpha,i)}=0$. 
\begin{claim}
\label{cl3}
Let $D$ be a non-principal ultrafilter on $\kappa$.\\
Then $d(\prod\limits_{\gamma<\kappa}B_{\bar{F}_\gamma}/D)\leq\theta$.
\end{claim}
\noindent Proof of the claim: Fix functions $h:\theta^+\times\theta^+
\longrightarrow\theta$ and $h^*:\theta^+\times\theta\longrightarrow\theta^+$
such that for $i\in (\theta,\theta^+)$, $\zeta\in\theta$:
\begin{quotation}
\noindent $j_1<j_2<i\ \ \Rightarrow\ \ h(i,j_1)\neq h(i,j_2)$, $\qquad
h^*(i,\zeta)<i$\quad and

\noindent $j<i\ \ \Rightarrow\ \ h^*(i,h(i,j))=j$.
\end{quotation}
For $\gamma<\kappa$ let $Z_\gamma\subseteq B_{\bar{F}_\gamma}$ be the set of
all meets $x_{a_0}\wedge\ldots\wedge x_{a_{n-1}}\wedge(-x_{b_0})\wedge\ldots
\wedge (-x_{b_{m-1}})$ such that:  
\begin{quotation}
\noindent $a_0,\ldots,a_{n-1},b_0,\ldots,b_{m-1}\in\lambda\times\theta^+$ are
with no repetition, 

\noindent for all $k,l<n$, $r<m$ and all $\ve<\theta$
\[F^\gamma_\ve(a_k)\neq a_l\quad \&\quad F^\gamma_\ve(b_r)\neq a_l\quad \&
\quad F^\gamma_\ve(a_l)\neq b_r\] 
\end{quotation}
Clearly $Z_\gamma$ is dense in $B_{\bar{F}_\gamma}$ and $Z\stackrel{\rm def}
{=}\prod\limits_{\gamma<\kappa} Z_\gamma/D$ is dense in $\prod\limits_{\gamma<
\kappa}B_{\bar{F}_\gamma}/D$. For $e\in\prod\limits_{\gamma<\kappa} Z_\gamma$
and $\gamma<\kappa$ let:  
\begin{itemize}
\item $e(\gamma)=\bigwedge\limits_{l<n(e,\gamma)}x_{a(e,l,\gamma)}\wedge
\bigwedge\limits_{l<m(e,\gamma)}-x_{b(e,l,\gamma)}$,
\item $a(e,l,\gamma)=(\alpha(e,l,\gamma),i(e,l,\gamma))$,
\item $b(e,l,\gamma)=(\beta(e,l,\gamma),j(e,l,\gamma))$,
\item $\base^\gamma(e)=\{a(e,l,\gamma): l<n(e,\gamma)\}\cup\{b(e,l,\gamma):
l<m(e,\gamma)\}$,
\item $\base(e)=\bigcup\limits_{\gamma<\kappa}\base^\gamma(e)$,
\item $u_0(e)=\{i<\theta^+:(\exists\alpha<\lambda)((\alpha,i)\in\base(e))\}$,
\item $u_1(e)$ be the (topological) closure of $u_0(e)$,
\item $u_2(e)$ be the closure of $u_1(e)$ under the functions $h$, $h^*$,
\item $\zeta(e)$ be the first $\ve<\theta$ such that
\[(\forall\gamma<\kappa)(\forall(\alpha,i)\in\base(e)\cap\dom(F^\gamma_\ve))
(\sup (u_1(e)\cap i)< F^\gamma_{\ve,2}(\alpha,i)).\]
\end{itemize}
[Note that $\|\base(e)\|\leq\kappa<\cf(\theta)=\theta$, so $\|u_0(e)\|,
\|u_1(e)\|, \|u_2(e)\|\leq \kappa$; looking at the definition of $\zeta(e)$
remember that $(\alpha,i)\in\dom(F^\gamma_\ve)$ implies $\cf(i)=\theta$.]

Next for each $\gamma<\kappa$, $(\alpha,i)\in S^+_\gamma$, $\zeta_0<\theta$
choose $\ve^\gamma_{\zeta_0}(\alpha,i)<\theta$ such that the sequence
$\langle\eta_{F^\gamma_\zeta(\alpha,i)}(\ve^\gamma_{\zeta_0}(\alpha,i)):
\zeta\leq\zeta_0\rangle$ is strictly increasing (it is enough to take
$\ve^\gamma_{\zeta_0}(\alpha,i)$ sufficiently large -- apply condition (7)
for $\bar{F}_\gamma$ remembering $\zeta_0<\theta$). Further, for $(\alpha,i),
(\beta,j)\in\lambda\times\theta^+$, $\gamma<\kappa$ and $\zeta_0<\theta$ such
that $(\beta,j)\notin\{F^\gamma_\ve(\alpha,i):\ve\leq\zeta_0\}$ choose
$\ve^\gamma_{\zeta_0}((\alpha,i),(\beta,j))<\theta$ such that for every
$\zeta\leq\zeta_0$: 
\begin{quotation}
\noindent either $\eta_{F^\gamma_\zeta(\alpha,i)}(\ve^\gamma_{\zeta_0}
((\alpha,i),(\beta,j))) \neq \eta_{\beta,j}(\ve^\gamma_{\zeta_0}
((\alpha,i),(\beta,j))) $

\noindent or $\eta_{F^\gamma_\zeta(\alpha,i)}(\ve^\gamma_{\zeta_0}
(\alpha,i)) \neq \eta_{\beta,j}(\ve^\gamma_{\zeta_0}(\alpha,i))$
\end{quotation}
(this is possible as the second condition may fail for at most one
$\zeta\leq\zeta_0$: the sequence $\langle\eta_{F^\gamma_\zeta(\alpha,i)}
(\ve^\gamma_\zeta(\alpha,i)):\zeta\leq\zeta_0\rangle$ is strictly increasing).
Next, for each $e\in\prod\limits_{\gamma<\kappa} Z_\gamma$ and $\gamma<\kappa$
choose a finite set ${\cal X}_\gamma(e)\subseteq\theta$ such that:
\begin{description}
\item[8. ] if $a,b\in\base^\gamma(e)$ are distinct then $\eta_a\rest {\cal
X}_\gamma(e)\neq \eta_b\rest {\cal X}_\gamma(e)$,
\item[9. ] if $a\in\base^\gamma(e)\cap S^+_\gamma$ then $\ve^\gamma_{\zeta(e)}
(a)\in {\cal X}_\gamma(e)$,
\item[10. ] if $a,b\in\base^\gamma(e)$ then $\ve^\gamma_{\zeta(e)}(a,b)\in
{\cal X}_\gamma(e)$ (if defined),
\end{description}
(remember that $\base^\gamma(e)$ is finite). Finally we define a function $H$
on $\prod\limits_{\gamma<\kappa} Z_\gamma$ such that for $e\in
\prod\limits_{\gamma<\kappa} Z_\gamma$ the value $H(e)$ is the sequence
consisting of the following objects: 
\begin{description}
\item[11. ] $\langle n(e,\gamma):\gamma<\kappa\rangle$,
\item[12. ] $\langle m(e,\gamma):\gamma<\kappa\rangle$,
\item[13. ] $\zeta(e)$,
\item[14. ] $\langle{\cal X}_\gamma(e): \gamma<\kappa\rangle$,
\item[15. ] $\langle (\gamma,l,\eta_{a(e,l,\gamma)}\rest {\cal
X}_\gamma(e)): \gamma<\kappa, l<n(e,\gamma)\rangle$,
\item[16. ] $\langle (\gamma,l,\eta_{b(e,l,\gamma)}\rest {\cal
X}_\gamma(e)): \gamma<\kappa, l<m(e,\gamma)\rangle$,
\item[17. ] $u_2(e)\cap\theta$,
\item[18. ] $\{(\otp(i\cap u_2(e)), \otp(i\cap u_1(e))): i\in u_2(e)\}$.
\end{description}
Since $\theta^\kappa=\theta$ we easily check that $\|\rng(H)\|\leq\theta$. For
$\Upsilon\in\rng(H)$ let 
\[Z_\Upsilon=\{e\in \prod_{\gamma<\kappa} Z_\gamma: H(e)=\Upsilon\}\quad
\mbox{ and }\quad Z^*_\Upsilon =\{e/D: e\in Z_\Upsilon\}\subseteq Z.\]
The claim will be proved if we show that
\medskip

\centerline{for each $\Upsilon\in\rng(H)$ the set $Z^*_\Upsilon$ is centered.}
\medskip

\noindent First note that if $e,e'\in Z_\Upsilon$ then $u_2(e)\cap u_2(e')$ is
an initial segment of both $u_2(e)$ and $u_2(e')$. Why? Suppose that $j<i\in
u_2(e)\cap u_2(e')$, $j\in u_2(e)$. If $j<\theta$ then $j\in u_2(e')$ since
$u_2(e)\cap\theta=u_2(e')\cap\theta$. Suppose that $\theta\leq j<\theta^+$.
Then $h(i,j)\in u_2(e)\cap\theta =u_2(e')\cap\theta$ and so $j=h^*(i,h(i,j))
\in u_2(e')$. This shows that $u_2(e)\cap u_2(e')$ is an initial segment of
$u_2(e)$. Similarly for $u_2(e')$. 

\noindent Applying to this fact condition (18) we may conclude that $u_1(e)\cap
u_1(e')$ is an initial segment of both $u_1(e)$ and $u_1(e')$ for $e,e'\in
Z_\Upsilon$. [Why? Assume not. Let $i<\theta^+$ be the first such that there
is $j\in u_1(e)\cap u_1(e')$ above $i$ but 
\[i\in (u_1(e)\setminus u_1(e'))\cup (u_1(e')\setminus u_1(e)).\]
By symmetry we may assume that $i\in u_1(e)\setminus u_1(e')$. Let $i^*$ be
the first element of $u_2(e)$ above $i$. Then necessarily $i^*\leq j$ (as
$j\in u_1(e)\cap u_1(e')\subseteq u_2(e)\cap u_2(e')$) and hence $i^*\in
u_2(e')$ (and $i^*$ is the first element of $u_2(e')$ above $i$). By the 
choice of $i,i^*$ we have
\[i,i^*\in u_2(e)\cap u_2(e'),\quad i\cap u_1(e)=i\cap u_1(e'),\quad\mbox{ and
} i^*\cap u_2(e)=i^*\cap u_2(e').\]
But now we may apply condition (18) to conclude that 
\[(\otp(i^*\cap u_2(e)),\otp(i^*\cap u_1(e)))=(\otp(i^*\cap u_2(e')),\otp(i^*
\cap u_1(e')))\]
and therefore
\[\otp(i^*\cap u_1(e'))=\otp(i^*\cap u_1(e))=\otp(i\cap u_1(e))+1= \otp(i\cap
u_1(e'))+1.\]
As there is no point of $u_1(e')$ in the interval $[i,i^*)$ (remember $u_1(e')
\subseteq u_2(e')$) we get a contradiction.]
\smallskip

\noindent For $e\in Z_\Upsilon$ we have: $n(e,\gamma)=n(\gamma)$, $m(e,\gamma)
=m(\gamma)$, $\zeta(e)=\zeta^*$, ${\cal X}_\gamma(e)={\cal X}_\gamma$. Let
$e_0,\ldots,e_{k-1}\in Z_\Upsilon$. We are going to show that 
\[\prod\limits_{\gamma<\kappa}B_{\bar{F}_\gamma}/D\models e_0/D\wedge\ldots
\wedge e_{k-1}/D\neq 0\]
and for this we have to prove that
\[I_{e_0,\ldots,e_{k-1}}\stackrel{\rm def}{=}\{\gamma<\kappa:
B_{\bar{F}_\gamma}\models\bigwedge\limits_{j<k}[\bigwedge\limits_{l<n(\gamma)}
x_{a(e_j,l,\gamma)}\wedge\bigwedge\limits_{l<m(\gamma)}-x_{b(e_j,l,\gamma)}]
\neq 0\}\in D.\]
First let us ask what can be the reasons for
\[\bigwedge\limits_{j<k}[\bigwedge\limits_{l<n(\gamma)}x_{a(e_j,l,\gamma)}
\wedge\bigwedge\limits_{l<m(\gamma)}-x_{b(e_j,l,\gamma)}]= 0.\]
There are essentially two cases here: either $x_a\wedge (-x_a)$ appears on the
left-hand side of the above equality or $x_a\wedge x_{F^\gamma_\zeta(a)}$ (for
some $\zeta<\theta$) appears there. Suppose that the first case happens. Then
we have distinct $j_1,j_2<k$ such that $a(e_{j_1},l_1,\gamma)=b(e_{j_2},l_2,
\gamma)$ for some $l_1,l_2$.  By (15) (and the definition of $Z_\Upsilon$) we
have $\eta_{a(e_{j_1},l_1,\gamma)}\rest {\cal X}_\gamma=\eta_{a(e_{j_2},l_1,
\gamma)}\rest {\cal X}_\gamma$ and by (8) we have $\eta_{a(e_{j_2},l_1,
\gamma)}\rest {\cal X}_\gamma\neq\eta_{b(e_{j_2},l_2,\gamma)}\rest {\cal
X}_\gamma$. Consequently $\eta_{a(e_{j_1},l_1,\gamma)}\rest {\cal X}_\gamma
\neq\eta_{b(e_{j_2},l_2,\gamma)}\rest {\cal X}_\gamma$, a contradiction. 
Consider now the second case and suppose additionally that $\zeta\leq\zeta^*$.
Thus we assume that for some $\zeta\leq\zeta^*$, for some distinct $j_1,j_2<k$
and some $l_1,l_2<n(\gamma)$ we have $F^\gamma_\zeta(a(e_{j_1},l_1,\gamma))=
a(e_{j_2},l_2,\gamma)$. Then by (15) we get $\eta_{a(e_{j_2},l_2,\gamma)}\rest
{\cal X}_\gamma=\eta_{a(e_{j_1},l_2,\gamma)}\rest{\cal X}_\gamma$. As $\zeta
\leq\zeta^*$ we have that [by the choice of $\ve^\gamma_{\zeta^*}(a(e_{j_1},
l_1,\gamma),a(e_{j_1},l_2,\gamma))$, $\ve^\gamma_{\zeta^*}(a(e_{j_1},l_1,
\gamma))$ -- note that $F^\gamma_\ve(a(e_{j_1}, l_1,\gamma))\neq a(e_{j_1},
l_2,\gamma)$ for all $\ve\leq\zeta^*$]: 
\begin{quotation}
\noindent either \ \ \ \ $\eta_{F^\gamma_\zeta(a(e_{j_1},l_1,\gamma))}(
\ve^\gamma_{\zeta^*}(a(e_{j_1},l_1,\gamma),a(e_{j_1},l_2,\gamma)))\neq$\\
$\eta_{a(e_{j_1},l_2,\gamma)}(\ve^\gamma_{\zeta^*}(a(e_{j_1},l_1,\gamma),
a(e_{j_1},l_2,\gamma)))$ 

\noindent or\ \ \ \ \ \ $\eta_{F^\gamma_\zeta(a(e_{j_1},l_1,\gamma))}(
\ve^\gamma_{\zeta^*}(a(e_{j_1},l_1,\gamma)))\neq
\eta_{a(e_{j_1},l_2,\gamma)}(\ve^\gamma_{\zeta^*}(a(e_{j_1},l_1,\gamma)))$
\end{quotation}
and $\ve^\gamma_{\zeta^*}(a(e_{j_1},l_1,\gamma),a(e_{j_1},l_2,\gamma)),
\ve^\gamma_{\zeta^*}(a(e_{j_1},l_1,\gamma))\in {\cal X}_\gamma$ (by (9), (10);
note that in the definition of $\ve^\gamma_\zeta(a,b)$ we allowed $a=b$ so no
problem appears if $l_1=l_2$). Hence $\eta_{F^\gamma_\zeta(a(e_{j_1},l_1,
\gamma))}\rest {\cal X}_\gamma\neq\eta_{a(e_{j_1},l_2,\gamma)}\rest {\cal
X}_\gamma$ and thus $\eta_{F^\gamma_\zeta(a(e_{j_1},l_1,\gamma))}\rest
{\cal X}_\gamma\neq\eta_{a(e_{j_2},l_2,\gamma)}\rest {\cal X}_\gamma$, a
contradiction. Consequently, the considered equality may hold only if $x_a
\wedge x_{F^\gamma_\zeta(a)}$ appears there for some $\zeta>\zeta^*$. 

Asume now that $I_{e_0,\ldots,e_{k-1}}\notin D$. From the above considerations
we know that for each $\gamma\in\kappa\setminus I_{e_0,\ldots,e_{k-1}}$ we find
distinct $j_1(\gamma),j_2(\gamma)<k$ and $l_1(\gamma),l_2(\gamma)<n(\gamma)$
and $\zeta_\gamma\in(\zeta^*,\theta)$ such that 
\begin{description}
\item[$(**)$] \qquad $F^\gamma_{\zeta_\gamma}(a(e_{j_1(\gamma)},l_1(\gamma),
\gamma))= a(e_{j_2(\gamma)},l_2(\gamma),\gamma)$
\end{description}
(note that $(**)$ implies $a(e_{j_1(\gamma)},l_1(\gamma),\gamma)\in{\rm
dom}(F^\gamma_{\zeta_\gamma})$, $i(e_{j_1(\gamma)},l_1(\gamma),\gamma)\in
S_\gamma$). We have assumed that $\kappa\setminus I_{e_0,\dots,e_{k-1}}\in D$
so we find $j_1,j_2<k$ such that
\[J\stackrel{\rm def}{=}\{\gamma\in\kappa\setminus I_{e_0,\dots,e_{k-1}}:
j_1(\gamma)=j_1, j_2(\gamma)=j_2\}\in D.\]
As we have remarked after $(**)$, $i(e_{j_1},l_1(\gamma),\gamma)\in S_\gamma$
(for $\gamma\in J$) and consequently there are no repetitions in the sequence
$\langle i(e_{j_1},l_1(\gamma),\gamma):\gamma\in J\rangle$ (and $J$ is
infinite). Choose $\gamma_n\in J$ (for $n\in\omega$) such that the sequence
$\langle i(e_{j_1},l_1(\gamma_n),\gamma_n): n\in\omega\rangle$ is strictly
increasing (so $i(e_{j_1},l_1(\gamma_n),\gamma_n)\in u_1(e_{j_1})\cap
i(e_{j_1},l_1(\gamma_{n+1}),\gamma_{n+1})$) and let $i=\lim\limits_n
i(e_{j_1},l_1(\gamma_n),\gamma_n)$. By the definition of $\zeta(e),\zeta^*$
and the fact that $\zeta_\gamma>\zeta^*$ for all $\gamma\in J$ (and by (5)) we
have that for $\gamma\in J$  
\[\begin{array}{ll}
i(e_{j_2},l_2(\gamma),\gamma)=&F^\gamma_{\zeta_\gamma,2}(a(e_{j_1},l_1(
\gamma),\gamma))\in\\
\ &i(e_{j_1},l_1(\gamma),\gamma)\setminus\sup(u_1(e_{j_1})\cap i(e_{j_1},l_1
(\gamma),\gamma)).
\end{array}\]
Applying this for $\gamma_{n+1}$ we conclude
\[i(e_{j_1},l_1(\gamma_n),\gamma_n)< i(e_{j_2},l_2(\gamma_{n+1}),
\gamma_{n+1})< i(e_{j_1},l_1(\gamma_{n+1}),\gamma_{n+1})\]
and $i=\lim\limits_n i(e_{j_2},l_2(\gamma_n),\gamma_n)$. Since $u_1(e_{j_1}),
u_1(e_{j_2})$ are closed we conclude that $i\in u_1(e_{j_1})\cap
u_1(e_{j_2})$. From the remark we did after the definition of $Z_\Upsilon$ we
know that the last set is an initial segment of both $u_1(e_{j_1})$ and
$u_1(e_{j_2})$. But this gives a contradiction: $i(e_{j_2},l_2(\gamma_{n+1}),
\gamma_{n+1})\in u_1(e_{j_2})\setminus u_1(e_{j_1})$ and it is below $i\in
u_1(e_{j_1})\cap u_1(e_{j_2})$.  The claim is proved.
\medskip

Similarly as in claim~\ref{cl3} (but much easier) one can prove that really
$d(B_{\bar{F}_\gamma})=\theta^+$. \QED
\medskip

We want to finish this section with posing two questions motivated by
\ref{ex2.4} and \ref{th2.5}:

\begin{problem}
Are the following theories consistent?
\begin{enumerate}
\item ZFC + there is a cardinal $\kappa$ such that for each Boolean
algebra $B$,
\[d_n(B)\leq\kappa\ \ \Rightarrow\ \ d_{n+1}(B)<2^\kappa.\]
\item ZFC + there is a cardinal $\theta$ such that
$\theta^{\aleph_0}=\theta$ and for each Boolean algebra $B$ and a
non-principal ultrafilter $D$ on $\omega$
\[d(B)\leq\theta\ \ \Rightarrow\ \ d(B^\omega/D)< 2^\theta.\]
\end{enumerate}
\end{problem}

\section{Hereditary cofinality and spread}

\subsection{The invariants}
The hereditary cofinality of a Boolean algebra $B$ is the cardinal
\[\hcof(B)=\min\{\kappa: (\forall X\subseteq B)(\exists C\subseteq X)(\|C\|
\leq\kappa\ \&\ C \mbox{ is cofinal in } X)\}.\]
It can be represented as a def.u.w.o.car.~invariant if we use the following
description of it (see \cite{M}): 
\begin{description}
\item[$(\otimes_{\hcof})$] $\hcof(B)=\sup\{\|X\|: X\subseteq B\ \& \ (X,<_B)$
is well-founded $\}$. 
\end{description}
Let the theory $T_{\hcof}$ introduce predicates $P_0,P_1$ on which it says
that: 
\begin{itemize}
\item $P_1$ is a well ordering of $P_0$,
\item $(\forall x_0,x_1\in P_0)(x_0<x_1\ \Rightarrow\ P_1(x_0,x_1))$
\end{itemize}
(in the above $<$ stands for the respective relation of the Boolean algebra).
Clearly $T_{\hcof}$ determines a def.u.w.o.car.~invariant and 
\[\Inv_{T_{\hcof}}(B)=\{\|X\|:X\subseteq B\ \&\ (X,<) \mbox{ is
well-founded}\}.\] 

The spread $s(B)$ of a Boolean algebra $B$ is
\[s(B)=\sup\{\|S\|: S\subseteq\Ult B\ \&\ S\mbox{ is discrete in the
relative topology}\}.\]
It can be easily described as a def.f.o.car.~invariant: the suitable theory
$T_s$ introduces predicates $P_0,P_1$ and it says that for each $x\in P_0$
the set $\{y: P_1(x,y)\}$ is an ultrafilter and the ultrafilters form a
discrete set (in the relative topology). Sometimes it is useful to remember
the following characterization of $s(B)$ (see \cite{M}):
\begin{description}
\item[$(\otimes_s)$] $s(B)=\sup\{\|X\|: X\subseteq B\mbox{ is
ideal-independent}\}.$ 
\end{description}
Using this characterization we can write $s(B)=s_\omega(B)$, where
\begin{definition}
\begin{enumerate}
\item $\phi_n^s$ is the formula saying that no member of $P_0$ can be
covered by union of $n+1$ other elements of $P_0$.
\item For $0<n\leq\omega$ let $T^n_s=\{\phi_k^s:k<n\}$.
\item For a Boolean algebra $B$ and $0<n\leq\omega$:
$s^{(+)}_n(B)=\inv_{T^n_s}^{(+)}(B)$ (so $s_n$ are def.f.o.car.~invariants).
\end{enumerate}
\end{definition}

The hereditary density of a Boolean algebra $B$ is the cardinal
\[ \hd (B)=\sup\{dS: S\subseteq\Ult B\}\]
where $dS$ is the (topological) density of the space $S$. The following
characterization of $\hd (B)$ is important for our purposes (see \cite{M}):
\begin{description}
\item[$(\otimes_{\hd })$] $\hd (B)=\sup\{\|\kappa\|:$ there is a strictly
decreasing sequence of ideals (in $B$) of the length $\kappa\ \}$.
\end{description}
We should remark here that on both sides of the equality we have $\sup$ but
the attainment does not have to be the same. If the $\sup$ of the left hand
side ($\hd (B)$) is obtained then so is the $\sup$ of the other side. If the
right hand side $\sup$ is obtained AND $\hd (B)$ is regular then the $\sup$ of
$\hd (B)$ is realized. An open problem is what can happen if $\hd (B)$ is
singular. 

The hereditary Lindel\"of degree of a Boolean algebra $B$ is
\[\hL(B)=\sup\{LS:S\subseteq\Ult B\},\]
where for a topological space $S$, $LS$ is the minimal $\kappa$ such that
every open cover of $S$ has a subcover of size $\leq\kappa$. The following
characterization of $\hL(B)$ is crucial for us (see \cite{M}): 
\begin{description}
\item[$(\otimes_{\hL})$] $\hL(B)=\sup\{\|\kappa\|:$ there is a strictly
increasing sequence of ideals (in $B$) of the length $\kappa\ \}$.
\end{description}
Note: we may have here differences in the attainment, like in the case of
$\hd$. 

\begin{definition}
\begin{description}
\item[1.\ ] Let the formula $\psi$ say that $P_1$ is a well ordering of
$P_0$ (denoted by $<_1$).
\item[2.\ ] For $n<\omega$ let $\phi^{\hd}_n$, $\phi^{\hL}_n$ be the
following formulas:
\item[$\phi^{\hd}_n\equiv$] $\psi\ \&\ (\forall x_0,\ldots,x_{n+1}\in P_0)(
x_0<_1\ldots<_1 x_{n+1}\ \Rightarrow\ x_0\not\leq x_1\vee\ldots\vee x_{n+1})$
\item[$\phi^{\hL}_n\equiv$] $\psi\ \&\ (\forall x_0,\ldots,x_{n+1}\in P_0)(
x_{n+1}<_1\ldots<_1 x_0\ \Rightarrow\ x_0\not\leq x_1\vee\ldots\vee x_{n+1})$.
\item[3.\ ] For $0<n\leq\omega$ we let $T^n_{\hd}=\{\phi^{\hd}_k:k<n\}$,
$T^n_{\hL}=\{\phi^{\hL}_k: k<n\}$. 
\item[4.\ ] For a Boolean algebra $B$ and $0<n\leq\omega$:
\[\hd^{(+)}_n(B)=\inv^{(+)}_{T^n_{\hd}}(B),\quad\quad \hL^{(+)}_n(B)=
\inv^{(+)}_{T^n_{\hL}}(B).\]
\end{description}
\end{definition}
So $\hd_n$, $\hL_n$ are def.u.w.o.car.~invariants and $\hd_\omega=\hd$,
$\hL_\omega=\hL$ (the sets $\Inv_{T^\omega_{\hd}}(B)$, $\Inv_{T^\omega_{\hL}}
(B)$ agree with the sets on the right-hand sides of $(\otimes_{\hd})$,
$(\otimes_{\hL})$, respectively). 

\subsection{Constructions from strong $\lambda$-systems.}
One of our tools for constructing examples of Boolean algebras is an
object taken from the pcf theory.

\begin{definition}
\label{systems}
\begin{enumerate}
\item {\em A weak $\lambda$-system} (for a regular cardinal $\lambda$) is a
sequence ${\cal S} =\langle \delta,\bar{\lambda},\bar{f}\rangle$ such that
\begin{description}
\item[a)\ ] $\delta$ is a limit ordinal, $\|\delta\|<\lambda$,
\item[b)\ ] $\bar{\lambda}=\langle\lambda_i: i<\delta\rangle$ is a strictly
increasing sequence of regular cardinals,
\item[c)\ ] $\bar{f}=\langle f_\alpha:\alpha<\lambda\rangle\subseteq
\prod\limits_{i<\delta}\lambda_i$ is a sequence of pairwise distinct functions,
\item[d)\ ] for every $i<\delta$, $\|\{f_\alpha\rest i:\alpha<\lambda\}\|
\leq\sup\limits_{i<\delta}\lambda_i$.
\end{description}
\item {\em A $\lambda$-system} is a sequence ${\cal S} =\langle \delta,
\bar{\lambda},\bar{f}, J\rangle$ such that ${\cal S}_0 =\langle \delta,
\bar{\lambda},\bar{f}\rangle$ is a weak $\lambda$-system and
\begin{description}
\item[e)\ ] $J$ is an ideal on $\delta$ extending the ideal $J^{\rm
bd}_\delta$ of bounded subsets of $\delta$,
\item[f)\ ] $\bar{f}$ is a
$<_J$-increasing sequence cofinal in $\prod\limits_{i<\delta}(\lambda_i,<)/J$,
\item[g)\ ] for every $i<\delta$, $\|\{f_\alpha\rest i:\alpha<\lambda\}\|
<\lambda_i$.
\end{description}
In this situation we say that the system ${\cal S}$ {\em extends}
the weak system ${\cal S}_0$.
\item ${\cal S}=\langle\delta,\bar{\lambda},\bar{f},J,(A_\zeta:
\zeta<\kappa)\rangle$ is {\em a strong $\lambda$-system for $\kappa$} if
$\langle\delta,\bar{\lambda},\bar{f},J\rangle$ is a $\lambda$-system and
\begin{description}
\item[h)\ ] $\cf(\delta)\leq\kappa$, $\sup\limits_{i<\delta}\lambda_i
\leq 2^\kappa$,
\item[i)\ ] $A_\zeta\subseteq\delta$, $A_\zeta\notin J$ (for
$\zeta<\kappa$) are pairwise disjoint.
\end{description}
\end{enumerate}
\end{definition}

In ZFC, there is a class of cardinals $\lambda$ for which there are (weak,
strong) $\lambda$-systems. We can even demand that, for (weak)
$\lambda$-systems, $\lambda$ is the successor of a cardinal $\lambda_0$
satisfying $\lambda^\omega_0=\lambda_0$ (what is relevant for ultraproducts,
see below). More precisely:

\begin{fact}
\label{existence}
\begin{enumerate}
\item If $\mu^{<\kappa}<\mu^\kappa=\lambda$ then there is a weak
$\lambda$-system ${\cal S}=\langle\delta,\bar{\lambda},\bar{f}\rangle$
such that $\sup\limits_{i<\delta}\lambda_i\leq\mu$, $\delta=\kappa$.
\item If $\kappa=\cf(\kappa)$ and
\begin{description}
\item[$(*)$] $\kappa>\aleph_0$, $\mu=\mu^{<\kappa}<\lambda$,
$\cf(\lambda)\leq\mu^\kappa$
\end{description}
or even
\begin{description}
\item[$(*)^-$] $\cf(\mu)=\kappa$,
\[(\forall\theta)(\exists\mu_\theta<\mu)(\forall\chi)(\mu_\theta<\chi<\mu\
\&\ \cf(\chi)=\theta\ \Rightarrow {\rm pp}_\theta(\chi)<\mu)\]
\end{description}
then there is a $\lambda$-system ${\cal S}=\langle\delta,\bar{\lambda},
\bar{f},J\rangle$ such that $\mu=\sup\limits_{i<\delta}\lambda_i$,
$\delta=\kappa$ (see \cite{Sh 371}).
\item If $\kappa=\aleph_0$, $\cf(\mu)=\aleph_0<\mu$ and
\begin{quotation}
\noindent either $\lambda^*={\rm cov}(\mu,\mu,\aleph_1,2)$

\noindent or $\lambda^*=\lambda^{\aleph_0}\ \&\ (\forall\chi<\mu)(
\chi^{\aleph_0}<\mu)$
\end{quotation}
then for many regular $\lambda\in(\mu,\lambda^*)$ there are
$\lambda$-systems ($\lambda=\mu_0^+$ really) (see \cite{Sh 430}.

\item There is a class of cardinals $\lambda$ for
which there are strong $\lambda$-systems (for some infinite $\kappa$),
even if we additionally demand that $\lambda$ is a successor cardinal (see
\cite{Sh 400}, \cite{Sh 410} or the proof of 4.4 of \cite{Sh 462}). \QED
\end{enumerate}
\end{fact}

\begin{theorem}
\label{hcof}    
Assume that there exists a strong $\lambda$-system for $\kappa$, $\lambda$ a
regular cardinal. Let $\theta$ be an infinite cardinal $\leq\kappa$.
Then there are Boolean algebras $B_\ve$ (for $\ve<\theta$) such that
$\inv_{T_{\hcof}}^+(B_\ve)\leq\lambda$ and for any ultrafilter $D$ on
$\theta$ containing all co-bounded sets we have
$s^+_\omega(\prod\limits_{\ve<\theta}B_\ve/D)>\lambda$.
\end{theorem}

\noindent PROOF: The algebras $B_\ve$'s are modifications of the algebra
constructed in Lemma 4.2 of \cite{Sh 462}. Let $\langle\delta,\bar{\lambda},
\bar{f}, J, (A_\zeta:\zeta<\kappa)\rangle$ be a strong $\lambda$-system for
$\kappa$. For distinct $\alpha,\beta<\lambda$ let $\rho(\alpha,\beta)=\min
\{i<\delta: f_\alpha(i)\neq f_\beta(i)\}$. 

Take a decreasing sequence $\langle w_\ve:\ve<\theta\rangle$ of subsets of
$\kappa$ such that $\|w_\ve\|=\kappa$ and $\bigcap\limits_{\ve<\theta}w_\ve=
\emptyset$.

Fix $\ve<\theta$.

\noindent For $i<\delta$ choose a family $\{F_{i,\zeta}:\zeta<\kappa\}$ of
subsets of $\{f_\alpha\rest i: \alpha<\lambda\}$ such that if $X_1,X_2\in
[\{f_\alpha\rest i: \alpha<\lambda\}]^{<\omega}$ then for some $\zeta<\kappa$
we have $X_1= F_{i,\zeta}\cap (X_1\cup X_2)$ (possible as $2^\kappa\geq
\lambda_i$). Next take a sequence $\langle (j_i,\zeta_i):i<\delta\rangle$ such
that $j_i\leq i$, $\zeta_i<\kappa$ and the set 
\[\{j<\delta: (\forall\zeta<\kappa)(\exists\xi\in w_\ve)(A_\xi\subseteq_J
\{i<\delta: j_i=j\ \&\ \zeta_i=\zeta\})\}\]
is unbounded in $\delta$ (possible as $\cf(\delta)\leq\kappa$, $\|w_\ve\|=
\kappa$). 

Now we define a partial order $\prec_\ve$ on $\lambda$:
\begin{quotation}
\noindent $\alpha\prec_\ve\beta$\ \ \ if and only if\ \ \
$i=\rho(\alpha,\beta)\in\bigcup\limits_{\xi\in w_\ve}A_\xi$ and  
\[f_\alpha\rest j_i\in F_{j_i,\zeta_i}\ \iff\
f_\alpha(i)<f_\beta(i).\]
\end{quotation}
The algebra $B_\ve$ is the Boolean algebra generated by the partial order
$\prec_\ve$. It is the algebra of subsets of $\lambda$ generated by sets
$Z_\alpha=\{\beta<\lambda: \beta\prec_\ve\alpha\}\cup\{\alpha\}$ (for
$\alpha<\lambda$).  
\begin{claim}
\label{cl8}
\begin{description}
\item[a)] If $\rho(\alpha,\beta)<\rho(\beta,\gamma)$, $\alpha,\beta,\gamma<
\lambda$ then $\beta\prec_\ve\alpha\iff\gamma\prec_\ve\alpha$ and $\alpha
\prec_\ve\beta\iff\alpha\prec_\ve\gamma$. 
\item[b)] If $\tau(x_0,\ldots,x_{n-1})$ is a Boolean term, $\alpha^l_k<
\lambda$ for $k<n$ are pairwise distinct ($l<2$), $i<\delta$ is such that
$\rho(\alpha^0_k,\alpha^1_k)\geq i$ (for $k<n$) but $\rho(\alpha^l_k,
\alpha^l_{k'})<i$ (for $l<2$, $k<k'<n$)\\
then denoting $X_l=\tau(Z_{\alpha^l_0},\ldots,Z_{\alpha^l_{n-1}})$ we have 
\end{description}
\begin{enumerate}
\item $X_0\cap\{\alpha<\lambda: (\forall k<n)(f_\alpha\rest i\neq
f_{\alpha^0_k}\rest i)\}= X_1\cap\{\alpha<\lambda: (\forall k<n)(f_\alpha\rest
i\neq f_{\alpha^0_k}\rest i)\}$ 
\item for each $k<n$,\\
either $X_l\supseteq\{\alpha<\lambda:f_\alpha\rest i=f_{\alpha^0_k}\rest i\}$
for $l<2$\\ 
or $X_l\cap\{\alpha<\lambda: f_\alpha\rest i=f_{\alpha^0_k}\rest i\}=
Z_{\alpha^l_k}\cap \{\alpha<\lambda:f_\alpha\rest i=f_{\alpha^0_k}\rest i\}$
for $l<2$\\ 
or $X_l\cap\{\alpha<\lambda: f_\alpha\rest i=f_{\alpha^0_k}\rest i\}=\{\alpha<
\lambda: f_\alpha\rest i=f_{\alpha^0_k}\rest i\}\setminus Z_{\alpha^l_k}$ for
$l<2$\\ 
or $X_l\cap\{\alpha<\lambda:f_\alpha\rest i=f_{\alpha^0_k}\rest i\}=\emptyset$
for $l<2$. 
\end{enumerate}
\end{claim}
\begin{claim}
\label{cl9}
Suppose that $\langle a_\alpha:\alpha<\lambda\rangle$ are distinct members
of $B_\ve$. Then there exist $\alpha<\beta<\lambda$ such that $a_\alpha\geq
a_\beta$.
\end{claim}
Proof of the claim:\ \ \ First we may assume that for some integers $n<m<
\omega$, a Boolean term $\tau(x_0,\ldots,x_{n-1},\ldots, x_{m-1})$, ordinals
$\alpha_n,\ldots,\alpha_{m-1}<\lambda$, an ordinal $i^*<\delta$ and a function
$\bar{\alpha}:\lambda\times n\longrightarrow\lambda\setminus \{\alpha_n,
\ldots,\alpha_{m-1}\}$ for all $\beta<\lambda$ we have  
\begin{quotation}
\noindent $a_\beta=\tau(Z_{\bar{\alpha}(\beta,0)},\ldots,
Z_{\bar{\alpha}(\beta,n-1)}, Z_{\alpha_n},\ldots,Z_{\alpha_{m-1}})$, and 

\noindent if $\beta'<\lambda$, $k,k'<n$, $(\beta,k)\neq (\beta',k')$ then
$\bar{\alpha}(\beta,k)\neq\bar{\alpha}(\beta',k')$, and 

\noindent $\{f_{\bar{\alpha}(\beta,k)}\rest i^*,
f_{\alpha_{k'}}\rest i^*: k<n, n\leq k'<m\}$ are pairwise distinct.
\end{quotation}
As we may enlarge $i^*$ we may additionally assume that 
\[(\forall\zeta<\kappa)(\exists\xi\in w_\ve)(A_\xi\subseteq_J \{i<\delta:
j_i=i^*\ \&\ \zeta_i=\zeta\}).\]
Furthermore, we may assume that $f_{\bar{\alpha}(\beta,k)}\rest i^* =
f_{\bar{\alpha}(0,k)}\rest i^*$ for all $\beta<\lambda$, $k<n$
(remember that $\|\{f_\alpha\rest i^*: \alpha<\lambda\}\|<\lambda$). 
Let $B$ be the set of all $i<\delta$ such that
\[(\forall \zeta<\lambda_i)(\exists^\lambda\beta<\lambda)(\forall k<n)(\zeta<
f_{\bar{\alpha}(\beta,k)}(i)).\]
Then the set $B$ is in the dual filter $J^c$ of $J$ (if not clear see Claim
3.1.1 of \cite{Sh 462}). Now apply the choice of $F_{i^*,\zeta}$'s to find
$\zeta<\kappa$ such that for $k<n$: 
\medskip

\noindent if $a_0\cap\{\alpha<\lambda: f_{\bar{\alpha}(0,k)}\rest i^*
= f_\alpha\rest i^*\}= Z_{\bar{\alpha}(0,k)}\cap\{\alpha<\lambda:
f_{\bar{\alpha}(0,k)}\rest i^* = f_\alpha\rest i^*\}$\\
then $f_{\bar{\alpha}(0,k)}\rest i^*\notin F_{i^*,\zeta}$ and

\noindent if $a_0\cap\{\alpha<\lambda: f_{\bar{\alpha}(0,k)}\rest i^*
= f_\alpha\rest i^*\}= \{\alpha<\lambda:
f_{\bar{\alpha}(0,k)}\rest i^* = f_\alpha\rest i^*\} \setminus
Z_{\bar{\alpha}(0,k)}$\\
then $f_{\bar{\alpha}(0,k)}\rest i^*\in F_{i^*,\zeta}$.
\medskip

\noindent Note that by claim~\ref{cl8} we can replace 0 in the above by any
$\beta<\lambda$. Take $\xi\in w_\ve$ such that $A_\xi\subseteq_J \{i<\delta:
j_i=i^*\ \&\ \zeta_i=\zeta\}$ and choose $i\in A_\xi\cap B$ such that $j_i=
i^*$, $\zeta_i=\zeta$. Since $\|\{f_\alpha\rest i:\alpha<\lambda\}\|<
\lambda_i$ and $i\in B$ we find $\beta_0<\beta_1<\lambda$ such that 
\[(\forall k<n)(\rho(\bar{\alpha}(\beta_0,k),\bar{\alpha}(\beta_1,k))=i)\
\&\ \max\limits_{k<n}f_{\bar{\alpha}(\beta_0,k)}(i)<\min\limits_{k<n}
f_{\bar{\alpha}(\beta_1,k)}(i).\] 
Now by the choice of $\zeta$, claim~\ref{cl8} and the property of
$\beta_0,\beta_1$ we get $a_{\beta_1}\subseteq a_{\beta_0}$. 
\begin{claim}
\label{cl10}
$\inv^+_{T_{\hcof}}(B_\ve)\leq\lambda$
\end{claim}
Proof of the claim:\ \ \ Directly from claim~\ref{cl9} noting that $\lambda
\longrightarrow (\lambda,\omega)^2$.
\begin{claim}
\label{cl11}
Suppose that $\alpha_0,\ldots,\alpha_n<\lambda$ are
pairwise $\prec_\ve$-incomparable. Then $Z_{\alpha_0}\not\subseteq
Z_{\alpha_1}\cup\ldots\cup Z_{\alpha_n}$.
\end{claim}

Suppose now that $D$ is an ultrafilter on $\theta$ containing all co-bounded
sets. 
\begin{claim}
\label{cl12}
$s^+_\omega(\prod\limits_{\ve<\theta}B_\ve/D)>\lambda$.
\end{claim}
Proof of the claim:\ \ \ We need to find an ideal--independent subset of
$\prod\limits_{\ve<\theta}B_\ve/D$ of size $\lambda$. But this is easy: 
for $\alpha<\lambda$ let $x_\alpha\in\prod\limits_{\ve<\theta}B_\ve/D$ be such
that $x_\alpha(\ve)=\{\beta<\lambda: \beta\preceq_\ve\alpha\}$. The set
$\{x_\alpha:\alpha<\lambda\}$ is ideal--independent since if $\alpha_0,\ldots,
\alpha_n<\lambda$ are distinct and $\ve$ is such that $\alpha_0,\ldots,
\alpha_n$ are pairwise $\prec_\ve$-incomparable then 
\[B_\ve\models x_{\alpha_0}(\ve)\not\leq x_{\alpha_1}(\ve)\vee\ldots\vee
x_{\alpha_n}(\ve)\]
(by claim~\ref{cl11}). Now note that if $\ve_0<\theta$ is such that
\[\bigcup_{\xi\in w_{\ve_0}}A_\xi\cap\{\rho(\alpha_l,\alpha_m):l<m<n\} =
\emptyset\]
then for all $\ve\geq\ve_0$ we have that $\alpha_0,\ldots,\alpha_n$ are
pairwise $\prec_\ve$-incomparable. Now by \L o\'s theorem we conclude
\[\prod\limits_{\ve<\theta}B_\ve/D\models x_{\alpha_0}\not\leq x_{\alpha_1}\vee
\ldots\vee x_{\alpha_n}. \QED\]
\medskip

\noindent{\bf Remark:}\ \ \ 1.\ \  For $\lambda$ such that there exists
a strong $\lambda$-system and $\lambda$ is a successor (and for the respective
$\theta,\kappa$'s) we have algebras $B_\ve$ (for $\ve<\theta$) such that
$\inv_T(B_\ve)<\lambda$ and for respective ultrafilters $D$ on $\kappa$ 
$\inv_T(\prod\limits_{\ve<\theta}B_\ve/D)\geq\lambda$, where $T$ is one of
the following:
\begin{quotation}
$T_{\hcof}$, $T_s^\omega$, $T_{\hd}^\omega$, $T_{\hL}^\omega$ or $T_{{\rm
inc}}$.
\end{quotation}
2.\ \ We do not know if (in ZFC) we can demand $\lambda=\lambda_0^+$ and
$\lambda_0^\omega=\lambda_0$; consistently yes.

\begin{theorem}
\label{newspread}
Assume that there exists a strong $\lambda$--system for $\kappa$,
$0<n<\omega$. Then there is a Boolean algebra $B$ such that $s^+_n(B)=\|B\|^+
=\lambda^+$ (so $\hd^+_n(B)=\hL^+_n(B)=\lambda^+$) but $s^+_\omega(B),\hd^+(B),
\hL^+(B)\leq\lambda$.
\end{theorem}

\noindent PROOF:\qquad The construction is slightly similar to the one of
\ref{hcof}. 

Let $\langle\delta,\bar{\lambda},\bar{f},J,(A_\zeta:\zeta<\kappa)
\rangle$ be a strong $\lambda$--system for $\kappa$, $\rho(\alpha,\beta)=\min
\{i<\delta: f_\alpha(i)\neq f_\beta(i)\}$ (for distinct $\alpha,\beta<
\lambda$) and let $F_{i,\zeta}\subseteq\{f_\alpha\rest i: \alpha<\lambda\}$
(for $i<\delta$, $\zeta<\kappa$) be such that if $X_1,X_2\in [\{f_\alpha\rest
i: \alpha<\lambda\}]^{<\omega}$ then there is $\zeta<\kappa$ with $X_1= F_{i,
\zeta}\cap (X_1\cup X_2)$. Like before, fix a sequence $\langle (j_i,\zeta_i):
i<\delta\rangle$ such that $j_i\leq i$, $\zeta_i<\kappa$ and the set 
\[\{j<\delta: (\forall\zeta<\kappa)(\exists\xi<\kappa)(A_\xi\subseteq_J
\{i<\delta: j_i=j\ \&\ \zeta_i=\zeta\})\}\]
is unbounded in $\delta$. 

Let $B$ be the Boolean algebra generated freely by $\{x_\alpha:\alpha<\lambda
\}$ except that 
\begin{description}
\item[$(\alpha)$] if $\alpha_0,\ldots,\alpha_{n+2}<\lambda$, $i<\delta$, $f_0
\rest i=\ldots=f_{\alpha_{n+2}}\rest i$, $f_{\alpha_0}(i)<f_{\alpha_1}(i)<
\ldots<f_{\alpha_{n+2}}(i)$ and $f_{\alpha_0}\rest j_i\in F_{j_i,\zeta_i}$
then $x_{\alpha_0}\leq x_{\alpha_1}\vee\ldots\vee x_{\alpha_{n+2}}$\qquad and
\item[$(\beta)$] if $\alpha_0,\ldots,\alpha_{n+2}<\lambda$, $i<\delta$, $f_0
\rest i=\ldots=f_{\alpha_{n+2}}\rest i$, $f_{\alpha_0}(i)<f_{\alpha_1}(i)<
\ldots<f_{\alpha_{n+2}}(i)$ and $f_{\alpha_0}\rest j_i\notin F_{j_i,\zeta_i}$
then $x_{\alpha_1}\wedge\ldots\wedge x_{\alpha_{n+2}}\leq x_{\alpha_0}$.
\end{description}

\begin{claim}
\label{cl20}
If $\alpha_0,\ldots,\alpha_n<\lambda$ are pairwise distinct then
\[B\models x_{\alpha_0}\not\leq x_{\alpha_1}\vee\ldots\vee x_{\alpha_n}.\]
Consequently $s^+_n(B)=\hd^+_n(B)=\hL^+_n(B)=\|B\|^+=\lambda^+$.
\end{claim}
Proof of the claim:\ \ \ Let $h:\lambda\longrightarrow 2$ be such that
$h(\alpha_0)=1$ and for $\alpha\in\lambda\setminus\{\alpha_0\}$ 
\[h(\alpha)=\left\{\begin{array}{ll}
0 &\mbox{ if }f_\alpha\rest j_i\notin F_{j_i,\zeta_i}\mbox{ or}\\
\ &\ \ \      f_\alpha\rest j_i\in F_{j_i,\zeta_i}\mbox{ and } f_\alpha\rest
(i+1)\in\{f_{\alpha_l}\rest (i+1):l=1,\ldots,n\},\\
1 &\mbox{ otherwise},
		   \end{array}\right.\]
where $i=\rho(\alpha_0,\alpha)$. We are going to show that the function $h$
preserves the inequalities imposed on $B$ in $(\alpha)$, $(\beta)$ above. To
deal with $(\alpha)$ suppose that $\beta_0,\ldots,\beta_{n+2}<\lambda$,
$f_{\beta_0}\rest i=\ldots=f_{\beta_{n+2}}\rest i$, $f_{\beta_0}(i)<\ldots<
f_{\beta_{n+2}}(i)$ and $f_{\beta_0}\rest j_i\in F_{j_i,\zeta_i}$. If
$h(\beta_0)=0$ then there are no problems, so let us assume that $h(\beta_0)=
1$. Since $f_{\beta_k}\rest(i+1)$ (for $k=1,\ldots,n+2$) are pairwise distinct
we find $k_0\in\{1,\ldots,n+2\}$ such that
\[f_{\beta_{k_0}}\rest (i+1)\notin\{f_{\alpha_l}\rest (i+1): l\leq n\}.\]
It is easy to check that then $h(\beta_{k_0})=1$, so we are done. Suppose now
that $\beta_0,\ldots,\beta_{n+2}<\lambda$, $f_{\beta_0}\rest i=\ldots=
f_{\beta_{n+2}}\rest i$, $f_{\beta_0}(i)<\ldots<f_{\beta_{n+2}}(i)$ but 
$f_{\beta_0}\rest j_i\notin F_{j_i,\zeta_i}$ and suppose $h(\beta_0)=0$
(otherwise trivial). If $\rho(\alpha_0,\beta_0)<i$ then clearly $h(\beta_k)=
h(\beta_0)=0$ for all $k\leq n+2$. If $\rho(\alpha_0,\beta_0)\geq i$ then for
some $k_0\in \{1,\ldots,n+2\}$ we have $\alpha_0\neq \beta_{k_0}$,
$\rho(\alpha_0,\beta_{k_0})=i$ and easily $h(\beta_{k_0})=0$.

\begin{claim}
\label{cl21}
$\hd^+(B),\hL^+(B)\leq\lambda$.
\end{claim}
Proof of the claim:\ \ \ Suppose that $\langle a_\beta:\beta<\lambda\rangle
\subseteq B$. After the standard cleaning we may assume that for some Boolean
term $\tau$, integers $m_0<m<\omega$, a function $\bar{\alpha}:\lambda\times
m\longrightarrow\lambda$, and an ordinal $i_0<\delta$ for all $\beta<\lambda$
we have: 
\begin{description}
\item[$(*)_1$] $a_\beta=\tau(x_{\bar{\alpha}(\beta,0)},\ldots,x_{\bar{\alpha}
(\beta,m-1)})$,
\item[$(*)_2$] $\langle f_{\bar{\alpha}(\beta,l)}\rest i_0: l<m\rangle$ are
pairwise distinct and $f_{\bar{\alpha}(\beta,l)}\rest i_0=f_{\bar{\alpha}(0,
l)}\rest i_0$ (for $l<m$) ,
\item[$(*)_3$] $\{\langle\bar{\alpha}(\beta,0),\ldots,\bar{\alpha}(\beta,m-1)
\rangle:\beta<\lambda\}$ forms a $\Delta$-system of sequences with the root
$\{0,\ldots,m_0-1\}$.
\end{description}
Moreover, as we are dealing with $\hd,\hL$, we may assume that the term $\tau$
is of the form
\[\tau(x_0,\ldots,x_{m-1})=\bigwedge_{l<m} x^{t(l)}_l,\]
where $t:m\longrightarrow 2$. Let $\zeta<\kappa$ be such that for each $l<m$ 
\[f_{\bar{\alpha}(0,l)}\rest i_0\in F_{i_0,\zeta}\iff t(l)=0.\]
Take $i_1>i_0$ such that $j_{i_1}=i_0$, $\zeta_{i_1}=\zeta$ and
\[(\forall\zeta<\lambda_{i_1})(\exists^\lambda \beta<\lambda)(\forall
l\in [m_0,m))(\zeta<f_{\bar{\alpha}(\beta,l)}(i_1))\]
(like in the proof of claim~\ref{cl9}). Now, as $\|\{f_\alpha\rest i_1:\alpha
<\lambda\}\|<\lambda_{i_1}$, we may choose distinct $\beta_0,\ldots,
\beta_{m\cdot(n+2)}<\lambda$ such that for $k\leq m\cdot (n+2)$, $l<m$
\[f_{\bar{\alpha}(\beta_0,l)}\rest i_1 =f_{\bar{\alpha}(\beta_k,l)}\rest
i_1 \stackrel{\rm def}{=} \nu_l\]
and for each $l\in [m_0,m)$
\[f_{\bar{\alpha}(\beta_0,l)}(i_1)<f_{\bar{\alpha}(\beta_1,l)}(i_1)<
\ldots<f_{\bar{\alpha}(\beta_{m\cdot(n+2)},l)}(i_1).\]
Note that we can demand any order between $\beta_0,\ldots,\beta_{m\cdot(n+2)}$
we wish what allows us to deal with both $\hd$ and $\hL$. We are going to
show that $a_{\beta_0}\leq\bigvee\limits_{k=1}^{m\cdot (n+2)} a_{\beta_k}$.  
Suppose that $l\in [m_0,m)$, $1\leq k_1<k_2<\ldots<k_{n+2}\leq m\cdot
(n+2)$. If $t(l)=0$ then, by the choice of $\zeta$ and $i_1$ we may apply
clause $(\alpha)$ of the definition of $B$ and conclude that
\[x_{\bar{\alpha}(\beta_0,l)}\leq x_{\bar{\alpha}(\beta_{k_1},l)}\vee\ldots
x_{\bar{\alpha}(\beta_{k_{n+2}},l)}.\]
Similarly, if $t(l)=1$ then 
\[x_{\bar{\alpha}(\beta_{k_1},l)}\wedge\ldots x_{\bar{\alpha}(\beta_{k_{n+2}},
l)}\leq x_{\bar{\alpha}(\beta_0,l)}.\]
Hence, for any distinct $k_1,\ldots,k_{n+2}\in \{1,\ldots,m\cdot (n+2)\}$ and
$l<m$ we have
\[x_{\bar{\alpha}(\beta_0,l)}^{t(l)}\leq x_{\bar{\alpha}(\beta_{k_1},
l)}^{t(l)}\vee\ldots x_{\bar{\alpha}(\beta_{k_{n+2}},l)}^{t(l)},\]
and therefore
\[\bigwedge_{l<m} x_{\bar{\alpha}(\beta_0,l)}^{t(l)}\leq \bigvee_{k=1}^{m\cdot
(n+2)}\bigwedge_{l<m} x_{\bar{\alpha}(\beta_k,l)}^{t(l)}.\qquad\QED\]

\noindent{\bf Remark:}\quad Theorem \ref{newspread} is applicable to
ultraproducts, of course, but we do not know if we can demand (in ZFC) that
$\lambda=\lambda^+_0$, $\lambda_0^\omega=\lambda$.\\
ZFC constructions (using $\lambda$--systems) parallel to \ref{newspread}
will be presented in a forthcoming paper \cite{Sh 620}. Some related
consistency results will be contained in \cite{RoSh 599}.

\begin{problem}
For each $0<n<\omega$ find (in ZFC) a Boolean algebra $B$ such that
$s_n(B)>s_{n+1}(B)$. Similarly for $\hL$, $\hd$. 
\end{problem}
\medskip

\subsection{Forcing an example}
\begin{theorem}
\label{diamond}
Assume that $\aleph_0\leq\kappa<\mu<\lambda=\mu^+=2^\mu$. Then there is a
forcing notion ${\Bbb P}$ which is $(<\lambda)$--complete of size $\lambda^+$
and satisfies the $\lambda^+$--cc (so it preserves cardinalities, cofinalities
and cardinal arithmetic) and such that in ${\bf V}^{\Bbb P}$:
\begin{quotation}
\noindent there exist Boolean algebras $B_\xi$ (for $\xi<\kappa$) such that
$\hd(B_\xi)$, $\hL(B_\xi)\leq\lambda$ (so $s(B_\xi)\leq\lambda$)
but for each ultrafilter $D$ on $\kappa$ containing co-bounded subsets of
$\kappa$ we have $\ind(\prod\limits_{\xi<\kappa}B_\xi/D)\geq\lambda^+$ (so
$\lambda^+\leq\hd (\prod\limits_{\xi<\kappa}B_\xi/D), \hL(\prod\limits_{\xi<
\kappa}B_\xi/D),s(\prod\limits_{\xi<\kappa}B_\xi/D)$).
\end{quotation}
\end{theorem}

\noindent PROOF: By Theorem 2.5(3) of \cite{Sh 462} there is a suitable
forcing notion ${\Bbb P}$ such that in ${\bf V}^{\Bbb P}$:
\smallskip

\noindent there is a sequence $\langle \eta_i: i<\lambda^+\rangle\subseteq
{}^\lambda\lambda$ with no repetition and functions $c,d$ such that:
\begin{description}
\item[(a)] $c:{}^{\lambda{>}}\lambda\longrightarrow\lambda$, 
\item[(b)] the domain $\dom(d)$ of the function $d$ consists of all pairs
$(\bar{x},h)$ such that $h:\zeta\longrightarrow\lambda\times\lambda\times
\lambda$ for some $\zeta<\mu$, and $\bar{x}:\mu\longrightarrow{}^\alpha
\lambda$ is one-to-one, $\alpha<\lambda$,
\item[(c)] for $(\bar{x},h)\in\dom(d)$, $d(\bar{x},h)$ is a function from
\[\{\bar{a}\in{}^\mu(\lambda^+): \bar{a}\mbox{ is increasing and }(\forall
i<\mu)(x_i\vartriangleleft\eta_{a_i})\}\]
to $\lambda$ such that $d(\bar{x},h)(\bar{a})=d(\bar{x},h)(\bar{b})$ implies
$\sup\bar{a}\neq\sup\bar{b}$ and denoting $t_i=\eta_{a_i}\wedge\eta_{b_i}$
for some $i^*<\mu$ we have:
\begin{description}
\item[$(\alpha)$] ${\rm level}(t_i)={\rm level}(t_{i^*})$ for $i>i^*$,
\item[$(\beta)$]  $(\forall\ve<\mu)(\exists^\mu i<\mu)(c(t_i)=\ve)$,
\item[$(\gamma)$] for $\mu$ ordinals $i<\mu$ divisible by $\zeta$ we have 
\begin{description}
\item[(i)] either there are $\xi_0<\xi_1<\lambda$ such that
\[\hspace{-0.5cm}(\forall\ve<\zeta)(\zeta\cdot\xi_0\leq\eta_{\bar{b}_{i+\ve}}
({\rm level}(t_{i+\ve}))<\zeta\cdot\xi_1\leq\eta_{\bar{a}_{i+\ve}}({\rm level}
(t_{i+\ve}))),\]
and
\[\hspace{-1.8cm}h=\langle(c(t_{i+\ve}),\eta_{\bar{b}_{i+\ve}}({\rm level}(
t_{i+\ve}))-\zeta\cdot\xi_0,\eta_{\bar{a}_{i+\ve}}({\rm level}(t_{i+\ve}))-
\zeta\cdot\xi_1):\ve<\zeta\rangle,\]
\item[(ii)] or a symmetrical condition interchanging $\bar{a}$ and $\bar{b}$.
\end{description}
\end{description}
\end{description}
From now on we are working in the universe ${\bf V}^{\Bbb P}$ using the objects
listed above. 

For distinct $i,j<\lambda^+$ let $\rho(i,j)=\min\{\xi<\lambda:
\eta_i(\xi)\neq\eta_j(\xi)\}$. For $\ve_0,\ve_1<\kappa$ we put
\[\begin{array}{lr}
R^\kappa_{\ve_0,\ve_1}=\{(i,j)\in\lambda^+\times\lambda^+\!:& 
i\neq j\ \mbox{ and }\ \eta_i(\rho(i,j))=\ve_0\ \mod\kappa\ \mbox{ and }\ \\
\ &\eta_j(\rho(i,j))=\ve_1\ \mod\kappa\},
  \end{array}\]
and now we define Boolean algebras $B_{\kappa,\bar{\ve}}$ for
$\bar{\ve}=\langle\ve_0,\ve_1,\ve_2,\ve_3\rangle$, $\ve_0<\ve_1<\ve_2<\ve_3
<\kappa$. $B_{\kappa,\bar{\ve}}$ is the Boolean algebra freely generated by
$\{x_i: i<\lambda^+\}$ except: 
\begin{quotation}
\noindent if $(i,j)\in R^\kappa_{\ve_0,\ve_1}$ then $x_i\leq x_j$,

\noindent if $(i,j)\in R^\kappa_{\ve_2,\ve_3}$ then $x_j\leq x_i$.
\end{quotation}
\begin{claim}
\label{cl4}
If $i,j<\lambda^+$, $(i,j)\notin
R^\kappa_{\ve_0,\ve_1}$, $(j,i)\notin R^\kappa_{\ve_2,\ve_3}$ then
$B_{\kappa,\bar{\ve}}\models x_i\not\leq x_j$. In particular, if
$i<j<\lambda^+$ then $B_{\kappa,\bar{\ve}}\models x_i\neq x_j$ and
$\|B_{\kappa,\bar{\ve}}\|=\lambda^+$.
\end{claim}
\noindent Proof of the claim:\ \ \ Fix $i<\lambda^+$. A function
$f:\{x_j:j<\lambda^+\}\longrightarrow {\cal P}(2)$ (where ${\cal P}(2)$ is
the Boolean algebra of subsets of $\{0,1\}$) is defined by $f(x_i)=\{0\}$
and for $j\in\lambda^+\setminus\{i\}$:
\begin{quotation}
\noindent if $(i,j)\in R^\kappa_{\ve_0,\ve_1}$ or $(j,i)\in
R^\kappa_{\ve_2,\ve_3}$ then $f(x_j)=\{0,1\}$

\noindent if $(i,j)\in R^\kappa_{\ve_2,\ve_3}$ or $(j,i)\in
R^\kappa_{\ve_0,\ve_1}$ then $f(x_j)=\emptyset$ and

\noindent otherwise $f(x_j)=\{1\}$.
\end{quotation}
We are going to show that $f$ respects all the inequalities we put on $x_j$'s
in $B_{\kappa,\bar{\ve}}$. So suppose that $(j_1,j_2)\in R^\kappa_{\ve_0,
\ve_1}$. If $f(x_{j_1})=\emptyset$ then there are no problems, so assume that
both $(i,j_1)\notin R^\kappa_{\ve_2,\ve_3}$ and $(j_1,i)\notin R^\kappa_{
\ve_0,\ve_1}$. Similarly, we may assume that $f(x_{j_2})\neq\{0,1\}$,
i.e.~that both $(i,j_2)\notin R^\kappa_{\ve_0,\ve_1}$ and $(j_2,i)\notin
R^\kappa_{\ve_2,\ve_3}$. Note that these two assumptions imply $j_1\neq i\neq
j_2$. Now we consider three cases:
\begin{description}
\item[--] if $\rho(j_1,j_2)>\rho(i,j_1)=\rho(i,j_2)$ then $f(x_{j_1})=f(
x_{j_2})$,
\item[--] if $\rho(j_1,j_2)=\rho(i,j_1)=\rho(i,j_2)$ then $f(x_{j_1})=\{1\}=
f(x_{j_2})$ (remember $(j_1,j_2)\in R^\kappa_{\ve_0,\ve_1}$, $(i,j_2),(j_1,i)
\notin R^\kappa_{\ve_0,\ve_1}$), 
\item[--] if $\rho(j_1,j_2)<\max\{\rho(j_1,i),\rho(j_2,i)\}$ then either
$\rho(j_1,j_2)=\rho(i,j_2)<\rho(i,j_1)$ and $(i,j_2)\in R^\kappa_{\ve_0,
\ve_1}$ (what is excluded already) or $\rho(j_1,j_2)=\rho(i,j_1<\rho(i,j_2$
and $(j_1,i)\in R^\kappa_{\ve_0,\ve_1}$ (what is against our assumption too).
\end{description}
This shows that $f(x_{j_1})\leq f(x_{j_2})$ whenever $(j_1,j_2)\in R^\kappa_{
\ve_0,\ve_1}$. Similarly one shows that $(j_1,j_2)\in R^\kappa_{\ve_2,\ve_3}$
implies $f(x_{j_2})\leq f(x_{j_1})$. Consequently the function $f$ respects
all the inequalities in the definition of $B_{\kappa,\bar{\ve}}$. Hence it
extends to a homomorphism $\bar{f}:B_{\kappa,\bar{\ve}}\longrightarrow{\cal
P}(2)$. But for each $j<\lambda^+$
\[((i,j)\notin R^\kappa_{\ve_0,\ve_1}\ \&\ (j,i)\notin
R^\kappa_{\ve_2,\ve_3})\ \Rightarrow\ (f(x_j)\in\{\emptyset,\{1\}\}\ \&\
f(x_i)=\{0\}).\] 
\begin{claim}
\label{cl5}
Suppose $\bar{i}:\lambda^+\times n\longrightarrow\lambda^+$, $\bar{t}:n
\longrightarrow 2$, $n<\omega$ are such that $(\forall\alpha<\lambda^+)(
\forall l_1<l_2<n)(\bar{i}(\alpha,l_1)<\bar{i}(\alpha,l_2))$. Then
\begin{description}
\item[($\oplus_1$)] \quad$(\exists\alpha<\beta<\lambda^+)(B_{\kappa,\bar{\ve}}
\models\bigwedge\limits_{l<n}(x_{\bar{i}(\alpha,l)})^{\bar{t}(l)}\leq
\bigwedge\limits_{l<n}(x_{\bar{i}(\beta,l)})^{\bar{t}(l)})$,
\item[($\oplus_2$)] \quad$(\exists\alpha<\beta<\lambda^+)(B_{\kappa,\bar{\ve}}
\models\bigwedge\limits_{l<n}(x_{\bar{i}(\alpha,l)})^{\bar{t}(l)}\geq
\bigwedge\limits_{l<n}(x_{\bar{i}(\beta,l)})^{\bar{t}(l)})$.
\end{description}
\end{claim}
\noindent Proof of the claim:\ \ \ To prove $(\oplus_1),(\oplus_2)$ it is
enough to show the following:
\begin{description}
\item[($\oplus_1^*$)] \quad$(\exists\alpha<\beta<\lambda^+)(\forall l<n)(
B_{\kappa,\bar{\ve}}\models(x_{\bar{i}(\alpha,l)})^{\bar{t}(l)}\leq (
x_{\bar{i}(\beta,l)})^{\bar{t}(l)})$,
\item[($\oplus_2^*$)] \quad$(\exists\alpha<\beta<\lambda^+)(\forall l<n)(
B_{\kappa,\bar{\ve}}\models(x_{\bar{i}(\alpha,l)})^{\bar{t}(l)}\geq(x_{
\bar{i}(\beta,l)})^{\bar{t}(l)})$.
\end{description}
By the definition of $B_{\kappa,\bar{\ve}}$ for $(\oplus_1^*)$ it is enough to
have 
\begin{description}
\item[($\oplus_1^{**}$)] \quad there are $\alpha<\beta<\lambda^+$ such that 
\[l<n\ \&\ \bar{t}(l)=0\quad \Rightarrow\quad (\bar{i}(\alpha,l),\bar{i}(
\beta,l))\in R^\kappa_{\ve_0,\ve_1}\mbox{ or }\bar{i}(\alpha,l)=\bar{i}(
\beta,l),\] 
\[l<n\ \&\ \bar{t}(l)=1\quad \Rightarrow\quad (\bar{i}(\alpha,l),\bar{i}(
\beta,l))\in R^\kappa_{\ve_2,\ve_3}\mbox{ or }\bar{i}(\alpha,l)=\bar{i}(\beta,
l),\] 
\end{description}
and similarly for $(\oplus_2^*)$. 

We will show how to get $(\oplus_1^{**})$ from the properties of $\langle
\eta_i:i<\lambda^+\rangle$. For this we start with a cleaning procedure in
which we pass from the sequence $\langle\langle\bar{i}(\alpha,l):l<n\rangle:
\alpha<\lambda^+\rangle$ to its subsequence $\langle\langle\bar{i}(\alpha,l):
l<n\rangle:\alpha\in A\rangle$ for some $A\subseteq\lambda^+$ of size
$\lambda^+$ (so we will assume $A=\lambda^+$). First note that if $i$ repeats
$\lambda^+$ times in $\langle\bar{i}(\alpha,l):\alpha<\lambda^+,l<n\rangle$
then for some $l<n$ we have $\|\{\alpha:\bar{i}(\alpha,l)=i\}\|=\lambda^+$ and
we may assume that for all $\alpha<\lambda^+$, $\bar{i}(\alpha,l)=i$.
Consequently $(\oplus_1^{**})$ holds trivially for this $l$ (and every $\alpha
<\beta<\lambda^+$). Thus we may assume that each value appears at most
$\lambda$ times in $\langle\bar{i}(\alpha,l):\alpha<\lambda^+,l<n\rangle$ and
hence we may assume that the sets $\{\bar{i}(\alpha,l):l<n\}$ are disjoint for
$\alpha<\lambda^+$ (so there are no repetitions in $\langle \bar{i}(\alpha,l):
\alpha<\lambda^+,l<n\rangle$). Further we may assume that 
\[\alpha<\beta<\lambda^+\ \Rightarrow\ \bar{i}(\alpha,0)<\ldots
<\bar{i}(\alpha,n-1)<\bar{i}(\beta,0)<\ldots<\bar{i}(\beta,n-1).\] 
For $l<n$, $\alpha<\lambda^+$ let $a_{n\alpha+l}=\bar{i}(\alpha,l)$. We find
$\xi<\lambda$ such that for $\lambda^+$ ordinals $\beta<\lambda^+$ divisible
by $\mu$ the sequence $\langle\eta_{a_{\beta+\ve}}\rest\xi:\ve<\mu\rangle$ is
with no repetitions and does not depend on $\beta$ (for these $\beta$). Since
$\lambda=\mu^+=2^\mu$ there are $\xi<\lambda$ and a one-to-one sequence
$\bar{x}:\mu\longrightarrow{}^\xi\lambda$ such that the set
\[\begin{array}{lr}
B=\{\beta<\lambda^+:&\beta\mbox{ is divisible by $\mu$ and }\ \\
\ &(\forall\ve<\mu)(\eta_{a_{\beta+\ve}}\rest\xi=x_\ve)\}
  \end{array}\]
is of size $\lambda^+$. Let $h:\kappa\longrightarrow\lambda^3$ be such that
for $l<2n$: 
\[h(l)=\left\{
\begin{array}{ll}
(0,\ve_0,\ve_1) &\mbox{ if }\ \bar{t}(l)=0,\ l<n,\\
(0,\ve_1,\ve_0) &\mbox{ if }\ \bar{t}(l-n)=0,\ n\leq l<2n,\\
(0,\ve_2,\ve_3) &\mbox{ if }\ \bar{t}(l)=1,\ l<n,\\
(0,\ve_3,\ve_2) &\mbox{ if }\ \bar{t}(l-n)=1,\ n\leq l<2n.\\
\end{array}
\right.\]
Consider the function $d(\bar{x},h)$. There are distinct $\beta_0,\beta_1\in
B$ such that $d(\bar{x},h)(\langle a_{\beta_0+\ve}:\ve<\mu\rangle)=d(\bar{x},
h)(\langle a_{\beta_1+\ve}: \ve<\mu\rangle)$. This implies that we find
$\delta<\mu$ divisible by $\kappa$ such that (possibly interchanging
$\beta_0,\beta_1$):\\
there are $\xi_0<\xi_1<\lambda$ such that for some $\gamma<\lambda$ for every
$\ve<\kappa$ 

$\rho(a_{\beta_0+\delta+\ve},a_{\beta_1+\delta+\ve})=\gamma$, \quad and

$\kappa\cdot\xi_0\leq\eta_{a_{\beta_0+\delta+\ve}}(\gamma)<\kappa\cdot\xi_1
\leq\eta_{a_{\beta_1+\delta+\ve}}(\gamma)$, \quad and

$h=\langle (c(\eta_{a_{\beta_0+\delta+\ve}}\rest\gamma),\eta_{a_{\beta_0+
\delta+\ve}}(\gamma)-\kappa\cdot\xi_0,\eta_{a_{\beta_1+\delta+\ve}}-\kappa
\cdot\xi_1):\ve<\kappa\rangle$.

\noindent Suppose that $\beta_0<\beta_1$ and look at the values
$\eta_{a_{\beta_0+\delta+l}}(\gamma)$, $\eta_{a_{\beta_1+\delta+l}}(\gamma)$
for $l<n$. By the definition of $h$ we have that 
\begin{description}
\item[--] if $t(l)=0$ then $\eta_{a_{\beta_0+\delta+l}}(\gamma)=\ve_0\ \mod
\kappa$ and $\eta_{a_{\beta_1+\delta+l}}(\gamma)=\ve_1\ \mod\kappa$ (so
$(a_{\beta_0+\delta+l},a_{\beta_1+\delta+l})\in R^\kappa_{\ve_0,\ve_1}$),\quad
and 
\item[--] if $t(l)=1$ then $\eta_{a_{\beta_0+\delta+l}}(\gamma)=\ve_2\ \mod
\kappa$ and $\eta_{a_{\beta_1+\delta+l}}(\gamma)=\ve_3\ \mod\kappa$ (so
$(a_{\beta_0+\delta+l},a_{\beta_1+\delta+l})\in R^\kappa_{\ve_2,\ve_3}$).
\end{description}
Consequently $\beta_0+\delta<\beta_1+\delta<\lambda^+$ are as required in
$(\oplus^{**}_1)$. If $\beta_1<\beta_0$ then we look at the values
$\eta_{a_{\beta_0+\delta+n+l}}(\gamma),\eta_{a_{\beta_1+\delta+n+l}}(\gamma)$
(for $l<n$) and similarly we conclude that $\beta_1+\delta+n<\beta_0+\delta+n<
\lambda^+$ witness $(\oplus^{**}_1)$. 

Similarly one can get $\oplus_2^*$.
\begin{claim}
\label{cl6}
$\hd(B_{\kappa,\bar{\ve}})\leq\lambda$,
$\hL(B_{\kappa,\bar{\ve}})\leq\lambda$. 
\end{claim}
\noindent Proof of the claim: \ \ \ Suppose that $\hL(B_{\kappa,\bar{\ve}})\
\geq\lambda^+$ (or $\hd(B_{\kappa,\bar{\ve}})\geq\lambda^+$). Then there is a
sequence $\langle y_\alpha:\alpha<\lambda^+\rangle\subseteq B_{\kappa,
\bar{\ve}}$ such that for each $\alpha<\lambda^+$ the element $y_\alpha$ is
not in the ideal generated by $\{y_\beta:\beta<\alpha\}$ ($\{y_\beta:\beta> 
\alpha\}$, respectively). Moreover we can demand that each $y_\alpha$ is of
the form $\bigwedge\limits_{l<n(\alpha)}(x_{\bar{i}(\alpha,i)})^{\bar{t}(
\alpha,l)}$ with $\bar{i}(\alpha,l_1)<\bar{i}(\alpha,l_2)$ for $l_1<l_2<n(
\alpha)$. Next we may assume that $n(\alpha)=n$, $\bar{t}(\alpha,l)=\bar{t}
(l)$ for $\alpha<\lambda^+$, $l<n$ and apply $(\oplus_1)$ ($(\oplus_2)$,
respectively) of claim~\ref{cl5} to get a contradiction.
\medskip

Now, for $\xi<\kappa$ let $B_\xi=B_{\kappa,\langle 4\xi,4\xi+1,4\xi+2,4\xi+3
\rangle}$ and $B=\prod\limits_{\xi<\kappa}B_\xi/D$, where $D$ is an
ultrafilter on $\kappa$ such that no its member is bounded in $\kappa$.
\begin{claim}
\label{cl7}
$\ind(B)\geq\lambda^+$.
\end{claim}
\noindent Proof of the claim:\ \ \ Let $f_i\in\prod\limits_{\xi<\kappa}B_\xi$
(for $i<\lambda^+$) be the constant sequence $f_i(\xi)=x_i$. Suppose $i_0<i_1<
\ldots<i_{n-1}<\lambda^+$ and look at the set 
\[\begin{array}{rl}
{\cal X}=\{\xi<\kappa:(\exists j<4)(\exists m<k<n) &
(\eta_{i_m}(\rho(i_m,i_k))=4\xi+j\ \mod \kappa\\
\mbox{ or } & \eta_{i_k}(\rho(i_m,i_k))=4\xi+j\ \mod \kappa)\}.
\end{array}\]
Obviously, the set ${\cal X}$ is bounded in $\kappa$. By claim~\ref{cl4} (or
actually by a stronger version of it, but with a similar proof) we have that
for $\xi\in\kappa\setminus{\cal X}$
\[B_\xi\models\mbox{``}f_{i_0}(\xi),\ldots,f_{i_{n-1}}(\xi)\mbox{ are
independent elements''}.\]
Therefore, we conclude $B\models\mbox{``}f_{i_0},\ldots,f_{i_{n-1}}$ are
independent''. \QED 

\section{Independence number and tightness}
\subsection{Independence.} In this section we are interested in the
cardinal invariants related to the independence number.
\begin{definition}
\begin{enumerate}
\item $\phi_n^{\ind}$ is the formula which says that any non-trivial Boolean
combination of $n+1$ elements of $P_0$ is non-zero (i.e. $\phi_n^{\ind}$ says
that if $x_0,\ldots,x_n\in P_0$ are distinct then $\bigwedge\limits_{l\leq n}
x_l^{t(l)}\neq 0$ for each $t\in {}^{n+1}2$).
\item For $0<n\leq\omega$ let $T^n_{\ind}=\{\phi_k^{\ind}:k<n\}$.
\item For a Boolean algebra $B$, $0<n\leq\omega$ we define
$\ind_n(B)=\inv_{T^n_{\ind}}(B)$ and $\ind^+_n(B)=\inv_{T^n_{\ind}}^+(B)$.
We will denote $\ind_\omega^{(+)}$ by $\ind^{(+)}$ too.
\item A subset $X$ of a Boolean algebra $B$ is {\em $n$--independent} if and
only if any non-trivial Boolean combination of $n$ elements of $X$ is
non-zero.  
\end{enumerate}
\end{definition}

\noindent{\bf Remark:}\ \ \ 1. Note that the theory $T^{n+1}_{\ind}$ consists
of formulas $\phi^{\ind}_0,\ldots,\phi^{\ind}_n$ and thus it says that the set
$P_0$ is $n+1$--independent. Consequently for each $n<\omega$:
\[\ind^{(+)}_{n+1}(B)=\sup\{\|X\|^{(+)}: X\subseteq B\mbox{ is
$n+1$--independent}\}.\] 
2. It should be underlined here that the cardinal invariants $\ind_n$ ({\em
the $n$--independence number}) were first introduced and studied by Monk in
\cite{Mo4}.

\begin{proposition}
\label{pr3.1}
Suppose that $\lambda$ is an infinite cardinal, $n$ is an integer greater than
1. Then there is a Boolean algebra $B$ such that 
\[\ind_n(B)=\lambda=\|B\|\ \ \&\ \ \ \ind_{n+1}(B)=\aleph_0.\]
\end{proposition}

\noindent PROOF: Surprisingly the example we give depends on the parity of
$n$.  
\medskip

{\sc CASE 1:}\ \ \ \ \ $n=2k$, \ \ \ $k\geq 1$.
\medskip

\noindent Let ${\cal X}=\{x\in{}^\lambda 2: \|x^{-1}[\{1\}]\| \leq k\}$ and
for $\alpha<\lambda$ let $Z_\alpha=\{x\in {\cal X}: x(\alpha)=1\}$. Let 
$B^k_0(\lambda)$ be the Boolean algebra of subsets of ${\cal X}$ generated by
$\{Z_\alpha:\alpha<\lambda\}$.
\begin{claim}
\label{cl13}
$\ind_n(B^k_0(\lambda))=\lambda$.
\end{claim}
Proof of the claim:\ \ \ For $\alpha>0$ put $Y_\alpha=Z_0\vartriangle
Z_\alpha$ ($\vartriangle$ stands for the symmetric difference). We are going
to show that the set $\{Y_\alpha:0<\alpha<\lambda\}$ is $n$-independent. For
this suppose that $t\in {}^n 2$, $0<\alpha_0<\ldots<\alpha_{n-1}<\lambda$.
Choose $x\in {\cal X}$ such that
\begin{quotation}
\noindent if $\|t^{-1}[\{0\}]\|\leq k$ then $x(0)=0$, $x(\alpha_l)=1-t(l)$ for
$l<n$,

\noindent if $\|t^{-1}[\{0\}]\|>k$ then $x(0)=1$, $x(\alpha_l)=t(l)$ for
$l<n$. 
\end{quotation}
Then easily $x\in\bigcap\limits_{l<n}Y^{t(l)}_{\alpha_l}$.
\begin{claim}
\label{cl14}
$\ind_{n+1}(B^k_0(\lambda))=\aleph_0$.
\end{claim}
Proof of the claim:\ \ \ It should be clear that $\ind(B^k_0(\lambda))\geq
\aleph_0$, so what we have to show is $\ind_{n+1}(B^k_0(\lambda))<\aleph_1$.
Suppose that $\langle Y_\alpha:\alpha<\omega_1\rangle\subseteq B^k_0(
\lambda)$. We may assume that 
\begin{itemize}
\item $Y_\alpha=\tau(Z_{\bar{i}(\alpha,0)},\ldots,Z_{\bar{i}
(\alpha,m-1)})$, where $m<\omega$, $\tau$ is a Boolean term, $\bar{i}:
\omega_1\times m\longrightarrow\lambda$ is such that $\bar{i}(\alpha,0),
\ldots,\bar{i}(\alpha,m-1)$ are pairwise distinct,

\item $\{\langle\bar{i}(\alpha,0),\ldots,\bar{i}(\alpha,m-1)\rangle: \alpha
<\omega_1\}$ forms a $\Delta$-system of sequences with the root
$\{0,\ldots,m^*-1\}$ (for some $m^*\leq m$).
\end{itemize}
Further we may assume that $\tau(x_0,\ldots,x_{m-1})=\bigvee\limits_{t\in A}
\bigwedge\limits_{i<m} x^{t(i)}_i$ for some $A\subseteq {}^m2$. If $m=m^*$
(i.e.~all the $Y_\alpha$'s are the same) the sequence is not
$n+1$-independent. If $m^*=0$ (i.e.~the sets $\{\bar{i}(\alpha,l):l<m\}$ are
disjoint for $\alpha<\omega_1$) then either $Y_0\wedge\ldots\wedge Y_n=0$ or
$(-Y_0)\wedge\ldots\wedge(-Y_n)=0$ (e.g.~the first holds if $1\hat{\ }\ldots
\hat{\ }1=\bar{\bf 1}\notin A$ and otherwise the second equality is true). So
we may assume that $0<m^*<m$.
\smallskip

\noindent Suppose that $\bar{\bf 1}\in A$. We claim that then $(-Y_0)\wedge
\ldots\wedge(-Y_k)\wedge Y_{k+1}\wedge\ldots\wedge Y_{2k}=0$. If not then we
find $x\in\bigcap\limits_{k<j<2k+1}Y_j\setminus\bigcup\limits_{j<k+1}Y_j$. For
$j<2k+1$ let $t_j\in {}^m2$ be defined by $t_j(l)=1-x(\bar{i}(j,l))$. Thus
$t_j\in A$ for $k<j<2k+1$ and $t_j\notin A$ for $j<k+1$. As $\|x^{-1}[\{1\}]\|
\leq k$ for some $j_0\leq k$ we necessarily have $(\forall l\in [m^*,m))
(t_{j_0}(l)=1)$. Since ${\bf 1}\in A$ and $t_{j_0}\notin A$, necessarily for
some $l_0<m^*$ we have $t_{j_0}(l_0)=0$. Now look at $t_j$ for $j\in [k+1,2k]$.
Since $t_j\rest m^*= t_{j_0}\rest m^*$ and $t_{j_0}\notin A$ (and $t_j\in A$,
remember $k+1\leq j\leq 2k$) we have
\[(\forall j\in [k+1,2k])(\exists l_j\in [m^*,m))(t_j(l_j)=0).\]
This implies that $x(\bar{i}(j,l_j))=1$ for $j\in [k+1,2k]$ and together with
$x(\bar{i}(j_0,l_0))=1$ we get contradiction to $\|x^{-1}[\{1\}]\|\leq k$.

Suppose now that $\bar{\bf 1}\notin A$. Symmetrically to the previous case we
show that then $Y_0\wedge\ldots\wedge Y_k\wedge (-Y_{k+1})\wedge\ldots\wedge
(-Y_{2k})=0$. The claim is proved.
\medskip

{\sc CASE 2:}\ \ \ \ \ $n=2k+1$, \ \ \ $k\geq 1$.
\medskip

\noindent In this case we consider
\[{\cal X}'=\{x\in{}^\lambda 2:\|x^{-1}[\{1\}]\|\leq k\mbox{ or }
\|x^{-1}[\{0\}]\|\leq k\}\]
and the Boolean algebra $B^k_1(\lambda)$ of subsets of ${\cal X}'$ generated
by sets $Z_\alpha^\prime=\{x\in{\cal X}':x(\alpha)=1\}$. Then the sequence
$\langle Z_\alpha^\prime:\alpha<\lambda\rangle$ is $n$-independent
(witnessing $\ind_n(B^k_1(\lambda))=\lambda$). Similarly as in 
claim~\ref{cl14} one can show that $\ind_{n+1}(B^k_1(\lambda))
=\aleph_0$ (after the cleaning consider $(-Y_0)\wedge\ldots\wedge(-Y_k)\wedge
Y_{k+1}\ldots\wedge Y_{2k+1}$). \QED 
\medskip

\noindent{\bf Remark:}\ \ \ Note that 
\[\ind^{(+)}_{2k+1}(B^k_0(\lambda)\times B^k_0(\lambda))=\lambda^{(+)}\]
as witnessed by the set $\{(Z_\alpha,-Z_\alpha): \alpha<\lambda\}$. 

\begin{corollary}
\label{co5.3}
Suppose that $\lambda$ is an infinite cardinal. Then there are Boolean
algebras $B_n$ (for $n<\omega$) such that $\ind(B_n)=\aleph_0$ but for every 
non-principal ultrafilter $D$ on $\omega$, $\ind(\prod\limits_{n<\omega}B_n/D)
=\lambda^{\aleph_0}$. \QED 
\end{corollary}
\medskip

A detailed study of the reasons why we did have to consider two cases in
Proposition~\ref{pr3.1} leads to interesting observations concerning the
invariant $\ind_n$ and products of Boolean algebras. First note that

\begin{fact}
\label{indprofac}
For any Boolean algebras $B_i$ ($i<\lambda$) we have
\begin{enumerate}
\item $\ind_{2n}^+(B_0\times B_0)\leq\ind_n^+(B_0)\leq\ind_n^+(B_0\times
B_0)$,
\item $\ind^+_{\sum\limits_{i<k}n_i}(B_0\times\ldots\times B_{k-1})\leq
\sum\limits_{i<k}\ind^+_{n_i}(B_i)$,
\item $\ind^+(\prod^w_{i<\lambda}B_i)=\sup\limits_{i<\lambda}
\ind^+(B_i)$. \QED
\end{enumerate}
\end{fact}

However there is no immediate bound on $\ind_{n+1}(B\times B)$ in this 
context. One can easily show that the algebra $B^k_1(\lambda)$ from the proof
of \ref{pr3.1} (case 2) satisfies
\[\ind_{2k+2}(B^k_1(\lambda)\times B^k_1(\lambda))=\aleph_0.\]
So we get an example proving:
\begin{corollary}
If $\lambda$ is an infinite cardinal, $n$ is an odd integer $>2$ then there is
a Boolean algebra $B$ such that $\ind_n(B)=\lambda$ and $\ind_{n+1}(B\times
B)=\aleph_0$. \QED 
\end{corollary}
The oddity of $n$ in the corollary is crucial. For even $n$ (and $\lambda$
strong limit) the situation is different. In the lemmas below $\mu$ is a
cardinal, $k$ is an integer $\geq 1$ and $B$ is a Boolean algebra.

\begin{definition}
For a cardinal $\mu$ and an integer $k\in\omega$ we define $\daleth_k(\mu)$
inductively by\footnote{Remember that $\daleth$ (daleth) is the second letter
after $\beth$ (beth) in the Hebrew alphabet}
\[\daleth_0(\mu)=\mu,\quad\quad\quad\daleth_{k+1}(\mu)=
(2^{\daleth_k(\mu)})^{++}.\]  
\end{definition}

\begin{lemma}
\label{lem1}
\begin{enumerate}
\item Suppose that
\begin{description}
\item[($\oplus$)] \ \ \ $\ind_{2k}(B)\geq\beth_{2k}(\mu)^+$
\end{description}
or at least
\begin{description}
\item[($\oplus^-$)] \ \ \ there exists a sequence $\langle x_i:i<\beth_{2k}
(\mu)^+\rangle\subseteq B$ such that if $i_0<i_1<\ldots<i_{2k-2}<i_{2k-1}
<\beth_{2k}(\mu)^+$ then $\bigwedge\limits_{l<k}x_{i_{2l}}\wedge
(-x_{i_{2l+1}})\neq 0$.
\end{description}
Then
\begin{description}
\item[$\spadesuit^{B,k}_{\mu^+}$\ \ \ ] there is a sequence $\langle y_j:
j<\mu^+\rangle\subseteq B$ such that for each $w\in [\mu^+]^k$ there is
an ultrafilter $D\in\Ult B$ with
\[(\forall j<\mu^+)(y_j\in D\iff j\in w).\]
\end{description}
\item If $\ind_{2k}(B)\geq\daleth_{k+1}(\mu)$ then we can conclude
$\spadesuit^{B,k}_\mu$. 
\end{enumerate}
[In 2) it is enough to assume a suitable variant of $(\oplus^-)$: see the
proof.] 
\end{lemma}

\noindent PROOF: {\em 1.}\ \ \ Assume $(\oplus^-)$. For each $i_0<\ldots<
i_{2k-1}<\beth_{2k}(\mu)^+$ fix an ultrafilter $D^{\{i_0,\ldots,i_{2k-1}\}}
\in \Ult B$ such that $\bigwedge\limits_{l<k}(x_{i_{2k}}\wedge
(-x_{i_{2l+1}}))\in D^{\{i_0,\ldots,i_{2k-1}\}}$. Let $F:[\beth_{2k}
(\mu)^+]^{2k+1}\longrightarrow {}^{2k+1}2$ be defined by
\[F(\{i_0,\ldots,i_{2k}\})(l)=1 \iff x_{i_l}\in D^{\{i_0,\ldots,
i_{2k}\}\setminus\{i_l\}}\]
(where $l<2k+1$, $i_0<\ldots<i_{2k}<\beth_{2k}(\mu)^+$). By the
Erd\"os---Rado theorem we find a homogeneous for $F$ set $I$ of the size
$\mu^+$. We may assume that the sequence $\langle x_i: i<\mu^+\rangle$
behaves uniformly with respect to $F$.

Put $y_j=x_{\omega\cdot j}\wedge (-x_{\omega\cdot j + 5})$ for $j<\mu^+$. We
claim that the sequence $\langle y_j:j<\mu^+\rangle$ has the required
property. For this suppose that $j_0<\ldots<j_{k-1}<\mu^+$ and let 
$i_{2l}=\omega\cdot j_l$, $i_{2l+1}=\omega\cdot j_l +5$ (for $l<k$). Then
$i_0<\ldots<i_{2k-2}<i_{2k-1}<\mu^+$ so we can take $D=D^{\{i_0,\ldots,
i_{2k-1}\}}$. Thus $y_{j_l}=x_{i_{2l}}\wedge (-x_{i_{2l+1}})\in D$ for $l<k$.
On the other hand suppose that $j\notin\{j_0,\ldots,j_{k-1}\}$ and look at
$i=\omega\cdot j$, $i'=\omega\cdot j+5$. Note that for each $l<k$ we have
\[i<i_{2l}\iff i<i_{2l+1}\iff i'<i_{2l+1}\iff i'<i_{2l}.\]
Since $F(\{i,i_0,\ldots,i_{2k-1}\})=F(\{i',i_0,\ldots,i_{2k-1}\})$ we get that 
\[x_i\in D\iff x_{i'}\in D\]
and hence $y_j=x_i\wedge(-x_{i'})\notin D$.
\medskip

\noindent {\em 2.}\ \ \ The proof is essentially the same as above but
instead of the Erd\"os--Rado theorem we use \ref{basic} which is a special
case of the canonization theorems of \cite{Sh 95}. We start with a sequence
$\langle x_{\alpha,\xi}:\alpha<\daleth_{k+1}(\mu),\xi<\mu\rangle\subseteq B$
such that if $\xi_0,\ldots,\xi_{k-1}<\mu$, $\alpha^0_l<\alpha^1_l<
\daleth_{k+1}(\mu)$ (for $l<k$) then $\bigwedge\limits_{l<k}(x_{\alpha^0_l,
\xi_l}\wedge (-x_{\alpha^1_l,\xi_l}))\neq 0$. Then we choose the respective
ultrafilters $D^{\alpha^0_0\alpha^1_0\ldots\alpha^0_{k-1}\alpha^1_{k-1}}_{
\xi_0,\ldots,\xi_{k-1}}\in\Ult B$ and we consider a function $F:
[\daleth_{k+1}(\mu)\times\mu]^{2k+1}\longrightarrow 2$ such that
\begin{quotation}
\noindent $F((\alpha^0_0,\xi_0),(\alpha^1_0,\xi_0),\ldots,
(\alpha^0_{k-1},\xi_{k-1}),(\alpha^1_{k-1},\xi_{k-1}),(\alpha,\xi))=1$ 

\noindent if and only if $x_{\alpha,\xi}\in D^{\alpha^0_0\alpha^1_0\ldots
\alpha^0_{k-1}\alpha^1_{k-1}}_{\xi_0,\ldots,\xi_{k-1}}$. 
\end{quotation}
By \ref{basic} a) we find $\alpha^0_\xi,\alpha^1_\xi<\daleth_{k+1}(\mu)$ (for
$\xi<\mu$) such that for each distinct $\xi_0,\ldots,\xi_k\in\mu$
\[F((\alpha^0_{\xi_0},\xi_0),(\alpha^1_{\xi_0},\xi_0),\ldots,
(\alpha^0_{\xi_{k-1}},\xi_{k-1}), (\alpha^1_{\xi_{k-1}},\xi_{k-1}),
(\alpha^0_{\xi_k},\xi_k))=\] 
\[=F((\alpha^0_{\xi_0},\xi_0),(\alpha^1_{\xi_0},\xi_0),\ldots,
(\alpha^0_{\xi_{k-1}},\xi_{k-1}), (\alpha^1_{\xi_{k-1}},\xi_{k-1}),
(\alpha^1_{\xi_k},\xi_k)).\] 
Finally put $y_\xi=x_{\alpha_\xi^0,\xi}\wedge(-x_{\alpha^1_\xi,\xi})$.\QED

\begin{lemma}
\label{lem2}
Suppose that there is a sequence $\langle y_j: j<\mu\rangle\subseteq B$
such that for every $w\in [\mu]^k$ there is an ultrafilter $D\in\Ult B$
such that 
\[(\forall j<\mu)(y_j\in D\iff j\in w).\]
Then
\[\ind_{2k+1}^+(B\times B)>\mu.\]
\end{lemma}
PROOF: Consider the sequence $\langle (y_j,-y_j): j<\mu\rangle\subseteq
B\times B$. To prove that it is $2k+1$-independent suppose that
$j_0<\ldots<j_{2k}<\mu$, $t\in{}^{2k+1}2$. Let $w_0=\{j_l: t(l)=0\}$,
$w_1=\{j_l: t(l)=1\}$. One of these sets has at most $k$ elements so we find
an ultrafilter $D\in\Ult B$ such that
\begin{quotation}
\noindent either $(\forall l<2k+1)(y_{j_l}\in D\iff t(l)=0)$

\noindent or $(\forall l<2k+1)(y_{j_l}\in D\iff t(l)=1)$.
\end{quotation}
In the first case $\bigwedge\limits_{l<2k+1}y^{t(l)}_{j_l}\in D$, in the
second case $\bigwedge\limits_{l<2k+1}(-y_{j_l})^{t(l)}\in D$. Consequently
$(\bigwedge\limits_{l<2k+1}y^{t(l)}_{j_l},\bigwedge\limits_{l<2k+1}
(-y_{j_l})^{t(l)})\neq 0$ and the lemma is proved. \QED

\begin{theorem}
\label{ind}
Let $k$ be an integer $\geq 1$, $B$ a Boolean algebra, $\lambda$ a cardinal.
Then 
\begin{enumerate}
\item $\ind_{2k}(B)\geq\max\{\beth_{2k}(\lambda)^+,
\daleth_{k+1}(\lambda^+)\}$ implies $\ind_{2k+1}(B\times B)\geq \lambda^+$.
\item If $\lambda$ is strong limit, $\ind_{2k}(B)\geq\lambda$
then $\ind_{2k+1}(B\times B)\geq\lambda$.
\item $\ind_{2k}(\prod^w_{i<\omega}B)<\beth_\omega(\ind_{2k+1}(\prod^w_{i<
\omega}B))$.
\end{enumerate}
\end{theorem}
PROOF: {\em 1.}\ \ It is an immediate consequence of lemmas~\ref{lem1} and
\ref{lem2}.
\medskip

\noindent {\em 2.}\ \ Follows from {\em 1.}
\medskip

\noindent {\em 3.}\ \ It follows from {\em 2} and the following observation.
\begin{claim}
For an integer $n>1$ and a Boolean algebra $B$ we have
\[\ind_n({\prod}^w_{i<\omega} B)=\ind_n({\prod}^w_{i<\omega} B\times
{\prod}^w_{i<\omega} B).\] 
\end{claim}
Proof of the claim: By \ref{indprofac}(1) we have
\[\ind_n({\prod}^w_{i<\omega} B)\leq \ind_n({\prod}^w_{i<\omega} B\times
{\prod}^w_{i<\omega} B).\]
For the other inequality assume that 
\[\kappa\stackrel{\rm def}{=}\ind_n({\prod}^w_{i<\omega} B)<\ind_n({\prod}^w_{
i<\omega} B\times {\prod}^w_{i<\omega} B).\]
Thus we find an $n$--independent set $X\subseteq {\prod}^w_{i<\omega} B\times
{\prod}^w_{i<\omega} B$ of size $\kappa^+$. For $x\in X$ let $a_x,b_x\in
{\prod}^w_{i<\omega} B$ and $m(x)<\omega$ be such that 
\[x=(a_x,b_x)\ \mbox{ and } \ (\forall m\geq m(x))(a_x(m)=a_x(m(x))\ \&\
b_x(m)=b_x(m(x))).\]
Take $m_0<\omega$ and $Y\in [X]^{\textstyle \kappa^+}$ such that $m(x)=m_0$
for $x\in Y$. For $x\in Y$ let 
\[c_x\stackrel{\rm def}{=} (a_x(0),\ldots,a_x(m_0),b_x(0),\ldots,b_x(m_0))\in
B^{2m_0+2}.\] 
The set $Z\stackrel{\rm def}{=}\{c_x: x\in Y\}$ is $n$-independent as $a_x(m)
=a_x(m_0)$, $b_x(m)=b_x(m_0)$ for $m\geq m_0$. As $\|Z\|=\kappa^+$ we conclude
that $\kappa^+\leq\ind_n(B^{2m_0+2})$. Now note that the algebras ${\prod}^w_{
i<\omega} B$ and $B^{2m_0+2}\times {\prod}^w_{i<\omega} B$ are isomorphic, so
(by \ref{indprofac}) 
\[\ind_n(B^{2m_0+2})\leq \ind_n({\prod}^w_{i<\omega} B),\]
and hence $\kappa^+\leq \ind_n({\prod}^w_{i<\omega} B)=\kappa$, a
contradiction.  \QED 

\begin{problem}
\begin{enumerate}
\item Can Lemma~\ref{lem1} be improved? Can we (consistently?) weaken the
variant of the assumption $(\oplus^-)$ for 2) to sequences shorter than
$\daleth_k(\mu)$ (we are interested in the reduction of the steps in the beth
hierarchy)? 
\item Describe (in ZFC) all dependences between $\ind_k(B^n)$ (for
$n,k<\omega$) [note that we may force them distinct].
\end{enumerate}
\end{problem}

\subsection{Tightness.}
The tightness $t(B)$ of a Boolean algebra $B$ is the minimal cardinal $\kappa$
such that if $F$ is an ultrafilter on $B$, $Y\subseteq{\rm Ult}B$ and $F
\subseteq\bigcup Y$ then there is $Z\in [Y]^{\leq\kappa}$ such that $F
\subseteq \bigcup Z$. To represent the tightness as a def.u.w.o.car.~invariant
we use the following characterization of it (see \cite{M}): 
\[t(B)=\sup\{\|\alpha\|:\mbox{ there exists a free sequence of the length
}\alpha\mbox{ in }B\}\] 
where a sequence $\langle x_\xi:\xi<\alpha\rangle\subseteq B$ is free if
\[(\forall\xi<\alpha)(\forall F\in [\xi]^{<\omega})(\forall G\in
[\alpha\setminus\xi]^{<\omega})[ \bigwedge_{\eta\in F}x_\eta\wedge
\bigwedge_{\eta\in G}-x_\eta\neq 0].\]
Now it is easy to represent $t(B)$ as def.u.w.o.car.~invariant. Together with
(finite versions of) the tightness we will define a def.f.o.car.~invariant
$\ut_k$ which is inspired by \ref{lem1}.

\begin{definition}
\begin{enumerate}
\item Let $\psi$ be the sentence saying that $P_1$ is a well ordering of $P_0$
(we denote the respective order by $<_1$). For $k,l<\omega$ let $\phi_{k,l}^t$
be the sentence asserting that 
\begin{quotation}
\noindent for each $x_0,\ldots,x_k,y_0,\ldots,y_l\in P_0$

\noindent if $x_0<_1\ldots<_1 x_k<_1 y_0<_1\ldots <_1 y_l$ then
$\bigwedge\limits_{i\leq k}x_i\not\leq\bigvee\limits_{i\leq l} y_i$,
\end{quotation}
and let the sentence $\phi_{k,l}^{\ut}$ say that
\[\mbox{for each distinct }x_0,\ldots,x_k,y_0,\ldots,y_l\in P_0\mbox{ we have }
\bigwedge\limits_{i\leq k}x_i\not\leq\bigvee\limits_{i\leq l} y_i.\]
\item For $n,m\leq\omega$ let $T^{n,m}_t=\{\phi_{k,l}^t:
k<n,l<m\}\cup\{\psi\}$ and $T^{n,m}_{\ut}=\{\phi^{\ut}_{k,l}: k<n,
l<m\}$ and for a Boolean algebra $B$:
\[t_{n,m}(B)=\inv_{T^{n,m}_t}(B)\quad\quad\&\quad\quad
\ut_{n,m}(B)=\inv_{T^{n,m}_{\ut}}(B).\]
\item {\em The unordered $k$-tightness} $\ut_k$ is the def.f.o.car.~invariant
$\ut_{k,\omega}$. 
\end{enumerate}
\end{definition}
{\bf Remark:}\ \ \ Note that $T^{n,m}_t =\{\psi\}$ if either $n=0$ or $m=0$
(and thus $t_{n,m}(B)=\|B\|$ whenever $n\cdot m = 0$). The theory $T^{n+1,m+
1}_t$ says that $P_1$ is a well ordering of $P_0$ and if $x_0<_1\ldots<_1 x_n
<_1 y_0<_1\ldots<_1 y_m$ then the meet $\bigwedge\limits_{i\leq n} x_i$ is not
covered by the union $\bigvee\limits_{i\leq m} y_i$. The invariant $t_{\omega,
\omega}(B)$ is just the tightness of $B$. Similarly for $T^{n,m}_{\ut}$.

\begin{corollary}
\label{co5.12}
For a Boolean algebra $B$ and $n,m\leq\omega$, $0<k<\omega$:
\begin{enumerate}
\item $\ind^{(+)}_{n+m}(B)\leq\ut^{(+)}_{n,m}(B)=\ut^{(+)}_{m,n}(B)\leq
t^{(+)}_{n,m}(B)$,
\item $\ut^{(+)}_{k}(B)=\sup\{\kappa^{(+)}:\spadesuit^{B,k}_\kappa\mbox{
holds true}\}$, where the condition $\spadesuit^{B,k}_\kappa$ is as defined
in lemma~\ref{lem1},
\item the condition $\spadesuit^{B,k}_\kappa$ is equivalent to:
\begin{quotation}
\noindent{\em the algebra $B^k_0(\kappa)$ of \ref{pr3.1} can be
embedded into a homomorphic image of $B$}, 
\end{quotation}
\item $\ut^{(+)}_k(B)\leq\ind^{(+)}_{2k}(B)$,\ \ \ $\ut^{(+)}_k(B)\leq
\ind^{(+)}_{2k+1}(B\times B)$.
\end{enumerate}
\end{corollary}

\noindent PROOF: {\em 1.} and {\em 2.} should be clear.
\medskip

\noindent{\em 3.} Assume $\spadesuit^{B,k}_\kappa$ and let $\langle y_j:
j<\kappa\rangle\subseteq B$ be a sequence witnessing it. Let $I$ be the
ideal of $B$ generated by the set
\[\{y_{j_0}\wedge\ldots \wedge y_{j_k}: j_0<\ldots<j_k<\kappa\}.\]
Then the algebra $B^k_0(\kappa)$ naturally embeds into the quotient algebra
$B/I$. Moreover, if $B'$ is a homomorphic image of $B$ and
$\spadesuit^{B',k}_\kappa$ then clearly $\spadesuit^{B,k}_\kappa$ so the
converse implication holds true too.
\medskip

\noindent{\em 4.} It follows from {\em 3.} and (the proof of)
Proposition~\ref{pr3.1} and the remark after the proof of \ref{pr3.1}. \QED  
\medskip

\noindent{\bf Remark}:\ \ \ Corollary \ref{co5.12}(3) is specially interesting
if you remember that $s^+(B)>\lambda$ if and only if the finite---cofinite
algebra on $\lambda$ can be embedded into a homomorphic image of $B$.
\medskip

From Lemma~\ref {lem1} we can conclude the following:

\begin{corollary} For $k>0$ and an algebra $B$: 
\begin{enumerate}
\item if either $\ut_{k,k}(B)>\beth_{2k}(\mu)$ or $\ut_{k,k}(B)\geq
\daleth_{k+1}(\mu^+)$ then $\ut_k(B)\geq\mu^+$, 
\item if $\lambda$ is strong limit, $\ut_{k,k}(B)\geq\lambda$ then
$\ut_k(B)\geq\lambda$,
\item $\ut_{k,k}(B)<\beth_\omega(\ut_k(B))$. \QED
\end{enumerate}
\end{corollary}

\begin{proposition}
Suppose $n,m<\omega$, $k=\min\{n,m\}$, $B$ is a Boolean algebra. Then
\[t_{n,m}(B)\leq\beth_{n+m}(\ut_k(B)+t(B)).\]
\end{proposition}

\noindent PROOF: Let $\mu=\ut_k(B)+t(B)$ and assume that $t_{n,m}(B)>
\beth_{n+m}(\mu)$. Then we have a sequence $\langle a_\alpha:
\alpha<\beth_{n+m}(\mu)^+\rangle\subseteq B$ such that 
\[(\forall \alpha_0<\ldots<\alpha_{n+m-1}<\beth_{n+m}(\mu)^+)[
(\bigwedge\limits_{l<n} a_{\alpha_l} \wedge 
\bigwedge\limits_{n\leq l<n+m}-a_{\alpha_l})\neq 0].\]
For each $\alpha_0,\ldots,\alpha_{n+m-1}$ as above fix an ultrafilter
$D^{\{\alpha_0,\ldots,\alpha_{n+m-1}\}}\in \Ult B$ containing the element
$\bigwedge\limits_{l<n} a_{\alpha_l} \wedge \bigwedge\limits_{n\leq
l<n+m}-a_{\alpha_l}$. Look at the function 
\[F:[\beth_{n+m}(\mu)^+]^{n+m+1}\longrightarrow {}^{n+m+1}2\]
defined by
\[F(\alpha_0,\ldots,\alpha_{n+m})(l)=1\iff a_{\alpha_l}\in
D^{\{\alpha_0,\ldots,\alpha_{n+m}\}\setminus\{\alpha_l\}}.\] 
By the Erd\"os-Rado theorem we may assume that $\mu^+$ is homogeneous for $F$
with the constant value $c\in {}^{n+m+1}2$.

If $c(l)=0$ for each $l\leq n+m$ then the sequence $\langle a_\alpha:\alpha<
\mu^+\rangle$ witnesses $\mu^+\leq\ut_n(B)$ giving a contradiction to the
definition of $\mu$ (remember $\ut_k(B)\geq\ut_n (B)$). In fact, given $n$
elements $\alpha_0<\ldots<\alpha_{n-1}$, choose $m$ additional elements
$\alpha_{n-1}<\alpha_n<\ldots<\alpha_{n+m-1}$. Suppose that $\beta\in\mu^+
\setminus\{\alpha_0,\ldots,\alpha_{n+m-1}\}$. Then by homogeneity $-a_\beta\in
D^{\{\alpha_0,\ldots,\alpha_{n+m-1}\}}$, proving the result.

If $c(l)=1$ for each $l$ then the sequence $\langle -a_\alpha:n\leq\alpha<
\mu^+\rangle$ exemplifies $\mu^+\leq\ut_m(B)$, once again a contradiction. In
fact, take any $m$ elements $n-1<\alpha_n<\ldots<\alpha_{n+m-1}$ and suppose
that $\beta\in \mu^+\setminus\{0,\ldots,n-1,\alpha_n,\ldots,\alpha_{n+m-1}\}$.
Then by homogeneity $a_\beta\in D^{\{0,\ldots,n-1,\alpha_n,\ldots,
\alpha_{n+m-1}\}}$, as desired.

Finally, suppose that there are $l_0,l_1\leq n+m$ such that $c(l_0)=0$ and
$c(l_1)=1$. \\
{\sc Case 1:}\ \ \ $l_1<l_0$\\
Let $\Gamma=\{\beta+\omega:\beta<\mu^+\}$. We claim that $\langle a_\alpha:
\alpha\in\Gamma\rangle$ witnesses $\mu^+\leq t(B)$, contradicting $\mu\geq
t(B)$. In fact, let $\alpha_0<\ldots<\alpha_p<\ldots<\alpha_{q-1}$ be elements
of $\Gamma$; we want to show that
\[\bigwedge_{l<p} a_{\alpha_l} \wedge \bigwedge_{p\leq l<q} -a_{\alpha_l}\neq
0.\]
Say $\alpha_p=\beta+\omega$. Define $\gamma_l=l$ for all $l<l_1$,
$\gamma_{l_1},\ldots,\gamma_{l_0-1}$ are consecutive values starting with
$\beta+1$, and $\gamma_{l_0},\ldots,\gamma_{m+n-1}$ are consecutive values
starting with $\alpha_{q-1}+1$ (none of the latter if $l_0=n+m$). Then
$a_{\alpha_l}\in D^{\{\gamma_0,\ldots,\gamma_{n+m-1}\}}$ for all $l<p$ and
$-a_{\alpha_l}\in D^{\{\gamma_0,\ldots,\gamma_{n+m-1}\}}$ for all $l\geq p$,
as desired.\\ 
{\sc Case 2:}\ \ \ $l_1\geq l_0$\\
This is similar, using $\langle -a_\alpha: \alpha\in\Gamma\rangle$. \QED
\medskip

Our next proposition is motivated by Theorem~\ref{ind} and the above
corollaries. 

\begin{proposition} 
Let $B$ be a Boolean algebra, $k$ a positive integer. Then 
\begin{enumerate}
\item $\ind_{2k}(\prod_{i<\omega}^w B)\leq\min\{\beth_{2k-1}(\ind_k(B)),
\beth_{2k-1}(\ut_k(B))\}$,
\item $\ut_{k+1}(\prod^w_{i<\omega}B)\leq\beth_k(\ut^+_{k+1}(B))$.
\end{enumerate}
\end{proposition}

\noindent PROOF: {\em 1.} Suppose that $\lambda_0=\beth_{2k-1}(\ind_k(B))<
\ind_{2k}(\prod^w_{i<\omega} B)$. Thus we find a sequence $\langle a_\alpha:
\alpha<\lambda_0^+\rangle\subseteq\prod^w_{i<\omega} B$ which is
$2k$-independent. Let $a_\alpha=\langle a_\alpha(i):i<\omega\rangle$ (for
$\alpha<\lambda_0^+$). Consider the function $F: [\lambda_0^+]^{2k}
\longrightarrow\omega$ given by $F(\alpha_0,\ldots,\alpha_{2k-1})=$
\[\min\{i\in\omega: B\models a_{\alpha_0}(i)\wedge(-a_{\alpha_1}(i))\wedge
\ldots\wedge a_{\alpha_{2k-2}}(i)\wedge(-a_{\alpha_{2k-1}}(i))\neq 0\},\] 
where $\alpha_0<\ldots<\alpha_{2k-1}<\lambda_0^+$. By the Erd\"os--Rado
theorem we find a set $I$ of the size $(\ind_k(B))^+$ homogeneous for $F$;
we may assume that $I=(\ind_k(B))^+$. Let $i_0$ be the constant value of $F$
(on $[(\ind_k(B))^+]^{2k}$). Look at the sequence $\langle a_\alpha(i_0):
\alpha<(\ind_k(B))^+\ \&\ \alpha\mbox{ limit}\rangle$. Any combination of
$k$ members of this sequence can be ``extended'' to a combination of $2k$
elements of $\langle a_\alpha(i_0):\alpha<(\ind_k(B))^+\rangle$ of the type
used in the definition of $F$. A contradiction.

Now suppose that $\lambda_1\stackrel{\rm def}{=}\beth_{2k-1}(\ut_k(B))<
\ind_{2k}({\prod}^w_{i<\omega}B)$. Like in \ref{lem1}, we take a sequence
$\langle a_\alpha: \alpha<\lambda_1^+\rangle\subseteq {\prod}^w_{i<\omega}B$
such that for some $n<\omega$, for each $\alpha<\lambda_1^+$, $a_\alpha\in
B^n$ (i.e. the support of $a_\alpha$ is contained in $n$) and
\[(\forall\alpha_0<\ldots<\alpha_{2k-1}<\lambda_1^+)(\bigwedge_{l<k}
a_{\alpha_{2l}} \wedge (-a_{\alpha_{2l+1}})\neq 0),\] 
and for each $\alpha_0<\ldots<\alpha_{2k-1}<\lambda_1^+$ we choose an
ultrafilter $D^{\{\alpha_0,\ldots,\alpha_{2k-1}\}}\in\Ult {\prod}^w_{i<\omega}
B$ such that 
\[\bigwedge_{l<k} a_{\alpha_{2l}}\wedge (-a_{\alpha_{2l+1}})\in D^{\{\alpha_0,
\ldots,\alpha_{2k-1}\}}.\]
Now we consider a colouring $F:[\lambda_1^+]^{\textstyle 2k+1} \longrightarrow
{}^{2k+1} (2\times n)$ given by
\[\begin{array}{ll}
F(\{\alpha_0,\ldots,\alpha_{2k}\})(l)=(1,m)\iff& a_{\alpha_l}\in
D^{\{\alpha_0,\ldots,\alpha_{2k}\}\setminus\{\alpha_l\}}\quad\mbox{ and}\\
\ &D^{\{\alpha_0,\ldots,\alpha_{2k}\}\setminus\{\alpha_l\}}\mbox{ is
concentrated on}\\
\ &\mbox{the $m^{\rm th}$ coordinate.}\\
  \end{array}\]
By Erd\H os--Rado theorem we may assume that the set of the first
$(\ut_k(B))^+$ elements of $\lambda_1^+$ is homogeneous for $F$. Now we finish
as in \ref{lem1} notifying that for some $m<n$, for all $\alpha_0<\ldots<
\alpha_{2k-1}<(\ut_k(B))^+$ the ultrafilter $D^{\{\alpha_0,\ldots,
\alpha_{2k-1}\}}$ is concentrated on the $m^{\rm th}$ coordinate. So we may
use elements of the form $a_{\alpha\cdot\omega}(m)\wedge (-a_{\alpha\cdot
\omega+5}(m))$ (for $\alpha<(\ut_k(B))^+$) to get a contradiction. 
\medskip

\noindent{\em 2.} Assume that $\ut_{k+1}(\prod_{i<\omega}^w B)>\beth_k(\mu)$,
$\mu=\ut_{k+1}^+(B)$. Then we find a sequence $\langle a_\alpha:\alpha<(
\beth_k(\mu))^+\rangle\subseteq\prod_{i<\omega}^w B$ such that for any $k+1$
distinct members of this sequence there is an ultrafilter containing all of
them and no other member of the sequence. We may assume that for some $n<
\omega$ we have $\langle a_\alpha:\alpha<(\beth_k(\mu))^+\rangle\subseteq
B^n$. For $\alpha_0,\ldots,\alpha_k<(\beth_k(\mu))^+$ let $D^{\alpha_0,\ldots,
\alpha_k}$ be the respective ultrafilter of $B^n$ (i.e.~it contains all
$a_{\alpha_l}$ (for $l\leq k$) and nothing else from the sequence) and let
$F(\alpha_0,\ldots,\alpha_k)<n$ be such that the ultrafilter $D^{\alpha_0,
\ldots,\alpha_k}$ is concentrated on that coordinate. By the Erd\H os--Rado
theorem we find a set $A\in [(\beth_k(\mu))^+]^{\mu^+}$ homogeneous for $F$.
Let $m$ be the constant value of $F$ on $A$. Look at the sequence $\langle
a_\alpha(m): \alpha\in A\rangle$ -- it witnesses $\spadesuit^{B,k+1}_{\mu^+}$
contradicting $\mu=\ut^+_{k+1}(B)$.  \QED 
\medskip

Finally note that for the algebra $B^k_0(\lambda)$ of \ref{pr3.1} we have:
\[\ut_k(B^k_0(\lambda))=t_{k,\omega}(B^k_0(\lambda))=\lambda\quad\quad
\mbox{ and}\]
\[\ut_{k+1}(B^k_0(\lambda))=\ut_{k+1,k+1}(B^k_0(\lambda))= t_{k+1,k+1}
(B^k_0(\lambda))=\aleph_0.\]
This gives us an example distinguishing $t_{k,\omega}$ and $t_{k+1,\omega}$
(and in corollary~\ref{co5.3} we may replace $\ind$ by $t$). But the
following problem remains open: 

\begin{problem}
Are the following inequalities possible?:
\[t_{k,\omega}(B)>\ut_k(B),\quad t_{\omega,k}(B)>\ut_k(B),\quad
t_{k,k}(B)>t_{k,k+1}(B).\]
\end{problem}

\subsection{Independence and interval Boolean algebras.}
Now we are going to reformulate (in a stronger form) and put in our
general setting the results of \cite{Sh 503}.

\begin{definition}
Let $B$ be a Boolean algebra.
\begin{enumerate}
\item For a filter $D$ on $[\lambda]^k$ we say that $B$ has the
{\em $D$-dependence property} if for every sequence $\langle a_i: i<\lambda
\rangle\subseteq B$ there is $A\in D$ such that for every
$\{\alpha_0,\alpha_1,\dots,\alpha_{k-1}\} \in A$ the set 
$\{a_{\alpha_0}, a_{\alpha_1},\ldots, a_{\alpha_{k-1}}\}$ is not independent.
  
\item For a filter $D$ on $[\lambda]^k$ and  a Boolean term $\tau(x_0,
x_1,\ldots, x_{k-1})$ we say that $B$ has the {\em $(D,\tau)$-dependence
property} if for every $\langle a_i:i<\lambda\rangle\subseteq B$, for some
$A\in D$, for every $\{\alpha_0, \alpha_1,\ldots,\alpha_{k-1}\} \in A$ with
$\alpha_0<\alpha_1<\ldots<\alpha_{k-1}$ we have
$B\models\tau(a_{\alpha_0},a_{\alpha_1},\dots,a_{\alpha_{k-1}})=0$. 
\end{enumerate}
\end{definition}
It should be clear that if $D$ is a proper filter on $[\lambda]^k$ and a
Boolean algebra $B$ has the $D$--dependence property then $\lambda\geq
\ind_k^+(B)$ (and so $\lambda\geq\ind^+(B)$). 

\begin{proposition}
\label{pr3.7A}
Let $\tau=\tau(x_0,x_1,\dots,x_{k-1})$ be a Boolean term and let $D$ be a
$\kappa$-complete filter on $[\lambda]^k$. Then any reduced product of $< 
\kappa$ Boolean algebras  having the $(D,\tau)$-dependence property has the
$(D,\tau)$-dependence property (this includes products and ultraproducts).
\QED 
\end{proposition}

\begin{proposition}
\label{pr3.7B}
Assume  $D$ is a proper filter on $[\lambda]^k$. Then there exists a sequence
$\langle\alpha_0,\alpha_1,\ldots,\alpha_{k-1}\rangle$ of ordinals $\leq
\lambda$  such that:   
\begin{description}
\item[(a)] $\{w\in [\lambda]^k:$ for each $\ell < k$ the $\ell$-th  member
of $w$ is $<\alpha_\ell\}\neq\emptyset\ \mod D$,
\item[(b)] if $\alpha_\ell^\prime\leq\alpha_\ell$ for all $\ell<k$, $n<k$ and
$\alpha_n^\prime<\alpha_n$ then 
\[\{w\in [\lambda]^k: \mbox{ for each $\ell<k$, the } \ell\mbox{-th member of
}w \mbox{ is }<\alpha_\ell^\prime\} = \emptyset\ \mod D.\]
\end{description}
[Note that necessarily $\langle\alpha_\ell:\ell<k\rangle$ is non-decreasing.] 
\end{proposition}

\noindent PROOF: Let $F$ be the set of all non-decreasing sequences
$\langle\alpha_\ell:\ell<k\rangle\subseteq\lambda+1$ such that the condition
{\bf (a)} holds. Then $F$ is upward closed (and $\langle\lambda,\ldots,\lambda
\rangle\in F$). Choose by induction $\alpha_0,\ldots,\alpha_{k-1}$ such that
for each $\ell< k$ 
\[\alpha_\ell=\min\{\beta: (\exists \bar{\alpha}\in F)(\bar{\alpha}
\rest\ell=\langle \alpha_0,\ldots,\alpha_{\ell-1}\rangle\ \&\
\alpha_\ell=\beta)\}. \QED\]

\begin{definition}
We call a filter $D$ on $[\lambda]^k$ {\em normal for} $\langle\alpha_0,
\alpha_1,\ldots,\alpha_{k-1}\rangle$ if condition {\bf (b)} of \ref{pr3.7B}
holds and 
\begin{description}
\item[(a)$^+$]  $\{w\in [\lambda]^k:$ for each $\ell<k$ the $\ell-$th member
of $w$ is  $<\alpha_\ell\}\in D$.
\end{description}
\end{definition}

\begin{proposition}  
Assume that
\begin{enumerate}
\item $D$ is a $\kappa$-complete filter on $[\lambda]^k$ which is normal for
$\langle\alpha_0,\alpha_1,\ldots,\alpha_{k-1}\rangle$, and
$\alpha_0,\ldots,\alpha_{k-1}$ are limit ordinals,  
\item $k(*) = k\cdot 2^k$, $i\mapsto (m_i,l_i): k(*)\longrightarrow 2^k\times
k$ is a one-to-one mapping such that $i_1<i_2$ implies that, 
lexicographically, $(\alpha_{\ell_{i_1}}, m_{i_1},\ell_{i_1})<(\alpha_{\ell_{
i_2}},m_{i_2},\ell_{i_2})$; for $(m,\ell)\in 2^k\times k$ the unique $i<k(*)$
such that $(m_i,\ell_i)=(m,\ell)$ is denoted by $i(m,\ell)$, 
\item $\kappa^*$ is a regular cardinal such that $(\forall \mu<\kappa^*)(
2^\mu<\kappa)$ (e.g. $\kappa^* =\aleph_0$),
\item for $X\in D$, $h:X\longrightarrow\mu$, $\mu<\kappa^*$:
\[A_{X,h}\stackrel{\rm def}{=}\{w\in [\lambda]^{k(*)}: (\forall m,m'<2^k)(w_m
\in X\ \& \ h(w_{m})=h(w_{m'}))\}\]
where for $m<2^k$, $w=\{\beta_0,\ldots,\beta_{k(*)-1}\}\in [\lambda]^{k(*)}$
(the increasing enumeration) the set $w_m$ is $\{\beta_{i(m,\ell)}:\ell<k\}$,
\item $D^*$ is the $\kappa^*$-complete filter on $[\lambda]^{k(*)}$ generated
by the family 
\[\{A_{X,h}: X\in D, h:X\longrightarrow\mu,\mu<\kappa^*\},\]
\item $\tau^*=\tau^*(x_0,x_1,\ldots,x_{k(*)-1})=\bigwedge\limits_{m<2^k}
\bigwedge\limits_{\ell<k} x^{f_m(\ell)}_{i(m,\ell)}$, where $\langle f_m:m<
2^k\rangle$ lists all the functions in ${}^k2$.
\end{enumerate}
\noindent Then 
\begin{description}
\item[(a)] $D^*$ is a proper $\kappa^*$-complete filter on $[\lambda]^{k(*)}$
which is normal for the sequence $\langle \alpha_{l_{i}}: i<k(*)\rangle$,
\item[(b)] if a Boolean algebra $B$ has the $D$-dependence property then it
has the $(D^*,\tau^*)$-dependence property.
\end{description}
\end{proposition}

\noindent PROOF: Assume that $X_j\in D$, $\mu_j<\kappa^*$, $h_j:X_j
\longrightarrow\mu_j$ for $j<\mu<\kappa^*$ and look at the intersection
$\bigcap\limits_{j<\mu}A_{X_j,h_j}$. Let $X^*=\bigcap\limits_{j<\mu}X_j$. Then
$X^*\in D$ as $\mu<\kappa$ and $D$ is $\kappa$--complete. Moreover for some
$\langle\xi_j:j<\mu\rangle\in\prod\limits_{j<\mu}\mu_j$ we have  
\[X^+\stackrel{\rm def}{=}\{w\in X^*: (\forall j<\mu)(h_j(w)=\xi_j)\}\neq
\emptyset\ \mod D,\] 
as $\prod\limits_{j<\mu}\mu_j<\kappa$ (remember $\kappa^*$ is regular and
$(\forall\mu<\mu^*)(2^\mu<\kappa)$). Let $r_0<r_1<\ldots<r_{\ell^*-1}<r_\ell^* 
=k-1$ be such that
\[\begin{array}{l}
\alpha_0=\ldots=\alpha_{r_0}<\alpha_{r_0+1}=\ldots=\alpha_{r_1}<
\alpha_{r_1+1}=\ldots\\
\ldots=\alpha_{r_{\ell^*-1}}<\alpha_{r_{\ell^*-1}+1}=\ldots=\alpha_{k-1}.
  \end{array}\] 
Now we choose inductively $\{\beta^m_0,\ldots,\beta^m_{k-1}\}\in X^+$ (for
$m<2^k$) such that 
\begin{quotation}
\noindent $\beta^m_n<\alpha_n$ for $n<k$, $m<2^k$;\ \ $\alpha_{r_u}<\beta^0_{
r_u+1}$ for $u<\ell^*$; 

\noindent $\beta^m_n<\beta^m_{n+1}$ for $n<k-1$, $m<2^k$,\ \ and

\noindent $\beta^m_{r_0}<\beta^{m+1}_0$,
$\beta^m_{r_{u+1}}<\beta^{m+1}_{r_u+1}$ for $u<\ell^*$, $m+1<2^k$. 
\end{quotation}
How? Since $D$ is normal for $\langle\alpha_0,\ldots,\alpha_{k-1}\rangle$ and
the $\alpha_i$'s are limit, the set 
\[\begin{array}{ll}
Y_0\stackrel{\rm def}{=} \{w\in [\lambda]^k:&\mbox{for each $n<k$ the $n$-th
member of $w$ is $<\alpha_n$ and}\\
\ &\mbox{for each $u<\ell^*$ the $r_u+1$-th element of $w$ is
$>\alpha_{r_u}$}\} 
  \end{array}\]
is in $D$. Thus we may choose $w_0=\{\beta^0_0,\ldots,\beta^0_{k-1}\}$ in
$X^+\cap Y_0$. Now suppose that we defined $\{\beta^m_0,\ldots,
\beta^m_{k-1}\}$. The set 
\[\begin{array}{ll}
Y_{m+1}\stackrel{\rm def}{=} \{w\in [\lambda]^k:&\mbox{for each $n<k$ the
$n$-th member of $w$ is $<\alpha_n$ and}\\
\ &\mbox{for each $u<\ell^*$ the $r_u+1$-th element of $w$ is
$>\beta^m_{r_{u+1}}$}\\ 
\ &\mbox{and the minimal element of $w$ is $>\beta^m_{r_0}$}\} 
  \end{array}\]
is in $D$ and we choose $w_{m+1}=\{\beta^{m+1}_0,\ldots,\beta^{m+1}_{k-1}\}$
in $X^+\cap Y_{m+1}$. Note that then $i_0<i_1\quad\Rightarrow\quad\beta^{
m_{i_0}}_{\ell_{i_0}}<\beta^{m_{i_1}}_{\ell_{i_1}}$ (for $i_0,i_1<k\cdot 2^k$)
and hence easily $w\stackrel{\rm def}{=}\{\beta^m_\ell: \ell<k,m<2^k\}\in
\bigcap\limits_{j<\mu} A_{X_j,h_j}$. Consequently the $\kappa^*$--complete
filter $D^*$ generated on $[\lambda]^{k(*)}$ by the sets $A_{X,h}$ is proper.
The filter $D^*$ is normal for $\langle\alpha_{\ell_i}:i<k(*)\rangle$ since:\\
if $X=\{\{\beta_0,\ldots,\beta_{k-1}\}\in [\lambda]^k: (\forall n<k)(\beta_n
<\alpha_n)\}$, $h$ is a constant function on $X$ then
\[A_{X,h}=\{\{\beta_0,\ldots,\beta_{k(*)-1}\}\in [\lambda]^{k(*)}: (\forall
i<k(*))(\beta_i<\alpha_{l_i})\}\in D^*;\]
if $i<k(*)$, $\alpha'<\alpha_{\ell_i}$ then the complement $X$ of the set
\[\{w\in [\lambda]^k:\mbox{ the }\ell_i\mbox{-th member of }w\mbox{ is less
then }\alpha'\}\]
is in $D$, and if $h$ is a constant function on $X$ then the set $A_{X,h}$
witnesses that
\[\{w\in [\lambda]^{k(*)}: \mbox{ the } i\mbox{-th member of } w \mbox{ is
less then }\alpha'\}=\emptyset\ \mod D^*.\]

\noindent It should be clear that the $D$-dependence property for $B$ implies
$(D^*,\tau^*)$-dependence property. \QED
\medskip

This is relevant to the product of linear orders. It was proved in \cite{Sh
503} that if $\kappa$ is an infinite cardinal, $B_\zeta$ (for $\zeta<\kappa$)
are interval Boolean algebras then $\ind(\prod\limits_{\zeta<\kappa}B_\zeta)
=2^\kappa$. The next result was actually hidden in the proof of Theorem
1.1 of \cite{Sh 503}.
  
\begin{theorem}
Let $\kappa$ be an infinite cardinal and let $\mu$ be a regular cardinal such
that for every $\chi<\mu$ we have $\chi^\kappa<\mu$ (e.g.~$\mu=(2^\kappa)^+$
in {\bf (1)} below or $\mu=(2^{2^\kappa})^+$ in {\bf (2)}). 
\begin{description}
\item[(1)] For a regressive function $f:\mu\longrightarrow\mu$ (i.e.~$f(
\alpha)<1+\alpha$), a two-place function $g:\mu^2\longrightarrow\chi$ for some
$\chi<\mu$ and for a closed unbounded subset $C$ of $\mu$ we put:     
\[\begin{array}{rl}
A_{C,f,g}=\{\{\alpha_0,\ldots,\alpha_5\}\in [\mu]^6: &
\alpha_0<\alpha_1<\dots<\alpha_5 \mbox{ are from }C,\\
\ & \mbox{each has cofinality}>\kappa,\\
\ & f(\alpha_0 )=f(\alpha_1)= \ldots  =f(\alpha_5)\ \mbox{ and}\\
g(\alpha_0,\alpha_1)= & g(\alpha_0,\alpha_2)=g(\alpha_3,\alpha_4)=
g(\alpha_3,\alpha_5)\}.
\end{array}\]
Let $D^6_{\mu,\kappa}$ be the filter on $[\mu]^6$ generated by all the sets
$A_{C,f,g}$. Finally, let $\tau_6$ be the following Boolean term: 
\[\tau_6(x_0,x_1,\dots,x_5)\stackrel{\rm def}{=} x_0\wedge (-x_1)\wedge x_2
\wedge (-x_3)\wedge x_4\wedge(-x_5).\] 
\end{description}
Then $D^6_{\mu,\kappa}$ is a proper $\kappa^+$-complete filter normal for
$\langle\mu,\mu,\mu,\mu,\mu,\mu\rangle$ and every interval Boolean algebra has
the $(D^6_{\mu,\kappa},\tau_6)$-dependence property.
\begin{description}
\item[(2)] Let $\mu_0$ be a cardinal such that $(\mu_0)^\kappa=\mu_0$ and
$(2^{\mu_0})^+\leq \mu$. For a closed unbounded set $C\subseteq\mu$, a
regressive function $f:\mu\longrightarrow\mu$ and a two-place function $g:
\mu^2\longrightarrow\mu_0$ we let:   
\[\hspace{-0.5cm}\begin{array}{rl}
A^*_{C,f,g}=\{\{\alpha_0,\alpha_1,\alpha_2,\alpha_3\}\in [\mu]^4: &
\alpha_0<\alpha_1<\alpha_2<\alpha_3 \mbox{ are from }C,\\
\ & \mbox{each has cofinality}>\kappa,\\
f(\alpha_0 )= & f(\alpha_1)=f(\alpha_2)=f(\alpha_3)\ \mbox{ and}\\
g(\alpha_0,\alpha_2)= & g(\alpha_0,\alpha_3)=g(\alpha_1,\alpha_2)=
g(\alpha_1,\alpha_3)\}.
\end{array}\]
Let $D^4_{\mu,\kappa}$ be the $\kappa$-complete filter on $[\mu]^4$
generated by all the sets $A^*_{C,f,g}$. Finally, let
\[\tau_4=\tau_4(x_0,x_1,x_2,x_3)\stackrel{\rm def}{=} x_0\wedge(-x_1 )\wedge
x_2\wedge(-x_3).\]   
\end{description}
Then the filter $D^4_{\mu,\kappa}$ is proper, $\kappa^+$-complete and normal
for $\langle\mu,\mu,\mu,\mu\rangle$ and every interval Boolean algebra has the
$(D^4_{\mu,\kappa},\tau_4)$-dependence property.  
\end{theorem}
PROOF: {\bf (1)}\ \ \ Let $\mu$ be a regular cardinal such that $(\forall\chi<
\mu)(\chi^\kappa<\mu)$ (so $\mu^\kappa=\mu$). First note that all the sets
$A_{C,f,g}$ are nonempty. [Why? Let $f:\mu\longrightarrow\mu$ be regressive,
$g:\mu^2\longrightarrow\chi$, $\chi<\mu$ and let $C\subseteq\mu$ be a club. 
Then for some $\rho$ the set  
\[S=\{\alpha\in C:\cf(\alpha)>\kappa\ \&\ f(\alpha)=\rho\}\]
is stationary (by Fodor lemma). Next for each $\alpha\in S$ take $h(\alpha)<
\chi$ such that the set $\{\alpha'\in S: \alpha<\alpha'\ \&\ g(\alpha,\alpha')
=h(\alpha)\}$ is stationary, and note that for some $\delta<\chi$ the set 
$Z=\{\alpha\in S: h(\alpha)=\delta\}$ is stationary. Take any $\alpha_0\in Z$
and then choose $\alpha_1<\alpha_2$ from $(\alpha_0,\mu)\cap S$ such that 
\[g(\alpha_0,\alpha_1)=g(\alpha_0,\alpha_2)=\delta.\]
Next choose $\alpha_3>\alpha_2$ from $Z$ and $\alpha_4,\alpha_5\in (\alpha_3,
\mu)\cap S$ such that 
\[g(\alpha_3,\alpha_4)=g(\alpha_3,\alpha_5)=\delta.\]
Clearly $\{\alpha_0,\alpha_1,\alpha_2,\alpha_3,\alpha_4,\alpha_5\}\in
A_{C,f,g}$.] 

Now suppose that $C_\zeta\subseteq\mu$, $f_\zeta:\mu\longrightarrow\mu$,
$g_\zeta:\mu^2\longrightarrow\chi_\zeta$, $\chi_\zeta<\mu$ (for
$\zeta<\kappa$) are as required in the definition of sets $A_{C_\zeta,f_\zeta,
g_\zeta}$. Let $\pi:{}^\kappa\mu\longrightarrow\mu$ be a bijection (remember
$\mu=\mu^\kappa$). Choose a club $C\subseteq\mu$ such that
$C\subseteq\bigcap\limits_{\zeta<\kappa}C_\zeta$ and 
\[\mbox{if } \alpha\in C, \beta<\alpha, F\in {}^\kappa\beta\ \mbox{ then
}\ \pi(F)<\alpha.\]
Let 
\[f:\mu\longrightarrow\mu:\alpha\mapsto\pi(\langle f_\zeta(\alpha):
\zeta<\kappa\rangle),\]
\[g:\mu^2\longrightarrow \prod_{\zeta<\kappa}\chi_\zeta:(\alpha,\beta)\mapsto
\langle g_\zeta(\alpha,\beta):\zeta<\kappa\rangle.\]
The function $f$ is regressive on $\{\alpha\in C:\cf(\alpha)>\kappa\}$, outside
this set we change the values of $f$ to 0. Since $\prod\limits_{\zeta<\kappa}
\chi_\zeta<\mu$ we have $A_{C,f,g}\in D^6_{\mu,\kappa}$. It should be clear
that $A_{C,f,g}\subseteq\bigcap\limits_{\zeta<\kappa}A_{C_\zeta,f_\zeta,
g_\zeta}$. Thus we have proved that the filter $D^6_{\mu,\kappa}$ generated by
the sets $A_{C,f,g}$ is proper $\kappa^+$--complete. To show that
$D^6_{\mu,\kappa}$ is normal for $\langle\mu,\mu,\mu,\mu,\mu,\mu\rangle$ note
that for $\alpha<\mu$, $\ell<6$, if we take $C=(\alpha,\mu)$, $f,g$ constant
functions then  
\[A_{C,f,g}\cap\{\{\alpha_0,\ldots,\alpha_5\}\in [\mu]^6: \alpha_\ell<\alpha\}
=\emptyset.\]

Suppose now that $(I,<_I)$ is a linear ordering. Let $-\infty$ be a new
element (declared to be smaller than all members of $I$) in the case that $I$
has no minimum element; otherwise $-\infty$ is that minimum element. Further,
let $\infty$ be a new element above all members of $I$. The interval Boolean
algebra $B(I)$ determined by the linear ordering $I$ is the algebra of subsets
of $I$ generated by intervals $[x,y)_I=\{z\in I: x\leq_I z<_I y\}$ for $x,y\in
I\cup\{-\infty,\infty\}$.  
\medskip

We are going to show that the algebra $B(I)$ has the $(D^6_{\kappa,\mu},
\tau_6)$--dependence property. Assume that $\langle a_\alpha:\alpha<\mu\rangle
\subseteq B(I)$. Since we can find a subset of $I$ of the size $\leq\mu$ which
captures all the dependences in the sequence we may assume that the linear
order $I$ is of the size $\mu$, so $I$ is a linear ordering on $\mu$. 

Fix a bijection $\phi:[\mu\cup\{-\infty,\infty\}]^{<\omega}\times
{}^{\omega>}4\longrightarrow\mu$.

For each $\alpha<\mu$ we have a (unique) $<_I$-increasing sequence 
\[\langle s^\alpha_i: i<2n(\alpha)\rangle\subseteq\mu\cup\{-\infty,\infty\},
\quad\quad n(\alpha)<\omega\]
such that $a_\alpha=\bigcup\limits_{i<n(\alpha)}[s^\alpha_{2i},
s^\alpha_{2i+1})_I$. Take a closed unbounded set $C\subseteq\mu$ such that for
each $\alpha\in C$:
\begin{description}
\item[(1)] if $w\in[\alpha\cup\{-\infty,\infty\}]^{<\omega}$, $c\in
{}^{\omega>}4$ then $\phi(w,c)<\alpha$,
\item[(2)] if $\phi(w,c)<\alpha$ then $w\subseteq\alpha\cup\{-\infty,
\infty\}$, 
\item[(3)] if $\beta<\alpha$ then $\{s^\beta_i:i<2n(\beta)\}\subseteq\alpha
\cup\{-\infty,\infty\}$.
\end{description}
For each $\alpha<\mu$ fix a finite set $w_\alpha\subseteq\alpha\cup\{-\infty,
\infty\}$ such that $-\infty,\infty\in w_\alpha$ and 
\begin{description}
\item[(4)] if $s^\alpha_i\in\alpha\cup\{-\infty,\infty\}$ then $s^\alpha_i\in
w_\alpha$ and 
\item[(5)] if $s,t\in\{s^\alpha_i: i<2n(\alpha)\}\cup\{-\infty,\infty\}$, $s<_I
t$ and $(s,t)_I\cap\alpha\neq\emptyset$ then $(s,t)_I\cap w_\alpha\neq
\emptyset$. 
\end{description}
Next let $c_\alpha:w_\alpha\longrightarrow 4$ (for $\alpha<\mu$) be such that
for $s\in w_\alpha$: 
\[c_\alpha(s)=\left\{
\begin{array}{ll}
0 & \mbox{if }(\exists x<_I s)([x,s)_I\subseteq a_\alpha),\\
1 & \mbox{if }(\exists x: s<_I x)([s,x)_I\subseteq a_\alpha),\\
2 & \mbox{if both of the above},\\
3 & \mbox{otherwise.}
\end{array}
\right.\]
We can think of $c_\alpha$ as a member of ${}^{\omega>}4$ and we put
$f(\alpha)=\phi(w_\alpha,c_\alpha)$ for $\alpha<\mu$. Note that the function
$f$ is regressive on $C$ (so we can modify it outside $C$ to get a really
regressive function). Now, if $\alpha_0<\alpha_1$, both in $C$,
$f(\alpha_0)=f(\alpha_1)$ then 
\[s^{\alpha_0}_i=s^{\alpha_1}_j\ \&\ i<2n(\alpha_0)\ \&\ j<2n(\alpha_1)\quad
\Rightarrow\quad s^{\alpha_1}_j\in w_{\alpha_1},\quad\mbox{ and}\]
\[w_{\alpha_0}\cap\{s^{\alpha_0}_i:i<2n(\alpha_0)\}= w_{\alpha_1}\cap
\{s^{\alpha_1}_i: i<2n(\alpha_1)\}.\]
[Why? For the first statement note that, by (3), $s^{\alpha_0}_i<\alpha_1$ (for
each $i<2n(\alpha_0)$) so we may use (4). For the second assertion suppose
that $s^{\alpha_0}_{2i}\in w_{\alpha_0}=w_{\alpha_1}$. Then necessarily
$c_{\alpha_0}(s^{\alpha_0}_{2i})=1=c_{\alpha_1}(s^{\alpha_0}_{2i})$. Checking
when the function $c_{\alpha_1}$ takes value $1$ and when $2$ we get that
$s^{\alpha_0}_{2i}=s^{\alpha_1}_{2j}$ for some $j<n(\alpha_1)$. Next, if
$s^{\alpha_0}_{2i+1}\in w_{\alpha_0}=w_{\alpha_1}$ then $c(s^{\alpha_0}_{2i+
1})=0$ and $s^{\alpha_0}_{2i+1}=s^{\alpha_j}_{2j+1}$ for some $j$. Similarly
if we start with $s^{\alpha_1}_i$.] Moreover, if $s,t\in w_{\alpha_0}$ are two
$<_I$-successive points of $w_{\alpha_0}$, $s\leq_I s^{\alpha_1}_i<_I
s^{\alpha_1}_{i+1}\leq_I t$, $i+1<2n(\alpha_1)$, then $(s^{\alpha_1}_i,
s^{\alpha_1}_{i+1})_I\cap\{s^{\alpha_0}_j:j<2n(\alpha_0)\}=\emptyset$. 

Let a function $g:\mu^2\longrightarrow {}^{\omega>}\omega$ be such that if
$\alpha<\beta$, $\alpha,\beta\in C$, $f(\alpha)=f(\beta)$ then
\[g(\alpha,\beta)=\langle \|w_\alpha\|, t^0,\ldots,t^{\|w_\alpha\|-1},
v^0,\ldots, v^{\|w_\alpha\|-1}\rangle\in {}^{\omega>}\omega,\]
where $\bar{t},\bar{v}$ are such that:

if $w_\alpha=\{w_\alpha(0),\ldots,w_\alpha(\|w_\alpha\|-1)\}$ (the
$<_I$-increasing enumeration), $\ell<\|w_\alpha\|$ then
\[t^\ell=0\iff\{s^\beta_j:j<2n(\beta)\}\cap(w_\alpha(\ell),w_\alpha(\ell+1))_I
=\emptyset,\]
and if $s^\beta_i\in (w_\alpha(\ell),w_\alpha(\ell+1))_I$ then 
\[v^\ell=0\ \Rightarrow\ (w_\alpha(\ell),s^\beta_i)_I\cap \{s^\alpha_j:
j<2n(\alpha)\}=\emptyset,\]
\[v^\ell>0\ \Rightarrow\ s^\alpha_{v^\ell-1}\in (w_\alpha(\ell),s^\beta_i)_I\
\&\ (s^\alpha_{v^\ell-1},s^\beta_i)_I\cap\{s^\alpha_j: j<2n(\alpha)\}=
\emptyset.\]
Suppose now that $\alpha_0<\ldots<\alpha_5$ from $C$ are such that $f(\alpha_0)
=\ldots=f(\alpha_5)$, $g(\alpha_0,\alpha_1)=g(\alpha_0,\alpha_2)=g(\alpha_3,
\alpha_4)=g(\alpha_3,\alpha_5)= \langle
k, t^0,\ldots,t^{k-1},v^0,\ldots,v^{k-1}\rangle$. Then $w_{\alpha_0}=\ldots
w_{\alpha_5}=w=\{w(0),\ldots,w(k-1)\}$ (the $<_I$-increasing enumeration).
We are going to show that for each $\ell<k-1$
\begin{description}
\item[($\circledcirc$)] \qquad $\tau_6(a_{\alpha_0},\ldots,a_{\alpha_5})\wedge
[w(\ell),w(\ell+1))_I=\emptyset$. 
\end{description}
Fix $\ell<k-1$. If $t^\ell=0$ then the interval $(w(\ell),w(\ell+1))_I$
contains no $s^{\alpha_1}_j,s^{\alpha_2}_j$ and therefore 
\[a_{\alpha_1}\wedge [w(\ell),w(\ell+1))_I=a_{\alpha_2}\wedge [w(\ell),w(\ell
+1))_I\in\{0,[w(\ell),w(\ell+1))_I\}\]
(remember $c_{\alpha_1}=c_{\alpha_2}$). Hence $(-a_{\alpha_1})\wedge
a_{\alpha_2}\wedge [w(\ell),w(\ell+1))_I=0$ and $(\circledcirc)$ holds. So
suppose that $t^\ell>0$. Then for each $k=1,2,4,5$ the interval
$(w(\ell),w(\ell+1)_I$ contains some $s^{\alpha_k}_j$. We know that if $j<j'$,
$k=1,2$, $s^{\alpha_k}_j, s^{\alpha_k}_{j'}\in (w(\ell),w(\ell+1))_I$ then
there is no $s^{\alpha_0}_i$ in $[s^{\alpha_k}_j,s^{\alpha_k}_{j'}]_I$ 
(and similarly for $\alpha_3$ and $k=4,5$. Assume that $v^\ell=0$ and for
$k=1,2,4,5$ let $j_k<2n(\alpha_k)$ be the last such that $s^{\alpha_k}_{j_k}
\in (w(\ell),w(\ell+1))_I$. By the definition of the functions $g$ and $f$ and
the statement before we conclude that
\begin{quotation}
\noindent either $a_0\wedge [w(\ell), s^{\alpha_k}_{j_k})_I=0$ (for $k=1,2$)
and $a_3\wedge [w(\ell), s^{\alpha_k}_{j_k})_I=0$ (for $k=4,5$)

\noindent or $a_0\wedge [w(\ell), s^{\alpha_k}_{j_k})_I=[w(\ell),s^{\alpha_k
}_{j_k})_I$ (for $k=1,2$) and\\
$a_3\wedge [w(\ell), s^{\alpha_k}_{j_k})_I=[w(\ell),s^{\alpha_k}_{j_k})_I$
(for $k=4,5$)  
\end{quotation}
and the parity of $j_k$'s is the same (just look at $c_{\alpha_k}(w(\ell+
1))$). Hence we conclude that either $a_{\alpha_0}\wedge (-a_{\alpha_1})\wedge
a_{\alpha_2}\wedge [w(\ell),w(\ell+1))_I=0$ or $(-a_{\alpha_3})\wedge
a_{\alpha_4}\wedge (-a_{\alpha_5})\wedge [w(\ell),w(\ell+1))_I=0$ (and in both
cases we get $(\circledcirc)$). Assume now that $v^\ell>0$. By similar
considerations one shows that if $v^\ell-1$ is even then  
\[(-a_{\alpha_3})\wedge a_{\alpha_4}\wedge (-a_{\alpha_5})\wedge [w(\ell),
w(\ell+1))_I=0\]
and if $v^\ell-1$ is odd then
\[a_{\alpha_0}\wedge(-a_{\alpha_1})\wedge a_{\alpha_2})\wedge [w(\ell),w(\ell
+1))_I=0.\]

Since $g$ can be thought of as a function from $\mu^2$ to $\omega<\mu$ the set
$A_{C,f,g}$ is in $D^6_{\mu,\kappa}$ and we have shown that it witnesses
$(D^6_{\mu,\kappa},\tau_6)$-dependence for the sequence $\langle a_\alpha:
\alpha<\mu\rangle$. 
\medskip

\noindent {\bf 2)} It is almost exactly like {\bf 1)} above. The only
difference is that showing that the sets $A^*_{C,f,g}$ are non-empty we use
the Erd\H os--Rado theorem (to choose $\alpha_0,\alpha_1,\alpha_2,\alpha_3$
suitably homogeneous for $g$), and then in arguments that $B(I)$ has the
dependence property we use triples $\alpha_0,\alpha_2,\alpha_3$ and
$\alpha_1,\alpha_2,\alpha_3$. \QED 

\subsection{Appendix: How one can use [Sh 95].}
For reader's convenience we recall here some of the notions and results of
\cite{Sh 95}. We applied them to reduce the number of steps in the beth
hierarchy replacing them partially by passing to successors. This reduction
is meaningful if the exponentiation function is far from GCH. Generally
we think that $\kappa^+$ (or even $\kappa^{++}$) should be considered as
something less than $2^\kappa$.

\begin{definition}
[see Definition 1 of \cite{Sh 95}]
\begin{enumerate}
\item For a sequence $\bar{r}=\langle n_0,\ldots,n_{k-1}\rangle \in{}^k\omega$
we denote:\ \ \  $n(\bar{r})=\sum\limits_{l<k}n_l$, $k(\bar{r})=k$,
$n_l(\bar{r})=n_l$. 
\item Let $B_\xi$ (for $\xi<\mu$) be disjoint well ordered sets, $\bar{r}=
\langle n_0,\ldots,n_{k-1}\rangle\in {}^k\omega$, $f:
[\bigcup\limits_{\xi<\mu} B_\xi]^{n(\bar{r})}\longrightarrow\chi$, $l\leq
n(\bar{r})$. We say that $f$ is $(\bar{r})^l$-canonical (on $\langle B_\xi:
\xi<\mu\rangle$) if 
\begin{quotation}
\noindent for every $\xi_0<\ldots<\xi_{k-1}<\mu$, $a_0<\ldots<a_{n_0-1}$ in
$B_{\xi_0}$, $a_{n_0}<\ldots<a_{n_0+n_1-1}$ in $B_{\xi_1}$ and so on,
the value $f(a_0,\ldots,a_{n(\bar{r})-1})$ depends on $a_0,\ldots,a_{n(
\bar{r})-1-l}$, $\xi_0,\ldots,\xi_{k-1}$ only (i.e.~it does not depend on
$a_{n(\bar{r})-l},\ldots,a_{n(\bar{r})-1}$).  
\end{quotation}
\item A sequence $\langle\lambda_\xi: \xi<\mu\rangle$ (of cardinals) has a
$\langle\kappa_\xi: \xi<\mu\rangle$--canonical form for
$\Gamma=\{(\bar{r}_i)^{l_i}_{\chi_i}: i<\alpha\}$ (where $l_i$'s are
integers, $l_i\leq n(\bar{r}_i)$, $\chi_i$'s are cardinals and $\bar{r}_i$'s
are finite sequences of integers) if 
\begin{quotation}
\noindent for each disjoint (well ordered) sets $A_\xi$,
$\|A_\xi\|=\lambda_\xi$ (for $\xi<\mu$) and functions $f_i:
[\bigcup\limits_{\xi<\kappa} A_\xi]^{n(\bar{r}_i)}\longrightarrow \chi_i$
(for $i<\alpha$) 

\noindent there are sets $B_\xi\subseteq A_\xi$, $\|B_\xi\|=\kappa_\xi$ such
that each function $f_i$ is $(\bar{r}_i)^{l_i}$-canonical on $\langle
B_\xi: \xi<\mu\rangle$ (for $i<\alpha$).
\end{quotation}
\end{enumerate}
\end{definition}

Several canonization theorems were proved in \cite{Sh 95}, we will quote
here two (the simplest actually) which we needed for our applications.

\begin{proposition}
[see Composition Claim 5 of \cite{Sh 95}]
\label{composition}
Let $\Gamma_1$ be 
\[
\begin{array}{rl}
\{(\langle n_0,\ldots,n_{k-1},\ldots,n_{m-1}
\rangle)^{p+q}_{2^\mu}: & (\langle n_0,\ldots,n_{k-1},\ldots,n_{m-1}\rangle
)^p_{2^\mu}\in\Gamma_3\ \&\\
(\langle n_0,\ldots,n_{k-2},n_{k-1}-s\rangle)^q_{2^\mu}\in\Gamma_2
& \&\ p=s+n_k+\ldots+n_{m-1}\ \&\\
\ & \ 0\leq s<n_{k-1}\}.
\end{array}
\] 
Suppose that the sequence $\langle\lambda^3_\xi:\xi<\mu\rangle$ has a
$\langle\lambda^2_\xi:\xi<\mu\rangle$--canonical form for $\Gamma_3$ and the
sequence $\langle\lambda^2_\xi:\xi<\mu\rangle$ has a
$\langle\lambda^1_\xi:\xi<\mu\rangle$--canonical form for $\Gamma_2$.\\
Then the sequence $\langle\lambda^3_\xi:\xi<\mu\rangle$ has a
$\langle\lambda^1_\xi:\xi<\mu\rangle$--canonical form for $\Gamma_1$. \QED
\end{proposition}

\begin{proposition}
[see Conclusion 8(1) of \cite{Sh 95}]
\label{step}
The sequence 
\[\langle (2^\mu)^{++}: \xi<\mu\rangle\]
has a $\langle\mu:\xi<\mu\rangle$--canonical form for $\{(\bar{r}\hat{\ }
\langle 1\rangle)^2_{2^\mu}: \bar{r}\in{}^k\omega, k<\omega\}$. \QED 
\end{proposition}

Recall that for a cardinal $\mu$ and an integer $k$ we have defined
$\daleth_k(\mu)$ by: \ \ \ $\daleth_0(\mu)=\mu$, $\daleth_{k+1}(\mu)=
(2^{\daleth_k(\mu)})^{++}$. 

\begin{proposition}
\label{basic}
Suppose that $\langle A_\xi: \xi<\mu\rangle$ is a sequence of disjoint sets,
$\|A_\xi\|=\daleth_{k+1}(\mu)$. Let $F:[\bigcup\limits_{\xi<\mu}A_\xi]^{
2k+1}\longrightarrow 2^\mu$. Then
\begin{description}
\item[a)\ ] there are $\alpha^0_\xi,\alpha^1_\xi\in A_\xi$ (for $\xi<\mu$),
$\alpha^0_\xi\neq\alpha^1_\xi$ such that for each pairwise distinct
$\xi_0,\ldots,\xi_k<\mu$ 
\end{description}
\begin{description}
\item[$(\oplus)$] $F(\alpha^0_{\xi_0},\alpha^1_{\xi_0},\ldots,\alpha^0_{\xi_{
k-1}}, \alpha^1_{\xi_{k-1}},\alpha^0_{\xi_k})=
F(\alpha^0_{\xi_0},\alpha^1_{\xi_0},\ldots,\alpha^0_{\xi_{k-1}},
\alpha^1_{\xi_{k-1}},\alpha^1_{\xi_k})$
\end{description}
and even more:
\begin{description}
\item[b)\ ] there are sets $B_\xi\in [A_\xi]^\mu$ (for $\xi<\mu$) such that
if $\xi_0,\ldots,\xi_k<\mu$ are distinct, and $a^0_{\xi_i}, a^1_{\xi_i}\in
B_{\xi_i}$ are distinct then $(\oplus)$ of {\bf a)} holds true.
\end{description}
\end{proposition}

\noindent PROOF: It follows from \ref{composition} and \ref{step} (e.g.
inductively) that $\langle \daleth_{k+1}(\mu): \xi<\mu\rangle$ has a
$\langle \mu: \xi<\mu\rangle$--canonical form for $\Gamma$, where $\Gamma$
consists of the following elements:
\[\begin{array}{l}
(\langle\underbrace{2\ldots 2}_k 1\rangle)^2_{2^\mu},
(\langle\underbrace{2\ldots 2}_k 11\rangle)^2_{2^\mu},
(\langle\underbrace{2\ldots 2}_{k-1} 1 2 1\rangle)^4_{2^\mu},\\
(\langle\underbrace{2\ldots 2}_{k-2} 1 2 2 1\rangle)^6_{2^\mu},\ldots,
(\langle1\underbrace{2\ldots 2}_{k} 1\rangle)^{2k+2}_{2^\mu}.
  \end{array}\] 
\QED

\eject
\shlhetal
\end{document}